\documentclass[a4paper,12pt]{amsart}
\usepackage{amssymb}
\usepackage{amsthm}
\usepackage{amscd}
\usepackage{amsmath}
\usepackage[dvipdfm]{graphicx}
\theoremstyle{plain}
\newtheorem{lemma}{Lemma}[subsection]
\newtheorem{prop}[lemma]{Proposition}
\newtheorem{thm}[lemma]{Theorem}
\newtheorem{cor}[lemma]{Corollary}
\newtheorem{aplemma}{Lemma~A.\hspace{-1.5mm}}
\newtheorem{approp}{Proposition~A.\hspace{-1.5mm}}
\newtheorem{apthm}{Theorem~A.\hspace{-1.5mm}}
\newtheorem{apcor}{Corollary~A.\hspace{-1.5mm}}
\newtheorem{intthm}{Theorem}

\usepackage[all]{xy}

\newtheorem{conj}[lemma]{Conjecture}

\theoremstyle{definition}

\newtheorem{rem}{Remark}
\newtheorem{rema}{Remark}

\newtheorem{remb}{Remark}

\newtheorem{defi}[lemma]{Definition}
\newtheorem{exa}[lemma]{Example}
\newtheorem{aprem}{Remark~A.\hspace{-1.5mm}}
\newtheorem{apdefi}{Definition~A.\hspace{-1.5mm}}
\newcommand{\bde}{\begin{defi}}
\newcommand{\ede}{\end{defi}\vspace{1mm}}
\newcommand{\ble}{\begin{lemma}}
\newcommand{\ele}{\end{lemma}}
\newcommand{\bpr}{\begin{prop}}
\newcommand{\epr}{\end{prop}}
\newcommand{\bt}{\begin{thm}}
\newcommand{\et}{\end{thm}}
\newcommand{\bco}{\begin{cor}}
\newcommand{\eco}{\end{cor}}
\newcommand{\bre}{\begin{rem}}
\newcommand{\ere}{\end{rem}}
\newcommand{\brea}{\begin{rema}}
\newcommand{\erea}{\end{rema}\vspace{1mm}}
\newcommand{\breb}{\begin{remb}}
\newcommand{\ereb}{\end{remb}\vspace{1mm}}
\newcommand{\bex}{\begin{exa}}
\newcommand{\eex}{\end{exa}}
\newcommand{\bpf}{\begin{proof}}
\newcommand{\epf}{\end{proof}\vspace{1mm}}

\newcommand{\bade}{\begin{apdefi}}
\newcommand{\eade}{\end{apdefi}}
\newcommand{\bale}{\begin{aplemma}}
\newcommand{\eale}{\end{aplemma}}
\newcommand{\bapr}{\begin{approp}}
\newcommand{\eapr}{\end{approp}}
\newcommand{\bat}{\begin{apthm}}
\newcommand{\eat}{\end{apthm}}
\newcommand{\baco}{\begin{apcor}}
\newcommand{\eaco}{\end{apcor}}
\newcommand{\bare}{\begin{aprem}}
\newcommand{\eare}{\end{aprem}}

\newcommand{\be}{\begin{enumerate}}
\newcommand{\ee}{\end{enumerate}}
\newcommand{\bcd}{\[\begin{CD}}
\newcommand{\ecd}{\end{CD}\]}
\newcommand{\bit}{\begin{itemize}}
\newcommand{\eit}{\end{itemize}}
\newcommand{\bq}{\begin{quote}}
\newcommand{\eq}{\end{quote}}
\newcommand{\ba}{\begin{array}}
\newcommand{\ea}{\end{array}}
\newcommand{\mcA}{\mathcal{A}}

\newcommand{\mcC}{\mathcal{C}}

\newcommand{\mcF}{\mathcal{F}}

\newcommand{\mcI}{\mathcal{I}}
\newcommand{\mcJ}{\mathcal{J}}

\newcommand{\mcL}{\mathcal{L}}

\newcommand{\mcN}{\mathcal{N}}
\newcommand{\mcO}{\mathcal{O}}

\newcommand{\mbA}{\mathbb{A}}

\newcommand{\mbG}{\mathbb{G}}

\newcommand{\mbV}{\mathbb{V}}

\newcommand{\mbZ}{\mathbb{Z}}
\newcommand{\mfA}{\mathfrak{A}}

\newcommand{\mfG}{\mathfrak{G}}

\newcommand{\mfI}{\mathfrak{I}}

\newcommand{\mfS}{\mathfrak{S}}

\newcommand{\mfV}{\mathfrak{V}}

\newcommand{\mfa}{\mathfrak{a}}

\newcommand{\mfc}{\mathfrak{c}}
\newcommand{\mfd}{\mathfrak{d}}
\newcommand{\mfe}{\mathfrak{e}}

\newcommand{\mfh}{\mathfrak{h}}

\newcommand{\mfm}{\mathfrak{m}}

\newcommand{\mfo}{\mathfrak{o}}
\newcommand{\mfp}{\mathfrak{p}}

\newcommand{\mfr}{\mathfrak{r}}
\newcommand{\mfs}{\mathfrak{s}}

\newcommand{\mfz}{\mathfrak{z}}
\newcommand{\migi}{\rightarrow}

\newcommand{\isom}{\stackrel{\sim}{\migi}}

\newcommand{\migiincl}{\hookrightarrow}

\newcommand{\migisurj}{\twoheadrightarrow}

\newcommand{\mr}{\mathrm}
\newcommand{\hidden}[1]{\,}

\begin{document}

\title[Categorical representation of superschemes]{Categorical representation of superschemes}
\author{Yasuhiro Wakabayashi}
\date{}
\markboth{Yasuhiro Wakabayashi}{}
\maketitle
\footnotetext{Y. Wakabayashi: Graduate School of Mathematical Sciences, The University of Tokyo, 3-8-1 Komaba, Meguro, Tokyo,  153-8914, Japan;}
\footnotetext{e-mail: {\tt wkbysh@ms.u-tokyo.ac.jp};}
\begin{abstract}
In the present paper, 
we prove that a locally noetherian superscheme $X^\circledS$ may be reconstructed (up to certain equivalence) category-theoretically from the category  of noetherian superschemes over   $X^\circledS$.
This result  is  a supergeometric generalization of 
the  result proved by Shinichi  Mochizuki concerning categorical reconstruction of  schemes.
% in ~\cite{Mzk1} (cf. ~\cite{Mzk1}, Theorem A) which concerns categorical reconstructibility of locally noetherian schemes to {\it super geometry}.

\end{abstract}
\tableofcontents 
%%%%%%%%%%%%%%%%%%%%%%%%%%%%%%%%--[ begin  section1]---%%%%%%
%\vspace{-3mm}
\section*{Introduction}

\vspace{3mm}
Superschemes (or, supermanifolds) were introduced and discussed in various works  from different point of views, especially in connection with the important physical applications, which stem from superstring theory.
Beside having such physical applications, the  theory of superschemes
will be  interesting on its own from purely mathematical viewpoint.
%A superscheme is an object in supergeometry that generalizes the notion of an ordinary scheme
In the present paper, we are interested in understanding the 
richness of algebraic supergeometry from  category-theoretic aspects.

As a main result of our study, 
%In the present paper,
 we shall give a supergeometric generalization of  the result  proved by S. Mochizuki  (cf. ~\cite{Mzk1}, Theorem A) concerning  categorical reconstructibility of locally noetherian schemes, as described below.
% the problem of reconstructing categorically , which is an analogue of the result in ? for super-schemes.
Let $X^\circledS = (X_b, \mcO_{X^\circledS})$ 
be a superscheme (cf. Definition ~\ref{d3} (i)), i.e., a  scheme $X_b$ together with a certain quasi-coherent 
sheaf of superalgebras $\mcO_{X^\circledS}$ on $X_b$.
%$\mcO_{X_b}$-superalgebra $\mcO_{X^\circledS}$.
% equipped with a structure of sheaves of superrings.
Suppose that $X^\circledS$ is
 locally noetherian  in the sense of   Definition \ref{d3LN}.
% whose underlying scheme $X_b$ is {\it locally noetherian}.
For each such $X^\circledS$, one  obtains the category
\begin{align} \label{DD067}
\mfS \mfc \mfh_{/X^\circledS}^\circledS
\end{align}
consisting  of noetherian superschemes over $X^\circledS$ (cf. (\ref{DD051}) for the precise definition of $\mfS \mfc \mfh_{/X^\circledS}^\circledS$).
%In this manner, we may associate, to  each  $X$, a category $\mr{SSch}_{X}$.
%Note that the  fiber products in $\mr{SSch}_X$ exist (cf. ?, Corollary 10.3.9). 
%To simplify terminology, 
%we simply  write $Y$ for 
%we shall often refer to the domain $Y$ of an arrow $Y \migi X$ which is an object of $\mr{SSch}_X$ as an ``{\it object of $\mr{SSch}_X$}".
 %  of super-schemes over $X^s$ (cf. the discussion following Definition ?).
The problem that we  consider in the present paper  
is to know  to what extent one can reconstruct the  superscheme-theoretic structure of   $X$  from the categorical structure of $\mfS \mfc \mfh_{/X^\circledS}^\circledS$.
Our  main result is the following assertion.

%We first verify  the following proposition. Here, for any noetherian scheme $X'$, we shall denote by $\mr{Sch}_{X'}$ the category whose objects are morphism of finite type $Y' \migi X'$, where $Y'$ is a noetherian scheme and whose morphisms (from an object $Y'_1 \migi X'$ to an object $Y'_2 \migi X$) are morphisms of finite type $Y'_1 \migi Y'_2$ lying over $X'$.

%---------------------------------------------------------------------[begin theorem]----------------------
%\vspace{3mm}
 \vspace{5mm}
\begin{intthm}
% [{\bf Generic \'{e}taleness of $\mfO \mfp_{\mfs \mfl_n, \hslash, \hslash \star \rho, g,r}^\ZZZ/\overline{\mfM}_{g,r}$}]
  \label{y0194A}
\leavevmode\\
\ \ \ 
%\vspace{-5mm}
Let
$X^\circledS$ and $X'^\circledS$ be two locally noetherian  superschemes.
% both of whose underlying schemes.
% are locally noetherian.
Then, 
\begin{equation}
X^\circledS \stackrel{f}{\sim} X'^\circledS  \ \ \text{if and only if}  \ \ \mfS \mfc \mfh^\circledS_{/X^\circledS} \cong \mfS \mfc \mfh^\circledS_{/X'^\circledS}.
\end{equation}
(Here,  $\stackrel{f}{\sim}$ denotes the equivalence relation defined in (\ref{EE07}).)
\end{intthm}
\vspace{5mm}
%------------------------------------------------------------------------------[end theorem]-------------

%The difference from the assertion in ~\cite{Mzk1} is
Theorem A implies, unlike the result of ~\cite{Mzk1}, 
 that isomorphism classes of  locally noetherian superschemes may not be  determined uniquely  from  the  categorical structure of   $\mfS \mfc \mfh^\circledS_{/X^\circledS}$.
 %, unlike the result of ~\cite{Mzk1}, Theorem A.
Indeed, 
suppose that $X^\circledS \stackrel{f}{\sim} X'^\circledS$, that is to say, $X'^\circledS$ is isomorphic to a fermionic twist of $X^\circledS$ (cf. Definition \ref{d33f4}).
By definition, $X'^\circledS$ may be obtained by twisting the fermionic portion of $X^\circledS$ by means of some element $a$ in  the first \'{e}tale cohomology group $H^1_{\text{\'{e}t}} (X_b, \mu_2)$.
By twisting various superschemes over $X^\circledS$ by means of $a$ in the same manner, we obtain the assignment from each object in $\mfS \mfc \mfh^\circledS_{/X^\circledS}$ to an object in $\mfS \mfc \mfh^\circledS_{/X'^\circledS}$; this assignment gives an equivalence of categories $\mfS \mfc \mfh^\circledS_{/X^\circledS} \isom \mfS \mfc \mfh^\circledS_{/X'^\circledS}$, and hence, shows one direction of the equivalence in Theorem A  (cf. Proposition \ref{p23} and the discussion in its proof).

On the other hand, the proof of the reverse direction (i.e., $\mfS \mfc \mfh^\circledS_{/X^\circledS} \cong \mfS \mfc \mfh^\circledS_{/X'^\circledS}$  implies 
%\Rightarrow 
$X^\circledS \stackrel{f}{\sim} X'^\circledS$) is technically much more difficult.
To complete the proof,  we reconstruct step-by-step various partial information of (the equivalence class of) $X^\circledS$ from the categorical structure of $\mfS \mfc \mfh_{/X^\circledS}^\circledS$, as discussed in \S\,2.
If $X^\circledS$ is a scheme in the classical sense, then any fermionic twists of $X'^\circledS$ is in fact isomorphic to $X^\circledS$ (in particular, $X'^\circledS$ is a scheme); in this case, Theorem A is exactly the same  as the result by S. Mochizuki.
%according to the definition of the equivalence relation ``\ $\stackrel{{\it fmt}}{\sim}$ \ ",  the super-scheme $Y$ obtained from $X$ by deforming,  in a certain manner, its fermionic structure  
%yields the category $\mr{SSch}_Y$ equivalent to $\mr{SSch}_X$.
%But, if $X$ satisfies  certain conditions, then
%one may prove more strong assertion which states that 
 %$X$ may be  completely determined from the category $\mr{SSch}_X$ as we shall assert below. 
%Let $X$ and $Y$ be as in Theorem A.
% and $f : Y \migi X$ a morphism of super-schemes (cf. Definition \ref{d3}, (ii))
%Denote  by $\mr{Isom} (Y, X)$ the set of isomorphisms of super-schemes from $Y$ to $X$.
%Also, denote by $\mr{Isom}( \mr{SSch}_{X}, \mr{SSch}_{Y})$ the set of isomorphism classes of (objects in the category of) equivalences $\mr{SSch}_{X} \isom \mr{SSch}_{Y}$.
%Since  the  fiber products in $\mr{SSch}_X$ exist (cf. ?, Corollary 10.3.9),
%one may define, for each morphism $f : Y \migi X$ of super-schemes,
%a functor
%\begin{equation}
%\mr{SSch}_f : \mr{SSch}_X \migi \mr{SSch}_Y
%\end{equation}
%induced by base-change by $f$. 
%(i.e., the set of isomorphism classes of objects in the category of e)

In the last section of the present paper, we shall prove further rigidity properties  concerning the category $\mfS \mfc \mfh_{/X^\circledS}^\circledS$ (cf. Propositions \ref{p756} and \ref{p4FRP}).

Finally, we want to remark that,  as a different  type of  reconstruction of  a superscheme,  one may find the  result in  ~\cite{D}, which asserts that  a superscheme may be reconstructed from the $(\mbZ/2 \mbZ)$-graded tensor triangulated category of perfect complexes on it (cf. Remark \ref{rg0df} of the present paper).

\vspace{6mm}
\hspace{-4mm}{\bf Acknowledgement} \leavevmode\\
 \ \ \ The author would like to thank
 %express his sincere gratitude to
  Professors Shinichi Mochizuki  for 
 his inspiring works concerning categorical reconstructibility of  geometric objects.
The author cannot express enough his sincere and  deep gratitude to  all those who give  the opportunity or   impart  the great joy of  studying 
%getting in touch with 
  mathematics to him.
The author means the present paper for a gratitude letter to them.
The author was partially  supported by the Grant-in-Aid for Scientific Research (KAKENHI No.\,15J02721) and the FMSP program at the Graduate School of Mathematical Sciences of the University of Tokyo.
% their  many helpful suggestions and advices.
%Without their philosophies and amazing insights, my study for mathematics would remain ``dormant".
%\vspace{10mm}

%%%%%%%%%%%%%%%%%%%%%%%%%%%%%%%%--[ begin  section1]---%%%%%%
\vspace{10mm}
\section{Superschemes} \vspace{3mm}

%----------------------------------------------------------------------[begin subsection]-------------
%\subsection{}
In this section, we  recall first the definition of a superscheme defined over $\mbZ [\frac{1}{2}]$ (cf. Definition \ref{d3}).
Then, we  introduce the notion  of a {\it fermionic twist} (cf. Definition \ref{d33ff4}), and  the  equivalence relation $\stackrel{f}{\sim}$ (cf. (\ref{EE07}))  appeared in the statement of Theorem A. 
%f$-deformation (cf. Definition ?), $f$-twisted deformation (cf. Definition ?), and $f$-equivalence (cf. Definition ?).
%The motivation for introducing these concepts is to deal with 
One direction of the equivalence in Theorem A (which is much  easier to prove than the reverse direction) will be proved in \S\,\ref{S13} (cf. Proposition \ref{p23u}).
% appeared in Theorem A.
%Since our primary interest throughout the present paper will be ``supergeomatric" objects,
%the objects and constructions of ordinary algebraic geometry will be distinguished by placing a ``bosonic" before them;  thus, as in ~\cite{AG},  a classical scheme will be referred to as a ``{\bf bosonic scheme}" for the remainder of the present paper.

Throughout the present paper, we denote, for any category  $\mcC$,
by $\mr{Ob} (\mcC)$ the set of objects of $\mcC$.
 Also, if both $A$ and $B$ are  objects of $\mcC$ (i.e., $A$, $B \in \mr{Ob}(\mcC)$), then we shall  denote by   $\mr{Map}_\mcC (A, B)$  the set of morphisms  (in $\mcC$) from $A$ to $B$.

\vspace{5mm}
%----------------------------------------------------------------------[begin subsection]-------------
\subsection{Superschemes} \label{sub1}
\vspace{-4mm}

%-----------------------------------------------------------------------[begin definition]------------------
\vspace{3mm}
\bde \label{d3}\leavevmode\\
 \vspace{-5mm}
\begin{itemize}
\item[(i)]
 A {\bf superscheme} is a pair  $X^\circledS := (X_b, \mcO_X)$ consisting of a scheme $X_b$ over $\mbZ [\frac{1}{2}]$ and a quasi-coherent sheaf of superalgebras $\mcO_{X^\circledS}$  over $\mcO_{X_b}$ 
 such that the natural morphism $\mcO_{X_b} \migi \mcO_{X^\circledS}$ is injective and its image coincides with 
 the bosonic (i.e., even) part of $\mcO_{X^\circledS}$.
 We shall write $\mcO_{X_f}$ for the fermionic (i.e., odd) part of $\mcO_{X^\circledS}$ and identify $\mcO_{X_b}$ with the bosonic part via the injection $\mcO_{X_b} \migi \mcO_{X^\circledS}$ (hence, $\mcO_{X^\circledS} = \mcO_{X_b} \oplus \mcO_{X_f}$).
%\item[(ii)]
%We shall say that a super-scheme $X := (X_b, \mcO_X)$ is {\bf of commutative type}
%if $\phi \cdot \phi' = 0$ for any local sections $\phi$, $\phi' \in  \mcO_{X_f}$.
 \item[(ii)]
 Let $X^\circledS := (X_b, \mcO_{X^\circledS})$ and $Y^\circledS := (Y_b, \mcO_{Y^\circledS})$ be two superschemes.
 A {\bf morphism of superschemes} from $Y^\circledS$  to $X^\circledS$  is a pair $f^\circledS := (f_b, f^\flat)$ consisting of a  morphism $f_b : Y_b \migi X_b$ of schemes and a morphism of superalgebras $f^\flat : f_b^*(\mcO_{X^\circledS}) \ (:= \mcO_{Y_b} \otimes_{f^{-1}_b (\mcO_{X_b})} f_b^{-1} (\mcO_{X^\circledS})) \migi \mcO_{Y^\circledS}$ over $\mcO_{Y_b}$. 
 \end{itemize}
 \ede
%\vspace{3mm}
%-----------------------------------------------------------------------[end definition]-------------------

In the following, let us fix a superscheme $X^\circledS  := (X_b, \mcO_{X^\circledS})$.

%-----------------------------------------------------------------------[begin definition]------------------
\vspace{3mm}
\bde \label{d3LN}\leavevmode\\
% \vspace{-5mm}
%\begin{itemize}
%\item[(i)]
 \ \ \ We shall say that 
%a super-scheme $X := (X_b, \mcO_X)$
$X^\circledS$ is {\bf locally noetherian} (resp., {\bf noetherian})
if $X_b$ is locally noetherian (resp., noetherian) and the $\mcO_{X_b}$-module $\mcO_{X_f}$ is coherent.
%\item[(ii)]
%We shall say that a super-scheme $X := (X_b, \mcO_X)$ is {\bf of commutative type}
%if $\phi \cdot \phi' = 0$ for any local sections $\phi$, $\phi' \in  \mcO_{X_f}$.
% \item[(ii)]
% \end{itemize}
 \ede
%\vspace{3mm}
%-----------------------------------------------------------------------[end definition]-------------------

We shall  denote by
\begin{equation} \label{DD051}
\mfS \mfc \mfh_{/X^\circledS}^\circledS
\end{equation}
 the category defined as follows:
 \begin{itemize}
 \item[$\bullet$]
 the  {\it objects} are morphisms of superschemes $Y^\circledS \ (= (Y_b, \mcO_{Y^\circledS})) \migi X^\circledS$ to $X^\circledS$ such that $Y^\circledS$ is noetherian  and the underlying morphism $Y_b \migi X_b$ of schemes is {\it of finite type};
\item[$\bullet$]
the {\it morphisms} (from an object $Y^\circledS_1 \migi X^\circledS$ to an object $Y^\circledS_2 \migi X^\circledS$) are morphisms of superschemes $Y^\circledS_1 \migi Y_2^\circledS$ lying over $X^\circledS$.
\end{itemize} 
% , and whose 
The fiber products and finite coproducts exist in $\mfS \mfc \mfh_{/X^\circledS}^\circledS$
%(which is well-defined since  the  fiber products in $\mr{SSch}_Y$ exist
 (cf. ~\cite{CR}, Corollary 10.3.9).

%-----------------------------------------------------------------------[begin remark]------------------
\begin{rema} \label{r4f} \leavevmode\\
 \ \ \ 
Let $X$ be a scheme (in the usual sense) over over $\mbZ [\frac{1}{2}]$.
Then, $X$ carries 
%In an evident fashion,  the  bosonic schemes  over $\mbZ [\frac{1}{2}]$ may be 
 a superschemes of the form $X^\circledS_{\mr{triv}} := (X, \mcO_{X^\circledS_{\mr{triv}}} \ (= \mcO_{X}\oplus \mcO_{X_f}))$ with $\mcO_{X_f} = 0$.
 (Conversely, any superscheme with vanishing fermionic part arises uniquely  from a scheme in this manner.)
In the rest of the present paper, we shall not distinguish between $X$ and $X^\circledS_{\mr{triv}}$.

 % the category whose {\it objects} are morphisms $Y \migi X^\circledS$ (hence, $Y$ lies  over $\mbZ [\frac{1}{2}]$)
% such that  $Y$ is a noetherian  scheme and the underlying morphism $Y \migi X_b$ of schemes is of finite type, and whose {\it morphisms} (from an object $Y_1 \migi X^\circledS$ to an object $Y_2 \migi X^\circledS$) are morphisms of finite type $Y_1 \migi Y_2$ lying over $X^\circledS$.
%The assignment $Y \mapsto Y_{\mr{triv}}^\circledS$ defines
% a fully faithful functor $\mfS \mfc \mfh_{/X^\circledS} \migi \mfS \mfc \mfh^\circledS_{/X^\circledS}$,
% which preserves the fiber product and finite coproducts.  
% By means of this functor,  we  identify $\mfS \mfc \mfh_{/X^\circledS}$ with a full subcategory of $\mfS \mfc \mfh^\circledS_{/X^\circledS}$.
%As in the definition of $\mr{SSch}_X$, we shall denote, for each  super-scheme $X$, by 
%\begin{equation}
%\mr{BSch}_X
%\end{equation}

  \end{rema}
%-----------------------------------------------------------------------[end remark]-------------------
%\vspace{5mm}

\vspace{5mm}
%----------------------------------------------------------------------[begin subsection]-------------
\subsection{Superschemes arising from a bilinear map} \label{ssub21}
\leavevmode\\
\vspace{-4mm}

%-----------------------------------------------------------------------[begin remark]------------------
%\begin{rema} \label{r4g} \leavevmode\\
 \ \ \ 
 Let $X^\circledS := (X_b, \mcO_{X^\circledS})$ be a superscheme.
 The multiplication morphism  $\mcO_{X^\circledS} \otimes \mcO_{X^\circledS} \migi \mcO_{X^\circledS}$ restricts to a  skew-symmetric $\mcO_{X_b}$-bilinear  map
 \begin{align} \label{EE010}
 m_{X^\circledS} :\mcO_{X_f}^{\otimes 2} \ (:= \mcO_{X_f} \otimes_{\mcO_{X_b}} \mcO_{X_f}) \migi \mcO_{X_b}.
 \end{align}
The associative property  of   the  multiplication gives rise to  the equality
 \begin{align} \label{EE066}
 m_{X^\circledS} \otimes  \mr{id}_{\mcO_{X_f}} =  \mr{id}_{\mcO_{X_f}} \otimes m_{X^\circledS} : \mcO_{X_f}^{\otimes 3} \migi \mcO_{X_f}^{}.
 \end{align}
 One verifies that  the  superscheme $X^\circledS$ is uniquely determined (up to isomorphism) by the triple
  \begin{align}
\mcA_{X^\circledS} :=  (X_b, \mcO_{X_f}, m_{X^\circledS}).
  \end{align}
% satisfying the equality (\ref{EE066}).
To make the discussion  precise, let us define 
\begin{align}
\mfA
\end{align}
to be the category, where
\begin{itemize}
\item[$\bullet$]
the {\it objects} are  triples $(Y, \mcF, \omega)$  consisting of a noetherian scheme $Y$ of finite type  over $\mbZ [\frac{1}{2}]$, a coherent  $\mcO_Y$-module $\mcF$, and a  skew-symmetric  $\mcO_{Y}$-bilinear map $\omega : \mcF \otimes \mcF \migi \mcO_{Y}$ on $\mcF$ satisfying the equality $\omega \otimes \mr{id}_\mcF = \mr{id}_\mcF \otimes \omega : \mcF^{\otimes 3} \migi \mcF$;
\item[$\bullet$]
the {\it morphisms} from $(Y, \mcF, \omega)$ to $(Y', \mcF', \omega')$ (where both $(Y, \mcF, \omega)$ and $(Y', \mcF', \omega')$ are objects of $\mfA$) are pairs $(f, f^\flat)$ consisting of a morphism $f : Y \migi Y'$ of schemes and an $\mcO_Y$-linear morphism $f^\flat : f^*(\mcF') \migi \mcF$ satisfying the equality 
\begin{align}
\omega \circ (f^\flat \otimes f^\flat) = f^*(\omega') : f^*(\mcF') \otimes f^* (\mcF') \ (= f^* (\mcF' \otimes \mcF'))\migi \mcO_{Y'}.
\end{align}
\end{itemize}
Then, the following proposition is verified.

%-----------------------------------------------------------------------[begin proposition]------------------
\vspace{3mm}
\bpr \label{p24pv}\leavevmode\\
 \ \ \ 
The assignment $X^\circledS \mapsto \mcA_{X^\circledS}$ defined above  is functorial, and
the resulting functor 
\begin{align} \label{R01}
(\mfS \mfc \mfh_{/\mr{Spec}(\mbZ [\frac{1}{2}])}^\circledS =:) \ \mfS \mfc \mfh_{/\mbZ [\frac{1}{2}]}^\circledS \migi  \mfA
\end{align}
is  an equivalence of categories. 
%In particular, superschemes over supers satisfies descent for \'{e}tale coverings.
\epr
%-----------------------------------------------------------------------[begin proof]-------------------
\begin{proof}
%First, we shall consider the former assertion.
Let us take an object $(Y, \mcF, \omega)$ of $\mfA$.
Then,
the direct sum $\mcO_Y \oplus \mcF$ admits a structure of $\mcO_{Y}$-superalgebra 
   (where the first and second factors are the bosonic and fermionic  parts respectively)
 with multiplication given by 
 \begin{align} \label{}
( \mcO_{Y} \oplus \mcF) \otimes (\mcO_{Y} \oplus \mcF) & \migi \hspace{10mm}\mcO_{Y} \oplus \mcF \\
(a,  \epsilon_a) \otimes (b,  \epsilon_b) \hspace{5mm}& \mapsto (a  b + \omega (\epsilon_a, \epsilon_b), a  \epsilon_b + b  \epsilon_a). \notag
 \end{align}  
The pair $Y^\circledS_{\mcF, \omega} := (Y, \mcO_{Y}\oplus \mcF)$ forms a superscheme and the resulting 
 assignment $(Y, \mcF, \omega) \mapsto Y^\circledS_{\mcF, \omega}$ is functorial in $\mfA$.
 This assignment  defines a functor $\mfA \migi \mfS \mfc \mfh^\circledS_{/\mbZ [\frac{1}{2}]}$ which  is the inverse to the functor (\ref{R01}).
 %$\mfS \mfc \mfh^\circledS_{/\mbZ [\frac{1}{2}]} \migi \mfA$ given by $X^\circledS \mapsto \mfA_{X^\circledS}$.
This completes the proof of \ref{p24pv}.
%The latter assertion follows from the former assertion and a well-known generalities of descent theory.
\end{proof}
%\vspace{3mm}
%-----------------------------------------------------------------------[end proposition]-------------------

%are the  inverse maps of   each other.
%But, every such triple need not to arise from a superscheme.
%  \end{rema}
%-----------------------------------------------------------------------[end remark]-------------------
%\vspace{5mm}

\vspace{5mm}
%----------------------------------------------------------------------[begin subsection]-------------
\subsection{From superschemes to schemes} \label{sub21}
\leavevmode\\
\vspace{-4mm}

In the following, we shall fix a superscheme $X^\circledS : = (X_b, \mcO_{X^\circledS} \ (= \mcO_{X_b}\oplus \mcO_{X_f}))$.
%If $X : = (X_b, \mcO_X)$ is a super-scheme, then 
By considering  the morphism 
\begin{align} \label{e22}
\beta_{X^\circledS} : X^\circledS \migi X_b
\end{align}
 corresponding to  
the inclusion $\mcO_{X_b} \migi \mcO_{X^\circledS}$,
$X^\circledS$ may be thought of as a superscheme over  the scheme  $X_b$.
% yields a morphism 
The construction of $\beta_{X^\circledS}$ is evidently  functorial in $X^\circledS$, that is to say,
$\beta_{X^\circledS} \circ f^\circledS = f_b \circ \beta_{Y^\circledS}$ for any superscheme $Y^\circledS$ and any morphism
 $f^\circledS \ (:= (f_b, f^\flat))  : Y^\circledS \migi X^\circledS$ of superschemes.

Also, 
 denote by
 \begin{align}
 \mcN_{X^\circledS}
 \end{align}
the superideal of $\mcO_{X^\circledS}$ generated by $\mcO_{X_f}$.
%We shall write
%\begin{align} \label{EE090}
%\tau_{X^\circledS} : X_t \migi X^\circledS
%\end{align}
The quotient of $\mcO_{X^\circledS}$ by $\mcN_{X^\circledS}$
% th
%two-sided
%superideal sheaf  generated by $\mcO_{X_f}$
 determines a  scheme $X_t$ equipped with a morphism
\begin{equation} \label{f6}
\tau_{X^\circledS} : X_t \migi X^\circledS
\end{equation} 
of superschemes.
The composite
\begin{align}
\gamma_X := \beta_{X^\circledS} \circ \tau_{X^\circledS}  : X_t \migi X_b
\end{align}
is a closed immersion of schemes  corresponding to the quotient $\mcO_{X_b} \migisurj \mcO_{X_b}/\mcO^2_{X_f}$ ($= \mcO_{X^\circledS}/\mcN_{X^\circledS}$) by the nilpotent ideal $\mcO^2_{X_f} \subseteq \mcO_{X_b}$.

If $f^\circledS: Y^\circledS \migi X^\circledS$ is a morphism of superschemes, then it induces  a  morphism
\begin{equation} \label{e421}
f_t : Y_t \migi X_t
\end{equation}
 of schemes satisfying that  $\tau_{X^\circledS} \circ f_t = f^\circledS \circ \tau_{Y^\circledS}$.

%For each superscheme $X^\circledS$, w
Next, we  denote by 
\begin{align} \label{E0001}
\mfS \mfc \mfh_{/X^\circledS}
\end{align}
 the full subcategory of $\mfS \mfc \mfh_{/X^\circledS}^\circledS$ consisting of objects  of the form $Y \migi X^\circledS$, where $Y$ is a  scheme.
%Let $X^\circledS$ be a locally noetherian superscheme.
%If we are given a scheme $X$, then t
The assignment $Y^\circledS \mapsto Y_t$ ($Y^\circledS \in \mr{Ob} (\mfS \mfc \mfh_{/X^\circledS}^\circledS)$) defines a functor 
\begin{equation} \label{DD014}
\tau : \mfS \mfc \mfh_{/X^\circledS}^\circledS \migi \mfS \mfc \mfh_{/X_t}
\end{equation}
%given by assigning 
which turns out to be  a right adjoint functor of the functor 
  \begin{align} \label{R09}
  \mfS \mfc \mfh_{/X_t} \hspace{3mm}& \migi \hspace{5mm}\mfS \mfc \mfh^\circledS_{/X^\circledS} \\
 ``Z \migi X_t" &  \mapsto ``Z \migi X_t \stackrel{\tau_{X^\circledS}}{\migi} X^\circledS". \notag
  \end{align}
 That is to say, the  functorial map of sets
 \begin{equation} \label{e968}
 \mr{Map}_{\mfS \mfc \mfh_{/X^\circledS}^\circledS} (Z, Y^\circledS) \migi \mr{Map}_{\mfS \mfc \mfh_{/X_t}} (Z, Y_t)
 \end{equation}
is bijective, where $Y \in \mr{Ob} (\mfS \mfc \mfh_{/X_t})$ and $Z^\circledS \in \mr{Ob} (\mfS \mfc \mfh_{/X^\circledS}^\circledS)$.
In particular, we obtain an equivalence of categories $\mfS \mfc \mfh_{/X_t} \isom \mfS \mfc \mfh_{/X^\circledS}$ (given as in (\ref{R09})).

%In the following, we introduce some definitions concerning superschemes.

%-----------------------------------------------------------------------[begin definition]------------------
\vspace{3mm}
\bde \label{d333}\leavevmode\\
 \ \ \ Let $X^\circledS$ be a superscheme and $U \migi X_b$ be an \'{e}tale morphism.  
  %an open subset of the underlying topological space of $X^\circledS$.
Then, we shall write
%the restriction 
\begin{align} \label{E0002}
X^\circledS  |_U  := X^\circledS \times_{\beta_{X^\circledS}, X_b} U
%(X_b |_U, \mcO_{X^\circledS} |_U)
\end{align}
% of $X^\circledS$  forms a superscheme together  with a morphism $X^\circledS |_U \migi X^\circledS$.
By an  {\bf open subsuperscheme} (resp., a  {\bf quasi-compact open subsuperscheme}) of $X^\circledS$, we mean   a superscheme of the form $X^\circledS  |_U$ for some open subscheme (resp., quasi-compact subscheme) $U$ of $X_b$.
%If $U$ is an open subscheme  (a quasi-compact open subscheme), then we shall refer to the  superscheme 
% (associated with some open subset $U$) 
%as an {\bf open subsuperscheme} ().
%, which we  denote by $X |_U$.
 \ede
%\vspace{3mm}
%-----------------------------------------------------------------------[end definition]-------------------

\vspace{5mm}
%----------------------------------------------------------------------[begin subsection]-------------
\subsection{Fermionic twists} \label{S13}
\leavevmode\\
\vspace{-4mm}

Let us  define the notion of a fermionic twist of a given superscheme.
%an equivalence relation  in the set of superschemes.
In the following, let us fix a  locally noetherian superscheme $X^\circledS:= (X_b, \mcO_{X^\circledS})$.
%In particular, the scheme-theoretic support
%\begin{align}
%\overline{X}_b := \mr{Supp} (\mcO_{X_f})
%\end{align}
% of the  $\mcO_{X_b}$-module $\mcO_{X_f}$ forms a closed subscheme of $X_b$.
%We shall write
%\begin{align}
%\overline{X}^\circledS := X^\circledS \times_{X_b} \overline{X}_b \ (:= (\overline{X}_b, \mcO_{\overline{X}_f})),
%\end{align}
%which admits naturally a morphism $\overline{X}^\circledS \migi X^\circledS$.

% Let $Y^\circledS := (Y_b, \mcO_{Y^\circledS})$ be a superscheme.
 % For each super-scheme $Y := (Y_b, \mcO_Y)$, 
% we shall write 
 We shall define $(-1)_{X^\circledS}$ to be  the automorphism
\begin{equation}
(-1)_{X^\circledS} := (\mr{id}_{X_b}, (-1)^\flat_{X^\circledS}) : X^\circledS \isom X^\circledS
\end{equation}
of $X^\circledS$, where $(-1)^\flat_{X^\circledS}$ denotes  the automorphism of $\mcO_{X^\circledS} = \mcO_{X_b} \oplus \mcO_{X_f}$ given by assigning $(a, \epsilon_a) \mapsto (a, - \epsilon_a)$.
In particular, $(-1)_{X^\circledS} \circ (-1)_{X^\circledS} = \mr{id}_{X^\circledS}$, and if $X^\circledS$ is a scheme (i.e., $\mcO_{X_f} =0$), then we have $(-1)_{X^\circledS} = \mr{id}_{X^\circledS}$.
If, moreover, $Y^\circledS$ is a locally noetherian superscheme and $f^\circledS : Y^\circledS \migi X^\circledS$ is a morphism of superschemes, then we have the equality of morphisms
 $f^\circledS \circ (-1)_{Y^\circledS} = (-1)_{X^\circledS} \circ f^\circledS$.
   Hence,  the collection  of automorphisms $\{ (-1)_{Y^\circledS} \}_{Y^\circledS \in \mr{Ob} (\mfS \mfc \mfh^\circledS_{/\mbZ [\frac{1}{2}]})}$ defines a nontrivial center of $\mfS \mfc \mfh^\circledS_{/\mbZ[\frac{1}{2}]}$ (i.e., an automorphism of the identity functor $\mfS \mfc \mfh^\circledS_{/\mbZ[\frac{1}{2}]} \isom \mfS \mfc \mfh^\circledS_{/\mbZ[\frac{1}{2}]}$).

%-----------------------------------------------------------------------[begin definition]------------------
\vspace{3mm}
\bde \label{d33ff4}\leavevmode\\
 \ \ \ 
We shall refer to $(-1)_{X^\circledS}$ as the {\bf fermionic involution} of $X^\circledS$.
%, where $(-1)^\sharp_Y$ denotes the automorphism of $\mcO_Y = \mcO_{Y_b}\oplus \mcO_{Y_f}$ given by
%$(a, b) \mapsto (a, -b)$.
%Note that $(-1)_Y \circ (-1)_Y = \mr{id}_Y$ and for each morphism $f : Y \migi X$, we have the equality $f \circ (-1)_Y = (-1)_{X} \circ f$.
%If $Y$ is bosonic, then $(-1)_Y = \mr{id}_Y$.
\ede
%\vspace{3mm}
%-----------------------------------------------------------------------[end definition]-------------------

Write  $\mr{Aut}_{X_b}(X^\circledS)$ for  the \'{e}tale sheaf on $X_b$ consisting of locally defined automorphisms of $X^\circledS$ over $X_b$ 
(i.e., the sheaf which, to any \'{e}tale scheme $U$ over $X_b$, assigns the group of automorphisms of $X^\circledS |_U$ over $U$),
%(cf. the latter assertion of Proposition ?),
 and $ (\mu_2)_{X_b}$ for the constant \'{e}tale sheaf on $X_b$ with coefficients in the square  roots of unity $\mu_2 := \{ \pm 1\}$.
Then, we have a homomorphism
\begin{align}
\eta_{X^\circledS} : (\mu_2)_{X_b} \migi \mr{Aut}_{X_b}(X^\circledS)
\end{align}
determined by $\eta_{X^\circledS} (1) = \mr{id}_{X^\circledS}$ and $\eta_{X^\circledS} (-1) = (-1)_{X^\circledS}$.
By applying the functor $H^1_{\text{\'{e}t}} (X_b, -)$, we have a homomorphism
\begin{align}
H^1_{\text{\'{e}t}} ( \eta_{X^\circledS}) : H^1_{\text{\'{e}t}} (X_b, \mu_2) \migi H^1_{\text{\'{e}t}} (X_b, \mr{Aut}_{X_b}(X^\circledS))
\end{align}

%-----------------------------------------------------------------------[begin definition]------------------
\vspace{3mm}
\bde \label{d33f4}\leavevmode\\
 \ \ \ A {\bf fermionic twist of $X^\circledS$} is a superscheme 
defined to be the twisted form of $X^\circledS$ (over the \'{e}tale topology on $X_b$) corresponding to  $H^1_{\text{\'{e}t}} ( \eta_{X^\circledS}) (a) \in H^1_{\text{\'{e}t}} (X_b, \mr{Aut}_{X_b}(X^\circledS))$ for some $a \in H^1_{\text{\'{e}t}} (X_b, \mu_2)$. 
We shall refer to this superscheme as the {\bf fermionic twist of $X^\circledS$ associated with $a$} and denote it  by 
\begin{align}
{^a X}^\circledS.
\end{align}
%such that there exist an \'{e}tale covering $U$ of $X_b$ and an isomorphism
%\begin{align}
%\eta : Y^\circledS \times_{X_b} U \isom X^\circledS \times_{X_b} U
%\end{align}
% satisfying the following condition:
%if 
%$\eta'_U$ and $\eta''_U$ denote the isomorphisms 
%\begin{align}
%\eta_1, \eta_2 : Y^\circledS \times_{X_b} U  \times_{X_b} U \isom X^\circledS \times_{X_b} U  \times_{X_b} U
%\end{align}
%denote the isomorphisms
% obtained as the  base-changes of  $\eta$ via the first and second projections $U \times_{X_b} U \migi U$ respectively, then the restriction of  
%$\eta_2 \circ \eta_1^{-1} \in \mr{Aut} (X^\circledS \times_{X_b} U  \times_{X_b} U) $ to  each connected component  coincides with either the identity morphism or the fermionic involution.
  \ede
%\vspace{3mm}
%-----------------------------------------------------------------------[end definition]-------------------
%-----------------------------------------------------------------------[begin remark]------------------
\begin{rema} \label{r4gg0} \leavevmode\\
 \ \ \ 
By the definition of a fermionic twist, the set of isomorphism classes of  fermionic twists of $X^\circledS$ corresponds bijectively to the set $\mr{Im} (H^1_{\text{\'{e}t}} ( \eta_{X^\circledS}))$.
In particular, if $X_b$ (as well as $X_t$) is a  scheme of finite type over $k$ (where $k$ is a separably closed field or a finite field), then  there are only a  finite number of isomorphism classes of  fermionic twists of $X^\circledS$.
Also, if $H^1_{\text{\'{e}t}} (X_b, \mu_2) =0$ (e.g., $X_t$ is simply connected) or $X^\circledS$ is a scheme (i.e., $\mcO_{X_f} =0$), then all fermionic twists of $X^\circledS$ are isomorphic.
   \end{rema}
%-----------------------------------------------------------------------[end remark]-------------------
%\vspace{5mm}

 Consider  a relation ``$\stackrel{f}{\sim}$" in the set of locally noetherian superschemes defined as follows:
 \begin{align} \label{EE07}
 Y^\circledS \stackrel{f}{\sim} Z^\circledS \stackrel{\mr{def}}{\Longleftrightarrow} \text{$Y^\circledS$ is isomorphic to a fermionic twist of $Z^\circledS$}.
 \end{align}
One verifies immediately that this relation forms  an equivalence relation.
The following proposition is one direction of the equivalence in Theorem A.
%If $Y$ and $Z$ are schemes, then $Y \stackrel{f}{\sim} Z$ if and only if $Y \cong Z$.
% (cf. Remark \ref{r4gg}).

%-----------------------------------------------------------------------[begin proposition]------------------
\vspace{3mm}
\bpr \label{p23u}\leavevmode\\
 \ \ \ 
Let   $X^\circledS$ and $Y^\circledS$ be two  locally noetherian superschemes and suppose that $X^\circledS \stackrel{f}{\sim} Y^\circledS$.
%a fermionic twist of $X^\circledS$.
%superscheme with $Y^\circledS \stackrel{f}{\sim} X^\circledS$.
%which is isomorphic (as a superscheme) to a fermionic  twist of $X^\circledS$.
%$Y$ be an {\it fmt} of $X$ (i.e., $X  \stackrel{{\it fmt}}{\sim} Y$).
Then, there exists an equivalence of categories $\mfS \mfc \mfh_{/X^\circledS}^\circledS \isom  \mfS \mfc \mfh_{/Y^\circledS}^\circledS$.
\epr
%-----------------------------------------------------------------------[begin proof]-------------------
\begin{proof}
Let $a \in H^1_{\text{\'{e}t}} (X_b, \mu_2)$.
Suppose that we are given a morphism $f^\circledS :Y^\circledS \migi X^\circledS$  in $\mfS \mfc \mfh_{/X^\circledS}^\circledS$.
Then, the homomorphism $H^1_{\text{\'{e}t}} (X_b, \mu_2) \migi H^1_{\text{\'{e}t}} (Y_b, \mu_2)$ induced by $f_b$ sends $a$ to an element of $ H^1_{\text{\'{e}t}} (Y_b, \mu_2)$; we write, by abuse of notation, for ${^a Y}^\circledS$ the fermionic twist of $Y^\circledS$ associated with this element.
It follows from the functoriality of $(-1)_{X^\circledS}$ (with respect to $X^\circledS$) that
$f^\circledS$ induces  a morphism ${^a f}^\circledS : {^a Y}^\circledS \migi {^a X}^\circledS$ in $\mfS \mfc \mfh_{/{^a X}^\circledS}^\circledS$. 
The assignment $Y^\circledS \mapsto {^a Y}^\circledS$ is functorial, and hence, defines a functor 
\begin{align} \label{R02}
\mfS \mfc \mfh_{/X^\circledS}^\circledS \migi \mfS \mfc \mfh_{/{^a X}^\circledS}^\circledS.
\end{align}
% and denote by $X_a^\circledS$
Since $X^\circledS$ is fermionic twist of ${^a X}^\circledS$ associated with $-a$ (under the  identification $H^1_{\text{\'{e}t}} (X_b, \mu_2) = H^1_{\text{\'{e}t}} ({^a X}_b, \mu_2)$), the discussion just discussed gives rise to  a functor $ \mfS \mfc \mfh_{/{^a X}^\circledS}^\circledS \migi \mfS \mfc \mfh_{/X^\circledS}^\circledS$, which becomes the inverse to the functor (\ref{R02}).
This completes the proof of Proposition \ref{p23u}.
%$\mfS \mfc \mfh_{/X^\circledS}^\circledS \migi \mfS \mfc \mfh_{/{^a X}^\circledS}^\circledS$.
%Thus, we have proved the following proposition, which is one direction of the equivalence in Theorem A. 
\end{proof}
\vspace{3mm}
%-----------------------------------------------------------------------[end proposition]-------------------

%-----------------------------------------------------------------------[end definition]-------------------
\begin{rema} \label{rg0df} \leavevmode\\
 \ \ \ 
% Let $X^\circledS$ be a superscheme.
Let us consider an analogous assertion of Proposition \ref{p23u} where $\mfS \mfc \mfh_{X^\circledS}^\circledS$ is replaced with the category of $\mcO_{X^\circledS}$-supermodules.
 We shall define  
 \begin{align}
 \mcO_{X^\circledS} \text{-} \mfm \mfo \mfd
 \end{align}
   to be the category defined as follows:
 \begin{itemize}
 \item[$\bullet$]
 the {\it objects} are 
 % $\mbZ/2 \mbZ$-graded
  $\mcO_{X^\circledS}$-supermodules $\mcF := \mcF_b \oplus \mcF_f$; 
\item[$\bullet$]
the {\it morphisms} from $\mcF := \mcF_b \oplus \mcF_f$ to $\mcF' := \mcF_b \oplus \mcF'_f$ (where both $\mcF$ and $\mcF'$ are objects in this category) are $\mcO_{X^\circledS}$-linear morphisms $h : \mcF \migi \mcF'$ {\it preserving parity}, i.e., satisfying that $h (\mcF_b) \subseteq \mcF_b'$ and $h (\mcF_f) \subseteq \mcF'_f$.
 \end{itemize}
One verifies that $\mcO_{X^\circledS} \text{-} \mfm \mfo \mfd$ forms an abelian category.
Now, let us take $Y^\circledS := {^a X}^\circledS$ for some $a \in H^1_{\text{\'{e}t}}(X_b, \mu_2)$.
By applying a procedure similar to the procedure in the proof of Proposition \ref{p23u},  one may construct, from each 
%$\mbZ/2 \mbZ$-graded 
$\mcO_{X^\circledS}$-supermodule $\mcF$, an $\mcO_{Y^\circledS}$-supermodule ${^a \mcF}$.
For instance, if $\mcF$ is locally free of finite rank and $\mbV (\mcF)^\circledS$ denotes the superscheme over $X^\circledS$  representing  $\mcF$, then $\mbV ({^a \mcF})^\circledS$ is isomorphic to ${^a \mbV} (\mcF)^\circledS$.
The assignment $\mcF \mapsto {^a \mcF}$ is functorial,  and moreover, determines an equivalence of categories $\mcO_{X^\circledS} \text{-} \mfm \mfo \mfd \isom \mcO_{Y^\circledS} \text{-} \mfm \mfo \mfd$.
%\begin{align}
%\mcO_{X^\circledS} \text{-} \mfm \mfo \mfd \isom \mcO_{Y^\circledS} \text{-} \mfm \mfo \mfd.
%\end{align}
%If we are given an $\mbZ/2 \mbZ$-graded $\mcO_{X^\circledS}$-module $\mcF$, then the endomorphism of $\mcF  := \mcF_b \oplus \mcF_f$ given by  $(a, b) \mapsto (a, -b)$ induces an isomorphism
%$(-1)_{X^\circledS} \mcF \isom \mcF$
Consequently, we conclude the assertion that
\begin{align}
X^\circledS \stackrel{f}{\sim} Y^\circledS \ \ \text{implies that} \ \ \mcO_{X^\circledS} \text{-} \mfm \mfo \mfd \cong \mcO_{Y^\circledS} \text{-} \mfm \mfo \mfd,
\end{align}
which may be thought of as an analogue of Proposition \ref{p23u}.
%the relation $X^\circledS \stackrel{f}{\sim} Y^\circledS$ implies that $\mcO_{X^\circledS} \text{-} \mfm \mfo \mfd \cong \mcO_{Y^\circledS} \text{-} \mfm \mfo \mfd$.
If  
 $\mcO_{X^\circledS} \text{-} \mfm \mfo \mfd$ contained $\mcO_{X^\circledS}$-linear morphisms which does not preserve parity, then  
%need not to preserve parity, then 
there would not be a natural way of construction  of a  functor $\mcO_{X^\circledS} \text{-} \mfm \mfo \mfd \migi  \mcO_{Y^\circledS} \text{-} \mfm \mfo \mfd$ as above.
In particular, $\mcO_{X^\circledS} \text{-} \mfm \mfo \mfd$ may not be equivalent to $\mcO_{Y^\circledS} \text{-} \mfm \mfo \mfd$ even if $Y^\circledS$ is equivalent to $X^\circledS$ (i.e, $X^\circledS \stackrel{f}{\sim} Y^\circledS$).
In other wards, the category of $\mcO_{X^\circledS}$-supermodule in which the  morphisms need not to preserve parity (hence, which is $(\mbZ/2\mbZ)$-graded) may have information which allow us to distinguish $X^\circledS$ from  superschemes equivalent to $X^\circledS$. 
Indeed, the tensor triangulated categories used in the category-theoretic reconstruction of superschemes executed by 
U. V. Dubey and V. M. Malick in ~\cite{D} are assumed to admits a structure of $(\mbZ/2 \mbZ)$-gradation; this assumption will be essential in the reconstruction of the isomorphism classes  (not only the equivalence classes) of superschemes.
   \end{rema}
%-----------------------------------------------------------------------[end remark]-------------------
%\vspace{5mm}

\vspace{5mm}
%----------------------------------------------------------------------[begin subsection]-------------
\subsection{Fermionic twists in the Zariski topology} \label{S136}
\leavevmode\\
\vspace{-4mm}

%-----------------------------------------------------------------------[begin example]------------------
%\begin{exa}
%\leavevmode\\
 \ \ \ 
 Denote by $(\mbG_m)_{X_b}$ the \'{e}tale sheaf on $X_b$  represented by the multiplicative group $\mbG_m$.
 The Kummer sequence 
 \begin{align}
 0 \migi (\mu_2)_{X_b} \migi  (\mbG_m)_{X_b} & \migi (\mbG_m)_{X_b} \migi 0 \\
  a \hspace{5mm}& \mapsto \hspace{5mm} a^2
 \end{align}
 induces 
 an exact sequence
 \begin{align}
 0 \migi \mu_2 \migi \Gamma (X_b, \mcO_{X_b}^\times) & \migi \Gamma (X_b, \mcO_{X_b}^\times) \stackrel{\delta}{\migi} H^1_{\text{\'{e}t}} (X_b, \mu_2) \stackrel{\sigma}{\migi} & \mr{Pic} (X_b) & \migi \mr{Pic} (X_b) \\
 a  \hspace{7mm} & \mapsto \hspace{7mm} a^2 & [\mcL]  \hspace{4mm} & \mapsto   \hspace{2mm} [\mcL^{\otimes 2}]. \notag
 \end{align}
 Any element of  $H^1_{\text{\'{e}t}} (X_b, \mu_2)$ may be represented by a collection of data
 \begin{align} \label{DD02}
s :=  (\{ U_\alpha \}_{\alpha \in I}, \{ s_\alpha \}_{\alpha \in I}, \{ t_{\alpha, \beta} \}_{(\alpha, \beta) \in I_2}),
 \end{align} 
 where
 \begin{itemize}
 \item[$\bullet$]
 $I$ is an index set;
 \item[$\bullet$]
 $\{ U_\alpha \}_{\alpha \in I}$ is a {\it Zariski} open covering of $X_b$;
 \item[$\bullet$]
 each $s_\alpha$ ($\alpha \in I$) is an element of $\Gamma (U_\alpha, \mcO_{U_\alpha}^\times)$;
 \item[$\bullet$]
 $I_2 := \{ (\alpha, \beta ) \in I \times I \ | \ U_{\alpha, \beta} := U_\alpha \cap U_\beta \neq \emptyset)$;
 \item[$\bullet$]
  $\{ t_{\alpha, \beta} \}_{(\alpha, \beta) \in I_2}$  is a $1$-cocycle of $\{ U_\alpha \}_{\alpha \in I}$ with coefficients in $\mcO_{X_b}^\times$ 
%  ($(\alpha, \beta) \in I_2$)
%   is an element of $\Gamma (U_{\alpha, \beta}, \mcO_{U_{\alpha, \beta}}^\times)$ 
   satisfying the equality  $s_\beta |_{U_{\alpha, \beta}} \cdot t^2_{\alpha, \beta} = s_\alpha |_{U_{\alpha, \beta}}$ for any $(\alpha, \beta) \in I_2$.
 \end{itemize}
The homomorphism $\delta$ (resp., $\sigma$) is given by assigning $a \mapsto (\{ X_b \}, \{ a \}, \{ 1 \})$ for any $a \in \Gamma (X_b, \mcO_{X_b}^\times)$ (resp., $\overline{s} \mapsto \overline{(\{ U_\alpha \}_{\alpha}, \{ t_{\alpha, \beta}\}_{\alpha, \beta} )}$  for any $s$ as in (\ref{DD02})).

%\vspace{5mm}
%----------------------------------------------------------------------[begin subsection]-------------
%\subsection{Fermionic twists in the Zariski topology} \label{S137}
%\leavevmode\\
%\vspace{-4mm}

Now,  let $u \in \Gamma (X_b, \mcO_{X_b}^\times)$.
%If there is no fear of confusion, then w
We shall write
\begin{align}
{^u X}^\circledS := {^{\delta (u)} X}^\circledS.
\end{align}
by abuse of notation.
%Let 
%\begin{align}
%{^s X}^\circledS
%\end{align}
One verifies that it is  a unique (up to isomorphism) superscheme such that the triple $\mcA_{{^u X}^\circledS}$ associated with it (cf. Proposition \ref{p24pv})  coincides with $(X_b, \mcO_{X_f}, u \cdot m_{X^\circledS})$.
%Then, ${^s X}^\circledS$ forms a fermionic twist of $X^\circledS$.
(In particular, $\mcO_{{^u X}^\circledS} = \mcO_{X^\circledS}$ as an $\mcO_{X_b}$-module.)
Indeed, 
let us write $Y^\circledS$ for the superscheme corresponding to $(X_b, \mcO_{X_f}, u \cdot m_{X^\circledS})$ (hence, $Y_b = X_b$). 
Also, let us take an \'{e}tale covering $U \migi X_b$ such  that there exists $v  \in \Gamma (U, \mcO_{U}^\times)$ with $v^2 = u$.
%(The pair $(\{ U \}, \{ u\})$ represents the element $\delta (u) \in H^1_{\text{\'{e}t}} (X_b, \mu_2)$.)
%to verify this fact, we may assume (after possibly replacing $X^\circledS$ with $X^\circledS \times_{X_b} U$ for some 
The automorphism of 
the $\mcO_{U}$-module  $\mcO_{U} \oplus \mcO_{X_f} |_{U}$ given by assigning $(a, \epsilon_a) \mapsto (a, v \cdot \epsilon_a)$ determines an isomorphism
%\begin{align}
%\end{align}
$X^\circledS |_{U} \isom Y^\circledS |_{U}$ that induces the identity morphism of $X_b$.
This implies that  $Y^\circledS$ is the  fermionic  twist of $X^\circledS$ associated with $\delta (u)$, as desired.

Conversely,  any  fermionic  twist of $X^\circledS$ is, {\it Zariski locally on $X_b$}, isomorphic to ${^u X}^\circledS$ (for some local section $u \in \mcO_{X_b}^\times$), as  described in the following proposition.
%Indeed,  the image of any element of $H^1_{\text{\'{e}t}} (X_b, \mu_2)$ via the homomorphism $H^1_{\text{\'{e}t}} (X_b, \mu_2) \migi \mr{Pic} (X_b)$ induced by the injection $(\mu_2)_{X_b} \migi (\mbG_m)_{X_b}$)
%\end{exa}
%-----------------------------------------------------------------------[end example]-------------------

%-----------------------------------------------------------------------[begin proposition]------------------
\vspace{3mm}
\bpr \label{p23}\leavevmode\\
 \ \ \ 
Let   $a$ be an element of $H^1_{\text{\'{e}t}} (X_b, \mu_2)$ 
(hence, we have a fermionic twist ${^a X}^\circledS$ of $X^\circledS$ associated with $a$).
Also, let
%We shall   choose a representative 
$(\{ U_{\alpha \in I} \}_{\alpha \in I}, \{ s_\alpha \}_{\alpha \in I}, \{ t_{\alpha, \beta} \}_{(\alpha, \beta) \in I_2})$ be a representative of $a$ as in (\ref{DD02}).
Then, there exists a collection of isomorphisms
%$\{ \xi^\circledS_{\alpha} \}_{\alpha \in I}$ of isomorphisms 
\begin{align}
\{ \xi^\circledS_{\alpha} : {^a X}^\circledS |_{U_\alpha} \isom {^{s_\alpha} X}^\circledS |_{U_\alpha} \}_{\alpha \in I}
\end{align} 
%of superschemes  
satisfying the following two conditions:
\begin{itemize}
\item[$\bullet$]
For each $\alpha \in I$,  the morphism $(\xi_\alpha)_b$ of schemes  underlying $\xi^\circledS_\alpha$  coincides with  the identity morphism of $U_\alpha$;
\item[$\bullet$]
For each $(\alpha, \beta) \in I_2$, the automorphism 
\begin{align}
\xi^\circledS_\beta \circ (\xi^{\circledS}_\alpha )^{-1} : {^{s_\alpha}X}^\circledS |_{U_{\alpha, \beta}}\isom  {^{s_\beta}X}^\circledS |_{U_{\alpha, \beta}}
\end{align}
  corresponds to the automorphism of the $\mcO_{U_{\alpha, \beta}}$-module $\mcO_{U_{\alpha, \beta}} \oplus \mcO_{X_f} |_{U_{\alpha, \beta}}$ given by assigning $(a, \epsilon_a) \mapsto (a, t_{\alpha, \beta} \cdot \epsilon_a)$.
\end{itemize}
%$Y^\circledS$ be a fermionic twist of $X^\circledS$.
%Then, there exists an equivalence of categories $\mfS \mfc \mfh_{/X^\circledS}^\circledS \isom  \mfS \mfc \mfh_{/Y^\circledS}^\circledS$.
\epr
%-----------------------------------------------------------------------[begin proof]-------------------
\begin{proof}
The assertion follows immediately from the definition of a fermionic twist and the above discussion. 
\end{proof}
\vspace{3mm}
%-----------------------------------------------------------------------[end proposition]-------------------

%-----------------------------------------------------------------------[begin remark]------------------
%\begin{rema} \label{r4gg1} \leavevmode\\
% \ \ \ 

%   \end{rema}
%-----------------------------------------------------------------------[end remark]-------------------
%\vspace{5mm}

\vspace{5mm}
%----------------------------------------------------------------------[begin subsection]-------------
\subsection{$\mbA^{0 \mid 1}$-twists } \label{S15}
\leavevmode\\
\vspace{-4mm}

%-----------------------------------------------------------------------[begin remark]------------------
%\begin{rema} \label{r1} \leavevmode\\
 %\ \ \ 
% Toward the proof of Theorem A, we shall introduce some notation.
For each pair $(n ,m)$ of nonnegative integers,   we shall denote by 
\begin{equation}
\mbA^{n\mid m}
\end{equation}
 the $(n|m)$-dimensional affine superspace over $\mbZ [\frac{1}{2}]$, i.e.,
  the superspectrum of the superring $\mbZ[\frac{1}{2}] [t_1, \cdots, t_n, \psi_1, \cdots, \psi_m]$,  where the  $t_1,  \cdots, t_n$ are  ordinary indeterminates and $\psi_1, \cdots, \psi_m$ are  odd indeterminates.
  %   and write $\mbA_X^{1\mid 1} := X \times \mbA^{1\mid 1}$.
%We have two  morphisms
%\begin{equation}
%\delta : \mbA^{0| 0} \migi \mbA^{0 | 1}, \ \ \ \pi : \mbA^{1\mid 1} \migi \mbA^{0 | 1},
%\end{equation}
%where the first morphism is given by $\psi \ (:= \psi_1)\mapsto 0$ and the second morphism denotes  the natural projection.
Also,  let us write 
 \begin{align}  \label{e43}
 \mbA^{n  | m}_{X^\circledS} := X^\circledS \times \mbA^{n | m}.
 \end{align}
 %and write
% \begin{align} \label{EE01}
% \delta_{X^\circledS} := \mr{id}_{X^\circledS} \times \delta : X^\circledS \migi \mbA_{X^\circledS}^{0 |1},  \  \ \pi_{X^\circledS} : = \mr{id}_{X^\circledS} \times \pi : \mbA^{1|1}_{X^\circledS} \migi \mbA_{X^\circledS}^{0|1}.
% \end{align}
For any  $Y^\circledS \in \mr{Ob} (\mfS \mfc \mfh^\circledS_{/ X^\circledS})$ and any nonnegative integers $n$, $m$, the superscheme
$\mbA^{n  | m}_{Y^\circledS}$ belongs to  $\mr{Ob} (\mfS \mfc \mfh^\circledS_{/ X^\circledS})$.
Also,
we have a sequence of 
 functorial (in $Y^\circledS$) bijections of sets:
\begin{align}  \label{e5}
& \ \ \ \ \,  \mr{Map}_{\mfS \mfc \mfh^\circledS_{/ Y^\circledS}} (Y^\circledS, \mbA_{Y^\circledS}^{1\mid 1}) \\
& \isom
\mr{Map}_{\mfS \mfc \mfh^\circledS_{/ Y^\circledS}} (Y^\circledS, \mbA_{Y^\circledS}^{1 \mid 0} \times_{Y^\circledS} \mbA^{0 \mid 1}_{Y^\circledS}) \notag  \\
& \isom \mr{Map}_{\mfS \mfc \mfh^\circledS_{/ Y^\circledS}} (Y^\circledS, \mbA_{Y^\circledS}^{1\mid 0})  \times \mr{Map}_{\mfS \mfc \mfh^\circledS_{/ Y^\circledS}} (Y^\circledS, \mbA_{Y^\circledS}^{0\mid 1})  \notag \\  
& \isom \Gamma (Y_b, \mcO_{Y_b}) \times \Gamma (Y_b, \mcO_{Y_f}) \notag \\
& \isom \Gamma (Y_b, \mcO_{Y^\circledS}), \notag
\end{align}
where the third bijection is given by $(h_1^\circledS, h_2^\circledS) \mapsto (h_1^\flat (t), h_2^\flat (\psi))$.
% \begin{equation}
% \Gamma (X_b, \mcO_{X^\circledS}) \isom \mr{Map}_{\mfS \mfc \mfh^\circledS_{/ X^\circledS}} (X^\circledS, \mbA_{X^\circledS}^{1\mid 1}).
% \end{equation}
% which sends the zero section  to $\delta_{X^\circledS}$. 
 %  between the set $\Gamma (X_b, \mcO_X)$ and the set of sections $X \migi X \times \mbA^{1\mid 1}$ of the projection  $X \times \mbA^{1\mid 1} \migi X$.
 The multiplication and addition in $\Gamma (Y_b, \mcO_{Y^\circledS})$ correspond, via (\ref{e5}),  to   morphisms 
 \begin{align}  \label{DD08}
 \mu_{Y^\circledS} : \mbA_{Y^\circledS}^{1\mid 1} \times_{Y^\circledS} \mbA_{Y^\circledS}^{1\mid 1} \migi \mbA_{Y^\circledS}^{1\mid 1}
 \ \  \text{and} \ \  
 \alpha_{Y^\circledS} : \mbA_{Y^\circledS}^{1\mid 1} \times_{Y^\circledS} \mbA_{Y^\circledS}^{1\mid 1} \migi \mbA_{Y^\circledS}^{1\mid 1}
 \end{align}
  respectively.
That is to say, the set $\mr{Map}_{\mfS \mfc \mfh_{/Y^\circledS}^\circledS} (Y^\circledS, \mbA_{Y^\circledS}^{1\mid 1})$   admits a structure of superring by means of $\mu_{Y^\circledS}$ and $\alpha_{Y^\circledS}$ (and the decomposition $\mbA_{Y^\circledS}^{1\mid 1} \isom \mbA_{Y^\circledS}^{1\mid 0} \times_{Y^\circledS} \mbA_{Y^\circledS}^{0\mid 1}$), and 
the composite bijection (\ref{e5}) becomes an isomorphism of superrings.
%Moreover, various $\mu_{Y^\circledS}$'s (resp.,  $\alpha_{Y^\circledS}$'s) (for quasi-compact open subsuperschemes $Y^\circledS$ of $X^\circledS$) may be glued together to a morphism
% \begin{align}  \label{DD089}
% \mu_{X^\circledS} : \mbA_{X^\circledS}^{1\mid 1} \times_{X^\circledS} \mbA_{X^\circledS}^{1\mid 1} \migi \mbA_{X^\circledS}^{1\mid 1}
% \ \  (\text{resp}.,  \  
% \alpha_{X^\circledS} : \mbA_{X^\circledS}^{1\mid 1} \times_{X^\circledS} \mbA_{X^\circledS}^{1\mid 1} \migi \mbA_{X^\circledS}^{1\mid 1}), 
% \end{align}
%which corresponds to the multiplication (resp., addition) in $\Gamma (X_b, \mcO_{X^\circledS})$.
In particular,  each  element $a$  of $\Gamma (Y_b, \mcO_{Y_b})$ corresponds to  a morphism
\begin{equation} \label{e75}
\sigma^{[a]}_{Y^\circledS} : Y^\circledS \migi \mbA_{Y^\circledS}^{1\mid 0}. 
\end{equation}
Denote by  $\mcA u t_{Y^\circledS} (\mbA_{Y^\circledS}^{1\mid 0}, \sigma^{[0]}_{Y^\circledS})$
the Zariski sheaf on $Y_b$ which, to any  open subsuperscheme $U$ of $Y_b$, assigns the group of automorphisms of $\mbA_{Y^\circledS |_U}^{1\mid 0}$ over $Y^\circledS |_U$ which are  compatible with $ \sigma^{[0]}_{Y^\circledS} |_{Y^\circledS |_U}$.
% to $X^\circledS |_U$.
%If $(\mbG_m)$
%If we consider, in a natural manner,  $\mcO_{X_b}^\times$ as an \'{e}tale sheaf on $X_b$, then the homomorphism
The homomorphism
\begin{align}
\mcO_{Y_b}^\times \isom \mcA u t_{Y^\circledS} (\mbA_{Y^\circledS}^{1\mid 0}, \sigma^{[0]}_{Y^\circledS})
\end{align}
which, to any local section $a \in \mcO_{Y_b}^\times$, assigns the automorphism of $\mbA_{Y^\circledS}^{1\mid 0}$ over $Y^\circledS$ determined by $\psi \mapsto a \cdot \psi$ turns out to be  bijective.
By applying the functor $H^1_{\mr{Zar}} (Y_b, -)$. we have an isomorphism
\begin{align} \label{DD01}
\mr{Pic}(Y_b) \isom H^1_{\mr{Zar}} (Y_b, \mcA u t_{Y^\circledS} (\mbA_{Y^\circledS}^{1\mid 0}, \sigma^{[0]}_{Y^\circledS})).
\end{align}

%For each $a \in \mcO_{X_b}^\times$

%-----------------------------------------------------------------------[begin definition]------------------
\vspace{3mm}
\bde \label{d363}\leavevmode\\
\vspace{-5mm}
\begin{itemize}
\item[(i)]
   An {\bf $\mbA^{0 \mid 1}$-twist over  $Y^\circledS$} is
a twisted form of $(\mbA_{Y^\circledS}^{1\mid 0}, \sigma^{[0]}_{Y^\circledS})$ (over the Zariski  topology on $Y_b$) determined, via (\ref{DD01}),  by some $a \in \mr{Pic} (Y_b)$; it may be described as 
 a pair 
\begin{equation}
(Z^\circledS,  \sigma_{Z^\circledS/Y^\circledS} )
\end{equation}
consisting of a twisted form $Z^\circledS$ of $\mbA_{Y^\circledS}^{1\mid 0}$ over $Y^\circledS$ and a section $ \sigma_{Z^\circledS/Y^\circledS} : Y^\circledS \migi Z^\circledS$ of the structure morphism of $Z^\circledS$.
We shall refer to the pair 
$(Z^\circledS,  \sigma_{Z^\circledS/Y^\circledS} )$ as the {\bf  $\mbA^{0 \mid 1}$-twist over  $Y^\circledS$ associated with $a$}.
\item[(ii)]
Let $(Z^\circledS, \sigma_{Z^\circledS/Y^\circledS})$ and $(Z'^\circledS, \sigma_{Z'^\circledS/Y^\circledS})$ be two $\mbA^{0 \mid 1}$-twists over $Y^\circledS$.
An {\bf isomorphism of $\mbA^{0 \mid 1}$-twists} from $(Z^\circledS, \sigma_{Z^\circledS/Y^\circledS})$ to $(Z'^\circledS, \sigma_{Z'^\circledS/Y^\circledS})$ is an isomorphism $h^\circledS : Z^\circledS \isom Z'^\circledS$ of superschemes over $Y^\circledS$ with $h^\circledS \circ \sigma_{Z^\circledS/Y^\circledS} = \sigma_{Z'^\circledS/Y^\circledS}$.
\end{itemize}
 \ede
\vspace{3mm}
%-----------------------------------------------------------------------[end definition]-------------------

By (\ref{DD01}), there exists canonically  a bijective correspondence between $\mr{Pic} (Y_b)$ and the set of isomorphism classes of $\mbA^{0 \mid 1}$-twists over  $Y^\circledS$.

\vspace{5mm}
%----------------------------------------------------------------------[begin subsection]-------------
\subsection{The multiplication morphisms of fermionic twists} \label{sub5}
\leavevmode\\
\vspace{-4mm}

Let  $u \in \Gamma (Y_b, \mcO_{Y_b}^\times)$.
% and  $X^s$ denotes the fermion-multiplication deformation of $X$ by means of $s$.
Since $\mcO_{{^u Y}^\circledS} = \mcO_{Y^\circledS}$ as $\mcO_{Y_b}$-modules, the multiplication in $\mcO_{{^u Y}^\circledS}$ gives rise to
%The multiplication morphism $\mu_{{^u X}^\circledS} : \mbA_{{^u X}^\circledS}^{1\mid 1} \times_{{^u X}^\circledS} \mbA_{{^u X}^\circledS}^{1\mid 1} \migi \mbA_{{^u X}^\circledS}^{1\mid 1}$ of the fermionic twist ${^u X}^\circledS$ of $X^\circledS$ associated with $u$  
%corresponds, via the bijection (\ref{e5}), to
%may be thought of as 
a morphism
\begin{equation} \label{DD010}
\mu_{Y^\circledS \rightsquigarrow {^u Y}^\circledS} : \mbA_{Y^\circledS}^{1 \mid 1} \times_{Y^\circledS} \mbA_{Y^\circledS}^{1 \mid 1} \migi \mbA_{Y^\circledS}^{1 \mid 1}
\end{equation}
over $Y^\circledS$ under  the bijection (\ref{e5}).
The morphism $\mu_{Y^\circledS \rightsquigarrow {^u Y}^\circledS}$ corresponds to 
%the $\mcO_{X_b}$-linear 
the homomorphism of superalgebras over $\mcO_{Y^\circledS}$ described as follows:
\begin{align} \label{f77}
\mcO_{Y^\circledS}  [t, \psi] &\migi\mcO_{Y^\circledS}  [t, \psi]  \otimes_{\mcO_{Y^\circledS}} \mcO_{Y^\circledS}  [t, \psi]  \\
t \hspace{5mm}& \mapsto \hspace{8mm} t \otimes t + s \cdot \psi \otimes \psi  \notag  \\
\psi \hspace{5mm} &\mapsto \hspace{7mm} \psi \otimes t + t \otimes \psi . \notag 
\end{align}
%Then, the morphism corresponding to the morphism \ref{f77} coincides with the multiplication in $X^s$.

%$t \mapsto \phi \otimes t + t \otimes \phi$ and $\theta \mapsto t \otimes t + s \cdot \phi \otimes \phi$.

Next, let 
$a$ be an element of $H^1_{\text{\'{e}t}} (Y_b, \mu_2)$ and let $Z^\circledS := {^a Y}^\circledS$.
 We shall choose  a representative 
$(\{ U_\alpha \}_{\alpha \in I}, \{ s_\alpha \}_{\alpha \in I}, \{ t_{\alpha, \beta} \}_{(\alpha, \beta) \in I_2})$ of $a$ as in (\ref{DD02}) (where $X^\circledS$ is replaced with $Y^\circledS$).
%$Y^\circledS := (Y_b, \mcO_{Y^\circledS})$ be  a fermonic twist of $X$ 
%and fix a deformation data $\mcD_Y := (\iota, \{U_\alpha \}_{\alpha \in I}, \{ s_\alpha \}_{\alpha \in I}, \{ \xi_\alpha \}_{\alpha \in I})$  for $Y$.
%If  $U_{\beta \alpha} := U_\alpha \cap U_\beta \neq \emptyset$, 
%then the isomorphism $\xi_\beta \circ \xi_\alpha^{-1} : (X |_{U_{\beta \alpha}})^{s_\alpha |_{U_{\beta \alpha}}} \isom (X |_{U_{\beta \alpha}})^{s_\beta |_{U_{\beta \alpha}}}$
%corresponds to
%the automorphism of $\mcO_{U_{\beta \alpha}} \oplus \mcO_{X_f} |_{U_{\beta \alpha}}$ given by
%$(a, \epsilon_a) \mapsto (a, t_{\beta \alpha} \cdot \epsilon_a)$ for some $t_{\beta \alpha} \in \Gamma (U_{\beta \alpha}, \mcO_{U_{\beta \alpha}}^\times)$.
Write 
\begin{equation} \label{e57}
(\mbA_{Y^\circledS \leadsto Z^\circledS}^{0 \mid 1}, \sigma_{\mbA_{Y^\circledS \leadsto Z^\circledS}^{0 \mid 1}})
% \ (\text{resp.,} \ \mbA_{X \leadsto Y}^{1\mid 1})
\end{equation}
for the $\mbA^{0 \mid 1}$-twist over  $Y^\circledS$ determined by $\sigma (a) \in \mr{Pic} (Y_b)$,
%corresponding, via the bijection obtained in Remark \ref{r7}, to the collection of data $(\{ U_\alpha\}, \{t_{\beta \alpha} \}_{\alpha, \beta})$, 
and write  
\begin{equation}
\mbA_{Y^\circledS \leadsto Z^\circledS}^{1\mid 1} := \mbA_{Y^\circledS \leadsto Z^\circledS}^{0 \mid 1} \times \mbA^{1 \mid 0}. \end{equation}
%Now, we consider the automorphism
%$\xi_{\beta \alpha}^{0 \mid 1}$ of $\mbA^{0\mid 1}_{X |_{U_{\beta \alpha}}}$ over $X |_{U_{\beta \alpha}}$ (where, for each open subscheme $U$ of $X$, we shall write  $X |_{U} := X \times_{X_b} U$)
%corresponding to
%the $\mcO_X$-linear  automorphism of $\mcO_X  \oplus \mcO_X \psi$ determined by $\psi \mapsto t_{\beta \alpha} \cdot \psi$.
%we consider the automorphism of $U_{\alpha \beta} \times_\mbZ \mbA^{1\mid 1}_\mbZ$ given by 
%\begin{align}
%\mcO_X [t] \oplus \mcO_X [t] \psi \isom  \mcO_X [t] \oplus \mcO_X [t] \psi  \\
%(a, b \phi) \mapsto (a, t_{\alpha\beta}) \phi \notag
%\end{align}
%Then, the super-schemes $\{  \mbA^{0\mid 1}_{X |_{U_\alpha}}\}_{\alpha\in I}$ may be glued, by means of the isomorphisms $\{\xi_{\beta \alpha}^{0 \mid 1} \}_{\alpha \in I}$, to a super-scheme  $\mbA_{Y/X}^{0\mid 1}$  over $X$, which is Zariski locally isomorphic to $\mbA^{0\mid 1}_X$.
%It follow from the definition of an {\it mod} that there exist a super-scheme 
%\begin{equation}
%\mbA_Y^{0 \mid 1}
%\end{equation}
% over $Y$ such that it is Zariski  locally  isomorphic to  $Y \times_\mbZ \mbA_\mbZ^{0 \mid 1}$ and there exists an isomorphism $\mbA_Y^{0 \mid 1} \times_\mbZ \mbA \isom \mbA^{1\mid 1}_Y$.
The multiplication morphisms $\mu_{Y^\circledS |_{U_\alpha}\leadsto {^{s_\alpha}Y}^\circledS |_{U_\alpha}}$ ($\alpha \in I$)
may be glued together to a morphism
%  $\mu_Y$ of $Y$ determines  a morphism
%We shall denote by 
\begin{equation}
\mu_{Y^\circledS \leadsto Z^\circledS} : \mbA_{Y^\circledS \leadsto Z^\circledS}^{1 \mid 1} \times_{Y^\circledS} \mbA_{Y^\circledS \leadsto Z^\circledS}^{1\mid 1} \migi \mbA_{Y^\circledS \leadsto Z^\circledS}^{1\mid 1}
\end{equation}
%the morphism corresponding to the
over $Y^\circledS$.
This morphism does not depend on the choice of a representative of $a$.
Also, we  obtain (by glueing together  the morphisms $\alpha_{{^{s_\alpha}Y}^\circledS |_{U_\alpha}}$)  a morphism
\begin{equation} \label{e10}
\alpha_{Y^\circledS \leadsto Z^\circledS} : \mbA_{Y^\circledS \leadsto Z^\circledS}^{1 \mid 1} \times_{Y^\circledS } \mbA_{Y^\circledS \leadsto Z^\circledS}^{1\mid 1} \migi \mbA_{Y^\circledS \leadsto Z^\circledS }^{1\mid 1}
\end{equation}
over $Y^\circledS$.
%Note that 
%$\mbA^{0\mid 1}_{Y/X}$ (as well as $\mbA^{1\mid 1}_{Y/X}$)
%and
The morphism $\alpha_{Y^\circledS \leadsto Z^\circledS}$ depends only on  the $\mbA^{0 \mid 1}$-twist $\mbA_{Y^\circledS \leadsto Z^\circledS}^{0 \mid 1}$ (i.e., the class $\sigma (a) \in \mr{Pic} (Y_b)$).
%(where $\mbA_{Y/X}^{1 \mid 1} = \mbA_{Y/X}^{0 \mid 1} \times \mbA^{1\mid 0}$) over $X$, which is obtained by gluing the morphisms $\{ \mu_{U_\alpha^{s_\alpha}}\}_{\alpha \in I}$ together.
%which is verified to correspond to the multiplication in $\mcO_Y$.
%That is to say, there exists a natural bijection
Owing to the morphisms $\alpha_{Y^\circledS \leadsto Z^\circledS} $ and $\mu_{Y^\circledS \leadsto Z^\circledS}$, we have  an isomorphism of superrings
\begin{equation} \label{e09}
\Gamma (Z_b, \mcO_{Z^\circledS}) \isom \mr{Map}_{\mfS \mfc \mfh^\circledS_{/Y^\circledS}} (Y^\circledS, \mbA_{Y^\circledS \leadsto Z^\circledS}^{1 \mid 1})
\end{equation}
which is functorial with respect to $Y^\circledS \in \mr{Ob} (\mfS \mfc \mfh^\circledS_{/X^\circledS})$.
%base-change   morphisms  $Y^\circledS \migi X^\circledS$ in $\mfS \mfc \mfh_{/ \mbZ [\frac{1}{2}]}^\circledS$.

%  \end{rema}
%-----------------------------------------------------------------------[end remark]-------------------

%%%%%%%%%%%%%%%%%%%%%%%%%%%%%%%%--[ begin  section1]---%%%%%%
\vspace{6mm}
\section{Proof of Theorem A} \vspace{3mm}

This section  is devoted to prove the remaining portion of  Theorem A, i.e.,  that the equivalence class  defined by ``  $\stackrel{f}{\sim}$ " of a  locally noetherian superscheme $X$ may be reconstructed purely category-theoretically from the category $\mfS \mfc \mfh_{/X^\circledS}^\circledS$.
In the following discussion, we will often speak of various properties of objects and   morphisms in $\mfS \mfc \mfh_{/X^\circledS}^\circledS$ as being ``{\it characterized (or reconstructed) category-theoretically}".
By this, we mean that they are preserved by arbitrary equivalences of categories
$\mfS \mfc \mfh_{/X^\circledS}^\circledS \isom \mfS \mfc \mfh_{/X'^\circledS}^\circledS$
(where $X'^\circledS$ is another locally noetherian superscheme).
For instance, the set of monomorphisms in $\mfS \mfc \mfh_{/X^\circledS}^\circledS$ may be  {\it characterized category-theoretically} as the morphisms $f^\circledS : Z^\circledS \migi Y^\circledS$ such that,  for any  $W^\circledS \in \mr{Ob} (\mfS \mfc \mfh_{/X^\circledS}^\circledS)$, the map of sets $\mr{Map}_{\mfS \mfc \mfh_{/X^\circledS}^\circledS} (W^\circledS, Z^\circledS) \migi \mr{Map}_{\mfS \mfc \mfh_{/X^\circledS}^\circledS} (W^\circledS, Y^\circledS)$ given by composing with $f^\circledS$ is  injective.
%there exists a morphism $h^\circledS$ with  which the composites $f^\circledS \circ h^\circledS$,   $h^\circledS \circ f^\circledS$ of $f^\circledS$  coincide with  the identity morphisms.
To simplify notation, however, we omit explicit mention of this equivalence $\mfS \mfc \mfh_{/X^\circledS}^\circledS \isom \mfS \mfc \mfh_{/X'^\circledS}^\circledS$, of $X'$, and of the various ``primed" objects and  morphisms corresponding to the original objects and morphisms, respectively,  in $\mfS \mfc \mfh_{/X^\circledS}^\circledS$.

In this section, {\it let us fix a locally noetherian superscheme $X^\circledS$}.
%suppose that {\it the superscheme $X^\circledS$ is locally noetherian}. 

\vspace{5mm}
%----------------------------------------------------------------------[begin subsection]-------------
\subsection{} \label{sub27}

Our tactics for completing the proof of  Theorem A (i.e., recognizing the structure of superscheme of $X^\circledS$) is, as in ~\cite{Mzk1},  to reconstruct step-by-step various partial information of $X^\circledS$ from the categorical structure of $\mfS \mfc \mfh_{/X^\circledS}^\circledS$.
As the first step, we reconstruct the set of objects in $\mfS \mfc \mfh_{/X^\circledS}^\circledS$ which are isomorphic to spectrums of  fields (cf. Proposition \ref{p8}).
Of course,  these objects allow us to know the points in 
%One knows the underlying set of
 the topological space underlying $X^\circledS$.
 % from these objects.

%Let $R$ be a {\it commutative} ring, $M$ an $R$-module, and $Y$ a super-scheme over $R$.
For each superring $R$, we denote by  
\begin{align}
\mr{Spec} (R)^\circledS
\end{align}
  the superspectrum of $R$.
Let $k$ be a field and  $M$ a finite-dimensional $k$-vector space.
 We shall equip $k \oplus M$ with a structure of superalgebra over $k$ given as follows:
 \begin{itemize}
 \item[$\bullet$]
 The bosonic part is the first factor $k$ and the fermionic part is the second factor $M$;
 \item[$\bullet$]
 The multiplication is given by assigning $(a, \epsilon_a) \cdot (b, \epsilon_b) := (ab, a\epsilon_b + b \epsilon_b)$ for any $a$, $b \in k$ and $\epsilon_a$, $\epsilon_b \in M$.
 \end{itemize}
   % The exterior algebra $\bigwedge_k^\bullet M$ over $k$  associated with $M$  admits naturally a structure of superalgebra over $k$.
 We shall write 
\begin{equation}
\mbA_k^{0 \mid M} :=  \mr{Spec} (k \oplus M)^\circledS.
% \mr{Spec}({\bigwedge}^\bullet_k M)^\circledS,
%\mbA_\mbZ^{m \mid n} := \mr{SSpec} (\bigwedge_R R e_1 \oplus \cdots R e_n)
\end{equation}
%where $\mr{SSpec}(\bigwedge^\bullet_k M)$ denotes  
%i.e.,  the superspectrum  of the exterior algebra $\bigwedge_k^\bullet M$ over $k$ (equipped with a natural $(\mbZ/2 \mbZ)$-gradation) associated with $M$
In other wards,  $\mbA_k^{0 \mid M}$ is a unique (up to isomorphism) superscheme satisfying that
$\mcA_{\mbA_k^{0 \mid M}} := (\mr{Spec}(k), \mcO_{\mr{Spec} (k)}\otimes_k M, 0)$.
 In particular, $\mbA_k^{0 \mid k} = \mbA_{k}^{0 \mid 1}$  (cf. (\ref{e43})).
If $M_1$ and $M_2$ are finite-dimensional $k$-vector spaces, then any morphism $\mbA_k^{0 \mid M_1} \migi \mbA_k^{0 \mid M_2}$ of superschemes over $k$ coincides with  the morphism induced from a   $k$-linear morphism $M_2 \migi  M_1$ which is uniquely determined.
This observation shows the following lemma.

%-----------------------------------------------------------------------[begin proposition]------------------
\vspace{3mm}
\ble \label{Pp3}\leavevmode\\
 \ \ \ 
Let us  write $\mfV \mfe \mfc_k$ for  the {\it opposite} category of finite-dimensional $k$-vector spaces and write 
\begin{align} \label{DD021}
{^\circ \mfS} \mfc \mfh_{/k}^\circledS
\end{align}
 for the full subcategory of $\mfS \mfc \mfh_{/k}^\circledS$ consisting of  superschemes  which are isomorphic to $\mbA_k^{0 \mid M}$ for  some finite-dimensional $k$-vector space $M$.
 Then,  the functor
 % the above observations implies that the assignment $M \mapsto \mbA_k^{0 \mid M}$ defines an equivalence of categories
 \begin{align} \label{DD020}
 \mfV \mfe \mfc_k & \migi  {^\circ \mfS} \mfc \mfh_{/k}^\circledS \\
 M \hspace{2mm} & \mapsto  \hspace{3mm} \mbA_k^{0 \mid M} \notag 
 \end{align}
defines  an equivalence of categories.
   \ele
%-----------------------------------------------------------------------[begin proof]-------------------
%\begin{proof}
%\end{proof}
%\vspace{3mm}
%-----------------------------------------------------------------------[end proposition]-------------------

%-----------------------------------------------------------------------[begin proposition]------------------
\vspace{3mm}
\ble \label{Pp2}\leavevmode\\
 \ \ \ 
 Suppose that  $\mr{Spec}(k)$ 
is an object of $\mfS \mfc \mfh_{/X^\circledS}^\circledS$, in particular, admits a structure morphism $\mr{Spec}(k) \migi X^\circledS$.
%, i.e., is assumed to be equipped with a morphism $\mr{Spec} (k) \migi X^\circledS$
 (Hence, $\mbA_k^{0 \mid M}$ is an object of $\mfS \mfc \mfh_{/X^\circledS}^\circledS$ by taking account of  the  composite  $\mbA_k^{0 \mid M} \migi \mr{Spec}(k) \migi X^\circledS$).
%If $k$ is a field for which there exists, at least,  one point $\mr{Spec}(k) \migi X$ and $V$ is a finite 
%dimensional  $k$-vector space, 
There exists a natural bijection
%we have a natural bijection
\begin{align}\label{f1}
& \mr{Map}_{\mfS \mfc \mfh_{/X^\circledS}^\circledS}(\mbA_k^{0\mid  M}, Y^\circledS) \\
 \isom & \big\{ (s, h) \ \big| \ s \in \mr{Map}_{\mfS \mfc \mfh_{/X_b}} (\mr{Spec}(k), Y_b), h \in \mr{Hom}_k(s^* (\mcO_{Y_f}), M) \big\} \notag 
\end{align}
for any object $Y^\circledS$ of $\mfS \mfc \mfh_{/X^\circledS}^\circledS$.
%where if $\mcC$ is a category and $y$ and $y'$ are objects of $\mcC$, then we shall denote by
%$\mr{Map}_\mcC (y, y')$ the set of morphisms from $y$ to $y'$.
\ele
%-----------------------------------------------------------------------[begin proof]-------------------
\begin{proof}
The assertion follows directly from the definition of $\mbA_k^{0 \mid M}$.
\end{proof}
\vspace{3mm}
%-----------------------------------------------------------------------[end proposition]-------------------

%-----------------------------------------------------------------------[begin proposition]------------------
\vspace{3mm}
\bpr \label{p2}\leavevmode\\
 \ \ \ 
 A morphism $f^\circledS \ (:=(f_b, f^\flat)) : Z^\circledS \migi Y^\circledS$ in $\mfS \mfc \mfh_{/X^\circledS}^\circledS$ is a monomorphism (in $\mfS \mfc \mfh_{/X^\circledS}^\circledS$)
 if and only if  the induced morphism  $f_t : Z_t \migi Y_t$  is a monomorphism in $\mfS \mfc \mfh_{/X_t}$ and $f^\flat  : f_b^*(\mcO_{Y^\circledS}) \migi \mcO_{Z^\circledS}$ is surjective.
 
  \epr
%-----------------------------------------------------------------------[begin proof]-------------------
\begin{proof}
%One may verify immediately that a morphism $f := (f_b, f^\sharp)$  in $\mr{SSch}_{X}$
% such that $f_b$ is a monomorphism and $f^\sharp$ is surjective is a monomorphism.
%Hence, we consider the inverse direction.

Let $f^\circledS  \ (:=(f_b, f^\flat)) : Z^\circledS \migi Y^\circledS$ be a monomorphism in $\mfS \mfc \mfh_{/X^\circledS}^\circledS$.
Suppose that  $f^\flat$ is {\it not} surjective, equivalently, its restriction  $f^\flat |_{f_b^*(\mcO_{Y_f})} :  f_b^* (\mcO_{Y_f}) \migi \mcO_{Z_f}$  is not surjective.
By Nakayama's lemma (and the condition that $Z^\circledS$ is noetherian),
there exists a point $s^\circledS \  (:= (s_b, s^\flat)) : \mr{Spec} (k) \migi Z^\circledS$ of $Z^\circledS$ such that $(f_b \circ s_b)^* (\mcO_{Y_f}) \migi s^*_b (\mcO_{Z_f})$ is not surjective. 
Hence, the induced morphism between $k$-vector spaces
\begin{equation}
\mr{Hom}_k (s^*_b (\mcO_{Z_f}),  k) \migi \mr{Hom}_k ((f_b \circ s_b)^* (\mcO_{Y_f}),k)
\end{equation}
is not injective.
It follows from  Lemma \ref{Pp2}
%the bijection (\ref{f1})
 that
the map  
\begin{equation} \label{R07}
 \mr{Map}_{\mfS\mfc \mfh_{/X^\circledS}^\circledS}(\mbA_k^{0\mid  k}, Z^\circledS) \migi  \mr{Map}_{\mfS \mfc \mfh_{/X^\circledS}^\circledS}(\mbA_k^{0\mid  k}, Y^\circledS)
\end{equation}
given by composing with $f^\circledS$ is not injective, and we obtain a contradiction.
Thus, $f^\flat$ must be surjective.

Next, suppose that $f_t$ is {\it not} a monomorphism in $\mfS \mfc \mfh_{/X_t}$, equivalently, 
 there exists 
an object $W$ of $\mfS \mfc \mfh_{/X_t}$ whose associated map
\begin{equation} \label{R05}
 \mr{Map}_{\mfS \mfc \mfh_{/X_t}} (W, Z_t) \migi \mr{Map}_{\mfS \mfc \mfh_{/X_t}} (W, Y_t)
\end{equation}
 is not injective.
%Here, recall that
%On the other hand, observe that 
But, since $\tau$ (cf. (\ref{DD014})) is a right adjoint functor of the functor 
 $\mfS \mfc \mfh_{/X_t} \stackrel{(\ref{R09})}{\migi} \mfS \mfc \mfh_{/X^\circledS}^\circledS$,
 %  (cf. (\ref{e968})), 
 the map  (\ref{R05})
% $\mr{Map}_{\mfS \mfc \mfh_{/X_t}} (W, f_t)$
  may be identified with the map
 \begin{equation} \label{R06}
 \mr{Map}_{\mfS \mfc \mfh_{/X^\circledS}^\circledS} (W, Z^\circledS) \migi \mr{Map}_{\mfS \mfc \mfh_{/X^\circledS}^\circledS} (W, Y^\circledS).
 \end{equation}
 This contradicts the assumption that $f^\circledS$ is a monomorphism.
 Thus, $f_t$ must be a monomorphism.
 
 The reverse direction  may be verified immediately,  and consequently, we  complete the proof of Proposition \ref{p2}.
\end{proof}
\vspace{3mm}
%-----------------------------------------------------------------------[end proposition]-------------------

%-----------------------------------------------------------------------[begin definition]------------------
%\vspace{3mm}
\bde
%[{\bf Reduced, one-pointed, and minimal}]
 \label{4}\leavevmode\\
 \ \ \  
% Let $Y^\circledS$ be an object of $\mfS \mfc \mfh_{/X^\circledS}^\circledS$.
 \vspace{-5mm}
 \begin{itemize}
 \item[(i)]
 We shall say that 
 an object $Y^\circledS$ in $\mfS \mfc \mfh_{/X^\circledS}^\circledS$ is {\bf minimal (over $X^\circledS$)} if it is nonempty (i.e., not an initial object of $\mfS \mfc \mfh_{/X^\circledS}^\circledS$) and
 any monomorphism $Z^\circledS \migi Y^\circledS$ from a nonempty object $Z^\circledS \in \mr{Ob} (\mfS \mfc \mfh_{/X^\circledS}^\circledS)$  to $Y^\circledS$ 
  is necessarily an isomorphism.
 % a  morphism  $Z^\circledS \migi Y^\circledS$ of $\mfS \mfc \mfh_{/X^\circledS}^\circledS$ is {\bf minimal (over $Y^\circledS$)} if $Z^\circledS$ is nonempty (i.e., not an initial object of $\mfS \mfc \mfh_{/X^\circledS}^\circledS$) and 
% any monomorphism $W^\circledS \migi Z^\circledS$ from a nonempty object $W^\circledS$ of  $\mfS \mfc \mfh_{/X^\circledS}^\circledS$ to $Z^\circledS$ 
%  is necessarily an isomorphism.
\item[(ii)]
We shall say that  an object  $Y^\circledS$ in $\mfS \mfc \mfh_{/X^\circledS}^\circledS$ is {\bf terminally minimal (over $X^\circledS$)} if  it is minimal over $X^\circledS$ and any minimal object $Z^\circledS$ over $X^\circledS$ with  $Y^\circledS \times_{X^\circledS} Z^\circledS \neq \emptyset$  admits a morphism $Z^\circledS \migi Y^\circledS$. 
%  morphism  $Z^\circledS \migi Y^\circledS$ of $\mfS \mfc \mfh_{/X^\circledS}^\circledS$  is {\bf terminally minimal (over $Y^\circledS$)}
%if $Z^\circledS$ is minimal over $Y^\circledS$ and any minimal object $W^\circledS$ over $Y^\circledS$ satisfying that the fiber product $Z^\circledS \times_{Y^\circledS} W^\circledS$ is nonempty admits a morphism $W^\circledS \migi Z^\circledS$ in $\mfS \mfc \mfh_{/Y^\circledS}^\circledS$. 
 \end{itemize}
\ede
\vspace{3mm}
%-----------------------------------------------------------------------[end definition]-------------------

These properties on objects in $\mfS \mfc \mfh_{/X^\circledS}^\circledS$ 
%The following assertion, which
 give a category-theoretic characterization of  spectrums of fields, as follows. 
%characterization of (terminally) minimal objects 
%may be verified immediately.
% from the definition of a(n) .
%In particular, the one-pointed objects of  $\mr{SSch}_X$ may be characterized category-theoretically.
%-----------------------------------------------------------------------[begin proposition]------------------
\vspace{3mm}
\bpr[{\bf Characterization of  spectrums of fields}] \label{p8}\leavevmode\\
 \ \ \ 
%Let $Y^\circledS$ be an object of $\mfS \mfc \mfh_{/X^\circledS}^\circledS$.
%Then, t
The following assertions (i) and (ii) are satisfied.
\begin{itemize}
\item[(i)]
An object $Y^\circledS$ of
%A morphism $Z^\circledS \migi Y^\circledS$ of 
$\mfS \mfc \mfh_{/X^\circledS}^\circledS$ is minimal
% in $\mr{SSch}_{X}$
if and only if
$Y^\circledS$ is isomorphic to $\mr{Spec}(k)$ for some field $k$.
\item[(ii)]
An object $Y^\circledS$ of
%A morphism $Z^\circledS \migi Y^\circledS$ of 
$\mfS \mfc \mfh_{/X^\circledS}^\circledS$
%A morphism $Z^\circledS \migi Y^\circledS$ of $\mfS \mfc \mfh_{/X^\circledS}^\circledS$
 is terminally minimal if and only if it  is a  point of $X_t$,  considered as an object of $\mfS \mfc \mfh_{/X^\circledS}^\circledS$ via composition with $\tau_{X^\circledS} : X_t \migi X^\circledS$.
\end{itemize}
Consequently, the objects of $\mfS \mfc \mfh_{/X^\circledS}^\circledS$ consisting of (super)schemes which are isomorphic to   $\mr{Spec} (k)$ for some field (resp., consisting of points of $X_t$)  may be reconstructed category-theoretically from the category $\mfS \mfc \mfh_{/X^\circledS}^\circledS$.
\epr
%-----------------------------------------------------------------------[begin proof]-------------------
\begin{proof}
The assertions are  formal consequences of the definitions of being  minimal and  terminally minimal.
\end{proof}
%\vspace{3mm}
%-----------------------------------------------------------------------[end proposition]-------------------

%-----------------------------------------------------------------------[begin definition]------------------
%\vspace{3mm}
%\bde
%[{\bf Reduced, one-pointed, and minimal}]
% \label{d4}\leavevmode\\
% \ \ \  
%We shall say that an object $Y$ is {\bf one-pointed}  if its underlying topological space consists of precisely one element.
%\ede
%\vspace{3mm}
%-----------------------------------------------------------------------[end definition]-------------------
\vspace{5mm}
%----------------------------------------------------------------------[begin subsection]-------------
\subsection{} \label{sub28}

Next, we shall consider the category-theoretic reconstruction of the  superschemes $\mbA_k^{0 \mid k}$ ($= \mbA^{0|1}_k$) and $\mbA_k^{\varepsilon \mid 0}$ (introduced below) in $\mfS \mfc \mfh_{/X^\circledS}^\circledS$.
After reconstructing these objects, one may use them to understand the local structure of $X^\circledS$ (cf. Proposition \ref{p19BS} described later).

%-----------------------------------------------------------------------[begin definition]------------------
\vspace{3mm}
\bde
%[{\bf Reduced, one-pointed, and minimal}]
 \label{D4}\leavevmode\\
 \ \ \  
We shall say that  an object  $Y^\circledS$ of $\mfS \mfc \mfh_{/X^\circledS}^\circledS$  is {\bf one-pointed}  if its underlying topological space consists precisely  of  one element.
\ede
\vspace{3mm}
%-----------------------------------------------------------------------[end definition]-------------------

The following proposition may be immediately verified.

%-----------------------------------------------------------------------[begin proposition]------------------
\vspace{3mm}
\bpr[{\bf Characterization of  one-pointed superschemes}] \label{p15}\leavevmode\\
 \ \ \ 
The one-pointed objects  of $\mfS \mfc \mfh_{/X^\circledS}^\circledS$ may be characterized category-theoretically as the nonempty objects $Y^\circledS$  which satisfy the following condition:
\begin{itemize}
\item[$(A)_{Y^\circledS}$:]
For any two minimal objects $Z^\circledS_1 \migi Y^\circledS$, $Z^\circledS_2 \migi Y^\circledS$ over $Y^\circledS$, the fiber product $Z^\circledS_1 \times_{Y^\circledS} Z^\circledS_2$ is nonempty.
\end{itemize}
\epr
%-----------------------------------------------------------------------[begin proof]-------------------
%\begin{proof}
%The assertion is a formal consequence of the definition of a minimal object.
%\end{proof}
\vspace{3mm}
%-----------------------------------------------------------------------[end proposition]-------------------

For any  field $k$, we shall  write
\begin{equation} \label{DD015}
\mbA_k^{\varepsilon \mid 0} := \mr{Spec} (k[\varepsilon]/\varepsilon^2).
\end{equation}

%-----------------------------------------------------------------------[begin proposition]------------------
\vspace{3mm}
\bpr[{\bf Characterization of  $ \mbA^{0 \mid 1}_k$}] \label{p10}\leavevmode\\
 \ \ \ 
%Let $k$ be a field for which 
Suppose that a morphism $\mr{Spec}(k) \migi X^\circledS$ (where  $k$ denotes  a field) is an object of $\mfS \mfc \mfh_{/X^\circledS}^\circledS$. (Hence, the category $\mfS \mfc \mfh_{/k}^\circledS$  may be characterized category-theoretically from the data $(\mfS \mfc \mfh_{/X^\circledS}^\circledS, \mr{Spec}(k))$, i.e., a pair consisting of a category and a minimal   object of  it.)
Then, the following assertions (i) and (ii) are satisfied.
\begin{itemize}
\item[(i)]
The  set  consisting of two objects
\begin{equation}
\{ \mbA^{0 \mid 1}_k,  \mbA^{\varepsilon \mid 0}_k \}
\end{equation}
of $\mfS \mfc \mfh_{/X^\circledS}^\circledS$  may be characterized 
(up to isomorphism in an evident sense)
category-theoretically as the image (via the  functor $\mfS \mfc \mfh_{/k}^\circledS \migi  \mfS \mfc \mfh_{/X^\circledS}^\circledS$ given by composing with $\mr{Spec}(k) \migi X^\circledS$) of  the set $\{ S^\circledS, T^\circledS \}$ of two one-pointed objects of  $\mfS \mfc \mfh_{/k}^\circledS$ which
% are not isomorphic and 
satisfies the following two conditions $(B)_{S^\circledS, T^\circledS}$ and $(C)_{S^\circledS, T^\circledS}$:
\vspace{2mm}
\begin{itemize}
\item[$(B)_{S^\circledS, T^\circledS}$:]
$S^\circledS$ is not isomorphic to $T^\circledS$, and 
$\mr{Spec}(k)$ is isomorphic to neither 
$S^\circledS$ nor  $T^\circledS$;
\vspace{2mm}
\item[$(C)_{S^\circledS,T^\circledS}$:]
Let $V^\circledS$ be 
a one-pointed object $V^\circledS$  of $\mfS \mfc \mfh_{/k}^\circledS$
satisfying the following two conditions: 
\begin{itemize}
\item[$\bullet$]
$V^\circledS$ is not isomorphic to $\mr{Spec}(k)$;
\item[$\bullet$]
Any  terminally minimal object over $V^\circledS$ (which is uniquely determined up to isomorphism)
 is isomorphic to the terminal object $\mr{Spec} (k)$.
\end{itemize}
Then, 
% which satisfies the condition $(D)_{V^\circledS}$ described in assertion (ii) below and which  is not isomorphic to $\mr{Spec}(k)$, 
 there exists either a monomorphism $S^\circledS \migiincl V^\circledS$ from $S^\circledS$ or a monomorphism $T^\circledS \migiincl V^\circledS$ from $T^\circledS$.
%\item[(3)]
\end{itemize}
\vspace{2mm}
%\item[(ii)]
%For $U^\circledS  \in \{ \mbA^{0 \mid 1}_k,  \mbA^{\varepsilon \mid 0}_k \}$,
%we shall denote by
% $\mcC_{U^\circledS}$ the full subcategory of $\mfS \mfc \mfh_{/k}^\circledS$  consisting of the one-pointed objects $V^\circledS$
% satisfying the following two conditions:
% \vspace{2mm}
% \begin{itemize}
% \item[$(D)_{V^\circledS}$:]
% Any  terminally minimal object over $V^\circledS$ (which is uniquely determined up to isomorphism)
% is isomorphic to the terminal object $\mr{Spec} (k)$;
% \vspace{2mm}
% \item[$(E)_{V^\circledS}$:]
% $V^\circledS$ admits a monomorphism $U^\circledS \migiincl V^\circledS$ from $U^\circledS$ and do not admit any monomorphism
%$U'^\circledS  \migiincl V^\circledS$ from $U'^\circledS$, where  $U'^\circledS$ denotes  the unique  object in    $ \{ \mbA^{0 \mid 1}_k,  \mbA^{\varepsilon \mid 0}_k \} \setminus \{ U^\circledS \}$.
% \end{itemize}
%\vspace{2mm} 
%Then, $\mcC_{\mbA^{0 \mid 1}_k}$ coincides with ${^\circ \mfS \mfc \mfh}_{/k}^\circledS$ and $\mcC_{ \mbA^{\varepsilon \mid 0}_k}$ coincides with the full subcategory of $\mfS \mfc \mfh_{/k}$ ($\subseteq \mfS \mfc \mfh_{/k}^\circledS$) consisting of  finite local $k$-schemes.
\item[(ii)]
 Let $U^\circledS$ be either  $\mbA^{0 \mid 1}_k$ or $\mbA^{\varepsilon \mid 0}_k$,  and denote by
 $U'^\circledS$ the unique  object in    $ \{ \mbA^{0 \mid 1}_k,  \mbA^{\varepsilon \mid 0}_k \} \setminus \{ U^\circledS \}$.
 Then, 
 %An  object 
 $U^\circledS$ 
 %of $\{ \mbA^{0 \mid 1}_k,  \mbA^{\varepsilon \mid 0}_k \}$
  coincides with   $\mbA^{0 \mid 1}_k$ if and only if 
 for any morphism $U'^\circledS \times_k U'^\circledS \migi U^\circledS \times_k U^\circledS$ 
 %(where  $U'^\circledS$ denotes  the unique  object in    $ \{ \mbA^{0 \mid 1}_k,  \mbA^{\varepsilon \mid 0}_k \} \setminus \{ U^\circledS \}$)
  factors through a terminally minimal morphism  over $U^\circledS \times_k U^\circledS$.
 %for  any two objects $S^\circledS$, $T^\circledS$ in $\mcC_{U^\circledS}$,
% there necessarily exists either a monomorphism from $S^\circledS$ to $T^\circledS$ or a monomorphism from $T^\circledS$ to $S^\circledS$.
 In particular, the object  $\mbA^{0 \mid 1}_k$ (resp., $\mbA_k^{\varepsilon \mid 0}$) in $\mfS \mfc \mfh_{/X^\circledS}^\circledS$ may be reconstructed  category-theoretically (up to isomorphism) from
  the minimal object $\mr{Spec} (k)$ in $\mfS \mfc \mfh_{/X^\circledS}^\circledS$.
 % object $\mr{Spec} (k) \migi X^\circledS$ in 
%  minimal morphism  over  a terminal object (i.e., $X^\circledS$) in $\mfS \mfc \mfh_{/X^\circledS}^\circledS$.
 %  the data $(\mfS \mfc \mfh_{/X^\circledS}^\circledS, \mr{Spec} (k))$, i.e., a pair consisting of a category and a minimal object over  its terminal object.
 % there does not exist any inductive system
%\begin{align}
%(\{ V^\circledS_i \}_{i \geq 0};  \{ \phi^\circledS_{i, i'} : V^\circledS_i \migi V^\circledS_{i'}\}_{i' \geq i})
%\end{align}
%  in $\mcC_{U^\circledS}$
%whose  inductive limit $\varinjlim_i V^\circledS_i$ exists 
% in $\mfS \mfc \mfh_{/k}^\circledS$ and   has at least two points (i.e., admits  a pair of two  minimal objects $(W_1 \migi \varinjlim_i V_i,  W_2 \migi \varinjlim_i V_i)$ over $\varinjlim_i V_i$ such that
%the fiber product $W_1 \times_{\varinjlim_i V_i} W_2$ is  isomorphic to an initial object of $\mr{SSch}_k$).
\end{itemize}
\epr
%-----------------------------------------------------------------------[begin proof]-------------------
\begin{proof}
Consider assertion (i).
Since  the set $\{ \mbA_k^{0 \mid 1}, \mbA_k^{\varepsilon \mid 0} \}$ is immediately verified to  satisfy both the conditions $(B)_{S^\circledS, T^\circledS}$ and $(C)_{S^\circledS, T^\circledS}$,
it suffices to prove its reverse direction.

 Note that  any one-pointed object of $\mfS \mfc \mfh_{/k}^\circledS$ is necessarily  isomorphic to the superspectrum of some (local) superalgebra over $k$.
%In the following,
For   a one-pointed object   $W^\circledS$ in $\mfS \mfc \mfh_{/k}^\circledS$,
% satisfying the condition $(D)_{W^\circledS}$ (hence, $W^\circledS$ is the superspectrum of a finite local $k$-superalgebra), then
 we shall write 
 \begin{align}
 \mr{dim}_k (W^\circledS) := \mr{dim}_k (\Gamma (W_b, \mcO_{W^\circledS})) \  (< \infty).
 \end{align}
 %, where $R_W$ denotes a  superalgebra over $k$ (which is  finite dimensional due to $(D)_{W^\circledS}$) with $W^\circledS \cong \mr{Spec}(R_W)^\circledS$.
%Moreover, if   a one-pointed object $w$ of the form $\mr{SSpec}(R)$ for some super-algebra $R$ over $k$ satisfies the equality $\mr{dim}_k R = 2$, then it is isomorphic to either $\mr{Spec}(k)^{+ k}$ or $\mr{Spec}(k[\varepsilon]/\varepsilon^2)$.

Now, let $\{ S^\circledS, T^\circledS \}$ be a set of two one-pointed objects of $\mfS \mfc \mfh_{/k}^\circledS$ which satisfies both the conditions  $(B)_{S^\circledS, T^\circledS}$ and $(C)_{S^\circledS, T^\circledS}$.
Suppose that one of the objects $S^\circledS$  in this set satisfies  the inequality  $\mr{dim}_k (S^\circledS) \geq 3$.
%as in the statement of assertion (i), and suppose that $u_1$ satisfies the inequality  $\mr{dim}_k (u_1) \geq 3$.
By Proposition ~\ref{p2}, there does not exist a monomorphism from   $S^\circledS$ to $\mbA_k^{0 \mid 1}$ since $\mr{dim}_k (\mbA_k^{0 \mid 1}) = 2$.
It follows from the condition $(C)_{S^\circledS, T^\circledS}$  that there exists a monomorphism 
from $T^\circledS$ to $\mbA^{0 \mid 1}_k$, and hence, that $\mr{dim}_k (T^\circledS) \leq 2$ (by Proposition ~\ref{p2} again).
Since $T^\circledS \ncong \mr{Spec}(k)$ and  there does not exist a monomorphism $\mbA_k^{\varepsilon \mid 0}$ from $\mbA_k^{0 \mid 1}$, 
%By the condition (b) and the fact that  there does not exist a monomorphism from $\mr{Spec}(k[\varepsilon]/\varepsilon^2)$
%to $\mr{Spec}(k)^{+ k}$, $u_2$
$T^\circledS$  must be isomorphic to $\mbA_k^{0 \mid 1}$.
One the other hand, by a similar argument where $\mbA_k^{0 \mid 1}$ is  replaced with $\mbA_k^{\varepsilon \mid 0}$, $T^\circledS$ must be isomorphic to $\mbA_k^{\varepsilon \mid 0}$, and we obtain a contradiction.
% This is a contradiction, and c
 Consequently, we have $\mr{dim}_k (S^\circledS) = \mr{dim}_k (T^\circledS) =2$.
This implies that $S^\circledS$ and $T^\circledS$ are respectively isomorphic to either  $\mbA_k^{0 \mid 1}$ or $\mbA_k^{\varepsilon \mid 0}$.
%Conversely, the set $\{ \mbA_k^{0 \mid 1}, \mbA_k^{\varepsilon \mid 0} \}$ satisfies both the conditions $(B)_{S, T}$ and $(C)_{S, T}$.
Thus, we  complete the proof of assertion (i).

%Hence, $u_2$ is isomorphic to $R_2$ for some super-algebra $R_2$ over $k$ with $\mr{dim}_k (R_2) \leq 2$, that is,  $u_2$ is (since it is not isomorphic to $\mr{Spec}(k)$ by the condition (b)) isomorphic to either $\mr{Spec}(k)^{+ k}$ or $\mr{Spec}(k[\varepsilon]/\varepsilon^2)$.
%But, since there does not exist a monomorphism from $\mr{Spec}(k[\varepsilon]/\varepsilon^2)$
%to $u_2$
% a set of two one-pointed objects of $\mr{SSch}_k$ which are not isomorphic

%The assertion follows from the observation that
%there does not exist a monomorphism between $\mr{Spec} (k[t]/t^3)$ and $\mr{Spec} (k [t, u]/ (t^2, tu, u^2))$.

Assertion (ii) follows directly from the fact that 
 \begin{align}
 \mbA_k^{\varepsilon \mid 0} \times_k \mbA_k^{\varepsilon \mid 0} \cong \mr{Spec} (k [\epsilon_1, \epsilon_2]/(\epsilon_1^2, \epsilon_1\epsilon_2, \epsilon_2^2))
 %,  \ \mbA_k^{0\mid 1} \times_k \mbA_k^{0\mid 1} \cong \mr{Spec} ({\bigwedge}^\bullet_k (k^{\oplus 2}))^\circledS
 \end{align}
  and 
  \begin{align}
   \mbA_k^{0\mid 1} \times_k \mbA_k^{0\mid 1} \cong \mr{Spec} ({\bigwedge}^\bullet_k (k^{\oplus 2}))^\circledS, \ \  (\mbA_k^{0\mid 1} \times_k \mbA_k^{0\mid 1})_t \cong \mr{Spec}(k)
   \end{align}
    (where $\bigwedge^\bullet_k (k^{\oplus 2})$ denotes the exterior algebra  over $k$  associated with $k^{\oplus 2}$, which  admits naturally a structure of superalgebra over $k$).
%Note that any object in $\mcC_{\mbA_k^{0 \mid 1}}$ is isomorphic to $\mbA_k^{0 \mid M}$ for some finite-dimensional $k$-vector space $M$ ($\neq 0$).
% Moreover, if $M_1$ and $M_2$ are finite-dimensional  $k$-vector spaces, then any morphism $\phi^\circledS : \mbA_k^{0 \mid M_1} \migi \mbA_k^{0 \mid M_2}$ of superschemes over $k$  
%  comes from 
%   the morphism  induced from a $k$-linear morphism $\phi^\sharp : M_2 \migi M_1$.
% Now, suppose that  we are given  an arbitrary inductive system  
% \begin{equation}
% \mfS := (\{ \mbA_k^{0 \mid M_i}\}_{i \geq 0}; \{ \phi_{i, i'} :  \mbA_k^{0 \mid M_i} \migi  \mbA_k^{0 \mid M_{i'}} \}_{i' \geq i})
% \end{equation}
%   in $\mcC_{\mbA_k^{0 \mid 1}}$.
% We denote by $\mfS^\sharp := (\{ M_i \}_{i \geq 0}; \{ \phi^\sharp_{i, i'} : M_{i'}\migi M_i \}_{i \geq i'})$ the projective system corresponding, via the discussion just above, to $\mfS$.
% Then, the inductive limit $\varinjlim_i \mbA_k^{0 \mid  M_i}$ of $\mfS$ (in $\mr{SSch}_k$) exists if and only if $\varprojlim_i M_i$ is finite-dimensional, and in this case, we have an isomorphism
% \begin{equation} \label{f9}
%\varinjlim_i \mbA_k^{0 \mid  M_i} \isom \mbA_k^{\varprojlim_i M_i}.
%\end{equation}
%But, $\mbA_k^{\varprojlim_i M_i}$ has only one point, and this implies  the proof of assertion (ii).
\end{proof}
%\vspace{3mm}
%-----------------------------------------------------------------------[end proposition]-------------------

\vspace{5mm}
%----------------------------------------------------------------------[begin subsection]-------------
\subsection{} \label{sub29}

Next, we consider reconstructing the schematic  structure of  $X_t$ from   $\mfS \mfc \mfh_{/X^\circledS}^\circledS$ (cf. Corollary \ref{p19} below), and consequently, a topological structure of the underlying space of $X^\circledS$ (cf. Proposition \ref{p12} below).
First, we observe that there exists, by means of Proposition \ref{p10}, the following category-theoretic criterion    for each  object  $Y^\circledS \in \mr{Ob} (\mfS \mfc \mfh_{/X^\circledS}^\circledS)$  to be a  scheme (i.e., $\mcO_{Y_f} =0$).
%, i.e., $ \mcO_{Y_f} = 0$.

%-----------------------------------------------------------------------[begin proposition]------------------
\vspace{3mm}
\bpr[{\bf Characterization of  schemes}] \label{p19BS}\leavevmode\\
 \ \ \ 
The objects   $Y^\circledS$ 
% $Y$   of $\mfS \mfc \mfh_{/X^\circledS}$
 of   $\mfS \mfc \mfh_{/X^\circledS}^\circledS$
 consisting of schemes (i.e., 
  contained in  the subcategory  $\mfS \mfc \mfh_{/X^\circledS}$)    %which are bosonic schemes
%which are objects of  $\mr{SSch}_X$ 
may be characterized category-theoretically as those objects which satisfy the following condition:
\begin{itemize}
\item[$(D)_{Y^\circledS}$:]
For any minimal object $W$ over $X^\circledS$ (hence $W \cong \mr{Spec}(k)$ for some field $k$), the map
\begin{equation}
\mr{Map}_{\mfS \mfc \mfh_{/X^\circledS}^\circledS} (W, Y^\circledS) \migi \mr{Map}_{\mfS \mfc \mfh_{/X^\circledS}^\circledS} (\mbA_W^{0 \mid 1}, Y^\circledS)
\end{equation}
 induced from the  morphism $\beta_{\mbA_W^{0 \mid 1}} : \mbA_W^{0 \mid 1} \migi W$ is bijective.
%If $s : \mr{Spec}(k) \migi Y$ (where $k$ is a field) is a minimal object over $Y$, then 
%any morphism $\mr{Spce}(k)^{+k} \migi Y$ whose restriction to $\mr{Spec}(k)$ ($\migiincl \mr{Spec}(k)^{+ k}$)  is  $s$ coincides with the composite  $\mr{Spec}(k)^{+ k}  \migi \mr{Spec}(k) \stackrel{s}{\migi} Y$, where the first arrow denotes the structure morphism of the super-scheme $\mr{Spec}(k)^{+ k}$ over $k$.
\end{itemize}
%the map $\mr{Map} (s^{+ 0 \mid 1})$
%any morphism $s^{+ 0 \mid 1} \migi Y$ factor 
In particular, the full subcategory $\mfS \mfc \mfh_{/X^\circledS}$ of $\mfS \mfc \mfh_{/X^\circledS}^\circledS$ may be reconstructed category-theoretically.
\epr
%-----------------------------------------------------------------------[begin proof]-------------------
\begin{proof}
The assertion is a formal consequence of  Nakayama' lemma and  Lemma \ref{Pp2}.
%the bijection (\ref{f1}).
\end{proof}
\vspace{3mm}
%-----------------------------------------------------------------------[end proposition]-------------------

%In particular, the truncation $\tau_b (Y)$ of $Y$ may be characterized category-theoretically, as we assert in the following.
%The following corollary follows from the 

Moreover, by  Proposition \ref{p19BS}, one may have,   for each $Y^\circledS \in \mr{Ob} (\mfS \mfc \mfh_{/X^\circledS}^\circledS)$,  a category-theoretic  reconstruction of the  schematic structure of  $Y_t$, as follows.
%and hence, of the various  open subschemes of $Y_t$, as well as the bosonic scheme $\mbA_{Y_t}^{1\mid 0}$ over $Y_t$. 

%-----------------------------------------------------------------------[begin proposition]------------------
\vspace{3mm}
\bco[{\bf Characterization of $Y_t$ for $Y^\circledS \in \mr{Ob} (\mfS \mfc \mfh_{/X^\circledS}^\circledS)$}] \label{p19}\leavevmode\\
 \ \ \ 
Let   $Y^\circledS$   be an object of $\mfS \mfc \mfh_{/X^\circledS}^\circledS$.
\begin{itemize}
\item[(i)]
The object $Y_t \in \mr{Ob} (\mfS \mfc \mfh_{/X^\circledS}^\circledS)$
%which are objects of  $\mr{SSch}_X$ 
may be characterized (up to  isomorphism) category-theoretically as
the object $Z^\circledS$ of $\mfS \mfc \mfh_{/X^\circledS}^\circledS$ which is a   scheme (i.e., satisfies the condition $(D)_{Z^\circledS}$ in Proposition \ref{p19BS}) and  satisfies the following condition:
% those super-schemes which satisfy the following property:
\vspace{0mm}
\begin{itemize}
\item[$(E)_{Z^\circledS}$:]
For any object $W$ in  $\mfS \mfc \mfh_{/Y^\circledS} \ (\subseteq \mfS \mfc \mfh_{/Y^\circledS}^\circledS)$, 
there exists uniquely a morphism $W \migi Z^\circledS$.
%If we are given  an object $W$ of $\mr{SSch}_Y$ which is a scheme,
% and a morphism $Z \migi $
%then there exists uniquely a morphism $W \migi Z$ in $\mr{SSch}_Y$.
\end{itemize}
\vspace{2mm}
%the map $\mr{Map} (s^{+ 0 \mid 1})$
%any morphism $s^{+ 0 \mid 1} \migi Y$ factor 
\item[(ii)]
The schematic structure of   $Y_t$ (i.e., a topological space together with a sheaf of rings on it), as well as  the topological structure of (the underlying space of) $Y_b$
%(in particular, those of  the set of open subschemas of $X$, as well as, of the bosonic scheme $\mbA^{1\mid 0}_{Y^\tau}$)
%Consequently, the category $\mr{Sch}_Y$ ($\cong \mr{Sch}_{Y^\tau}$)
 may be reconstructed (up to isomorphism) category-theoretically from the data $(\mfS \mfc \mfh_{/X^\circledS}^\circledS, Y^\circledS)$, i.e., a pair consisting of a category and an object of it.
\end{itemize}
\eco
%-----------------------------------------------------------------------[begin proof]-------------------
\begin{proof}
Assertion (i) follows from the functorial bijection (\ref{e968}).
% the discussion preceding Definition \ref{d3}.
Assertion (ii) follows from ~\cite{Mzk1}, Theorem A, and the fact that 
the morphism of topological spaces underlying $\gamma_Y : Y_t \migi Y_b$ is a homeomorphism.
Indeed,  we may reconstruct (un to equivalence) the category $\mfS \mfc \mfh_{/Y_t}$ ($\cong \mfS \mfc \mfh_{/Y^\circledS}$) from the data $(\mfS \mfc \mfh_{/X^\circledS}^\circledS, Y^\circledS)$ (by  Proposition \ref{p19BS} and assertion (i)).
%The assertion is a formal consequence of  Nakayama' lemma and the bijection (\ref{f1}).
\end{proof}
\vspace{3mm}
%-----------------------------------------------------------------------[end proposition]-------------------

%-----------------------------------------------------------------------[begin proposition]------------------
\vspace{3mm}
\bpr[{\bf Characterization of $X^\circledS |_U$ for an  open $U$}] \label{p12}\leavevmode\\
\ \ \ Let $Y^\circledS$ be an object of $\mfS \mfc \mfh_{/X^\circledS}^\circledS$ and  $\overline{U}$  a quasi-compact   open 
subscheme of $Y_t$. 
Denote by $U$ the  (quasi-compact) open subscheme of $Y_b$ with $\gamma_Y^{-1} (U) = \overline{U}$.
%subset of the underlying topological space of $Y$ (hence, of $Y_t$).
Then, the object $Y^\circledS |_U$  of $\mfS \mfc \mfh_{/Y^\circledS}^\circledS$ may be characterized (up to  isomorphism) category-theoretically as the object $Z^\circledS$ of $\mfS \mfc \mfh_{/Y^\circledS}^\circledS$ which satisfies the following condition:
\begin{itemize}
\item[$(F)_{Z^\circledS, \overline{U}}$:]
For any  object  $W^\circledS \stackrel{f^\circledS}{\migi} Y^\circledS$ of $\mfS \mfc \mfh_{/Y^\circledS}^\circledS$ such that the image of  $f_t : W_t \migi Y_t$   lies in $U$,  there exists uniquely a morphism $W^\circledS \migi Z^\circledS$ in $\mfS \mfc \mfh_{/Y^\circledS}^\circledS$.
% such that the composite
%\item[(ii)]
\end{itemize}
Consequently, the objects of $\mfS \mfc \mfh_{/X^\circledS}^\circledS$ consisting of 
quasi-compact open subsuperschemes of $X^\circledS$
may be characterized as the objects $V^\circledS$ such that for any $Y^\circledS \in \mr{Ob} (\mfS \mfc \mfh_{/X^\circledS}^\circledS)$, the fiber product $V^\circledS \times_{X^\circledS} Y^\circledS$ 
%is an open subsuperscheme of $Y^\circledS$
satisfies the condition $(F)_{V^\circledS \times_{X^\circledS} Y^\circledS, \overline{U}}$ for some open subscheme  $\overline{U}$ of $Y_t$.
% (resp., and moreover,  $V^\circledS \in \mr{Ob} (\mfS \mfc \mfh_{/X^\circledS}^\circledS)$). 
%quasi-compact open subsuperschemes of $Y^\circledS$ may be reconstructed category-theoretically (up to isomorphism) from the data  $(\mfS \mfc \mfh_{/X^\circledS}^\circledS, Y^\circledS)$.
\epr
%-----------------------------------------------------------------------[begin proof]-------------------
\begin{proof}
This is a formal consequence of the definition of a quasi-compact  open subsuperscheme.
\end{proof}
%\vspace{3mm}
%-----------------------------------------------------------------------[end proposition]-------------------

\vspace{5mm}
%----------------------------------------------------------------------[begin subsection]-------------
\subsection{} \label{sub36}
Next, we consider reconstructing (cf. Proposition \ref{p1811},  Lemma \ref{p1844},  and Lemma  \ref{p1833} below)
the ring object $\mbA^{1|0}_{X^\circledS}$ over $X^\circledS$ (more precisely, the objects $\mbA^{1|0}_{Y^\circledS}$  for various $Y^\circledS \in \mr{Ob} (\mfS \mfc \mfh_{/X^\circledS}^\circledS)$) corresponding to the ring structure of $\mcO_{X_b}$.
% the structure of the  scheme  $X_b$ (cf. Corollary \ref{p12open} below).

%-----------------------------------------------------------------------[begin proposition]------------------
\vspace{3mm}
\bpr
[{\bf Characterization of $\mbA_{Y^\circledS}^{1\mid 0}$ for $Y^\circledS \in \mr{Ob} (\mfS \mfc \mfh_{/X^\circledS}^\circledS)$}]
\label{p1811}\leavevmode\\
 \ \ \ 
Let $Y^\circledS$ be an object of $\mfS \mfc \mfh_{/X^\circledS}^\circledS$.
%\begin{itemize}
%\item[(i)]
Also, let  
\begin{align}
\mfz := (Z^\circledS, \sigma^{0 \circledS}, \sigma^{1 \circledS})
\end{align}
 be  a triple consisting of 
an object $Z^\circledS$ of $\mfS \mfc \mfh_{/Y^\circledS}^\circledS$ and two sections $Y^\circledS \migi Z^\circledS$ of the structure morphism $Z^\circledS \migi Y^\circledS$ of $Z^\circledS$.
Then, $\mfz$ is isomorphic to 
 $\mfa_Y := (\mbA_{Y^\circledS}^{1\mid 0}, \sigma^{[0]}_{Y^\circledS}, \sigma^{[1]}_{Y^\circledS})$  (more precisely,
 % (which means that
   there exists an isomorphism $h^\circledS : Z^\circledS \isom \mbA_{Y^\circledS}^{1\mid 0}$ over $Y^\circledS$  satisfying the equalities  $h^\circledS \circ \sigma^{0\circledS} = \sigma^{[0]}_{Y^\circledS}$ and $h^\circledS \circ \sigma^{1 \circledS} = \sigma^{[1]}_{Y^\circledS}$)  if and only if 
%The object $\mbA_Y^{1\mid 0}$ of $\mr{SSch}_Y$  may be characterized category-theoretically (up to isomorphism)
%as an object of $\mr{SSch}_Y$
% t triple   $(Z, \sigma^{[0]}, \sigma^{[1]})$ 
it  satisfies  the following three conditions $(G)_{\mfz}$-$(I)_{\mfz}$:
%\begin{itemize}
%\item[]
\begin{itemize}
\item[]
\begin{itemize}
\item[$(G)_{\mfz}$:]
The fiber product $Z^\circledS \times_{Y^\circledS} Y_t$ is isomorphic (over $Y_t$) to the scheme  $\mbA^{1 \mid 0}_{Y_t}$ (which may be reconstructed by  Corollary \ref{p19} (ii));
%$Z$ is a super-scheme over $Y$ such that $Z \times_Y \tau_b (Y)$ is isomorphic to $\tau_b (Y) \times_\mbZ \mbA_\mbZ^1$;
\vspace{2mm}
%\item[$(I)_{Z}$:]
%There exists two  sections $\sigma_i : Y \migi Z$  ($i =1, 2$) of the structure morphism of $Z$ such that the fiber product $Y \times_{\sigma_1, Z \sigma_2} Y$ is the initial object of $\mr{SSch}_X$;
%\vspace{2mm}
\item[$(H)_{\mfz}$:]
Suppose that we are given an arbitrary  commutative square diagram
\begin{equation}
\begin{CD}
W_0^\circledS @>>> Z^\circledS
\\
@VVV @VVV
\\
W_1^\circledS @>>> Y^\circledS
\end{CD}
\end{equation}
in $\mfS \mfc \mfh_{/Y^\circledS}^\circledS$ such that  $W_1^\circledS$ is  one-pointed and $W_0^\circledS$ is  terminally minimal  over both $W_1^\circledS$ and  $Z^\circledS$.
Then, there exists a morphism $W_1^\circledS \migi Z^\circledS$  over $Y^\circledS$, as well as under $W_0^\circledS$;
% (which makes sense due to the above square diagram).
%For any morphism $w : Y \migi W$ in $\mr{SSch}_X$  such that  $W$ is a scheme (i.e.,  satisfies the property $(d)_W$), there exists uniquely a morphism
%$w_b : Y_b \migi W$ satisfying the equality $w_b \circ \tau_Y = w$.
\vspace{2mm}
\item[$(I)_{\mfz}$:]
The fiber product $Y^\circledS \times_{\sigma^{0\circledS}, Z^\circledS, \sigma^{1 \circledS}} Y^\circledS$ is empty.
\end{itemize}
\vspace{2mm}
\end{itemize}
%\end{itemize}
%\end{itemize}
\epr
%-----------------------------------------------------------------------[begin proof]-------------------
\begin{proof}
%First, we consider assertion (i).
One may verify immediately that  the triple  $\mfa_Y$
%$(\mbA_{Y^\circledS}^{1\mid 0}, \sigma^{[0]}_{Y^\circledS}, \sigma^{[1]}_{Y^\circledS})$
  satisfies the three conditions $(G)_{\mfa_Y}$, $(H)_{\mfa_Y}$, and $(I)_{\mfa_Y}$. 
Hence, it suffices to prove its reverse direction.

Let $\mfz : = (Z^\circledS, \sigma^{0\circledS}, \sigma^{1\circledS})$ be a triple satisfying the required three conditions.
To begin with, we shall prove the claim that
{\it $Z^\circledS$ is, Zariski locally on $Y_b$, isomorphic to  $\mbA_{Y^\circledS}^{1\mid 0}$}.
Let $y$ be a closed point of $Y_b$ and write $Y'^\circledS :=  Y^\circledS \times_{Y_b} \mr{Spec}(\mcO_{Y_b, y})$ and $Z'^\circledS := Z^\circledS \times_{Y_b} \mr{Spec}(\mcO_{Y_b, y})$.
By the condition $(G)_{\mfz}$, the  fiber of the natural morphism  $Z'^\circledS \migi Y'^\circledS$ at $y$  is isomorphic to $\mbA^{1 \mid 0}_y$.
Let us take 
a morphism $f^\circledS: Z'^\circledS \migi \mbA_{Y'^\circledS}^{1\mid 0}$ over $Y'^\circledS$ whose restriction to the  fibers  at $y$ is an isomorphism.
(Such a morphism necessarily exists due to the universal property of the polynomial ring $\mcO_{Y'^\circledS}[t]$ with coefficients in $\mcO_{Y'^\circledS}$.)
To complete the proof of the claim, it suffices to prove that $f^\circledS$ is an isomorphism.
Let $z$ be an arbitrary point of $Z'_b$ lying over $y$.
%$Z'^\circledS \times_{Y_b} y$.
%Let us write $Z' = \mr{SSpec} (R_y)$ for some super-algebra $R_Z$ over $R_y$ and write
%$f^\dagger : R_y [t] \migi R_Z$ for the homomorphism of superalgebras over $R_y$ corresponding to
%$f$.
Write  
\begin{equation}
f^\flat_z : (\mcO_{\mbA_{Y'^\circledS}^{1\mid 0}, f_b (z)},  \mfm_{\mbA_{Y'^\circledS}^{1\mid 0}, f_b (z)}) \migi (\mcO_{Z'^\circledS, z}, \mfm_{Z'^\circledS, z})
\end{equation}
(where $\mfm_{(-)}$ denotes the maximal ideal) for   the homomorphism  of local rings defined by $f^\circledS$ and  (for each  $i \geq 1$)  write
\begin{equation}
f^{\flat, i}_z: \mcO_{\mbA_{Y'^\circledS}^{1\mid 0}, f_b (z)}/ \mfm^i_{\mbA_{Y'^\circledS}^{1\mid 0}, f_b (z)} \migi \mcO_{Z'^\circledS, z}/\mfm^i_{Z'^\circledS, z}
\end{equation}
for the  induced homomorphism.
By the definition of $f^\circledS$
% (which implies that $f^{\sharp, 1}_z$ is an isomorphism)
and   Nakayama's lemma for noncommutative rings,
all $f^{\flat}_z$ and $f^{\flat, i}_z$ ($i =1, 2, \cdots$) are  surjective.
% for any $i$.
%this morphism is surjective.
%this homomorphism is surjective after tensoring with $R_{y, b}$ (where $R_{y, b}$ denotes the )
We shall show that $f^\flat_z$ is also injective.
Suppose that $f^\flat_z$ is not injective.
%$\mr{Ker} (f^\sharp_y) \neq 0$.
% , i.e., the ideal $\mr{Ker} (F^\sharp)$ is zero.
%Let $z : W_0 \migi Z'$ be a closed point of $Z'$ lying in $Z' \times_{Y'} y$ such that $\mr{Ker} (F^\sharp)$ is not zero at $z$.
One verifies, like as the case of commutative rings, that $\bigcap_{i \geq 1}  \mfm^i_{\mbA_{Y'^\circledS}^{1\mid 0}, f_b (z)} = 0$.
Hence, there exists $i \geq 1$
for which  $f^{\flat, i}_z$ is not injective.
%Since $\varprojlim_i \mcO_{Y_b, y}/\mfm_{Y_b, y}^i$ (where $\mfm_{Y_b, y}$ denotes the maximal ideal of the local ring $\mcO_{Y_b, y}$) is faithfully flat over $\mcO_{Y_b, y}$,
%there exists a positive integer $i$ for which the  morphism $\mcO_{\mbA^{1\mid 0}_Y, f(z)} /\mfm_{Y_b, y}^i \mcO_{\mbA^{1\mid 0}_Y, f(z)} \migi \mcO_{Z, z} / \mfm_{Y_b, y}^i \mcO_{Z', z}$ induced by $f^\sharp_y$  is not  injective.
By the condition  $(H)_{\mfz}$, there exists a homomorphism
\begin{equation}
g : \mcO_{Z'^\circledS, z} \migi \mcO_{\mbA^{1\mid 0}_{Y'^\circledS}, f_b(z)} /\mfm_{Y_b, y}^i \mcO_{\mbA^{1\mid 0}_{Y'^\circledS}, f_b(z)}
\end{equation}
which makes the following diagram
 \begin{equation}
 \xymatrix{
\mcO_{\mbA^{1\mid 0}_{Y'^\circledS}, f_b(z)} /\mfm_{Y_b, y} \mcO_{\mbA^{1\mid 0}_{Y'^\circledS}, f_b(z)}   & \mcO_{Z'^\circledS, z}  \ar[l] \ar[dl]_g \\
\mcO_{\mbA^{1\mid 0}_{Y'^\circledS}, f_b(z)} /\mfm_{Y_b, y}^i \mcO_{\mbA^{1\mid 0}_{Y'^\circledS}, f_b(z)} \ar[u] & \mcO_{Y'^\circledS, y} \ar[l] \ar[u]
 }
 \end{equation}
commute, where the upper horizontal arrow denotes the composite
of the quotient $\mcO_{Z'^\circledS, z} \migisurj \mcO_{Z'^\circledS, z}/\mfm_{Z'^\circledS, z}$ and the isomorphism $(f^{\flat, 1}_z)^{-1}$.
This homomorphism $g$ factors through the quotient $\mcO_{Z'^\circledS, z} \migisurj \mcO_{Z'^\circledS, z}/\mfm^i_{Z'^\circledS, z}$.
The resulting homomorphism
% which is over $\mcO_{\mbA^{1\mid 0}_Y, f(z)}$, as well as   under $\mcO_{Z', z} / \mfm_{Y_b, y}^i \mcO_{Z', z}$.
%It induces, via tensoring with  a homomorphism
%Since it factors through
%a natural quotient $ \mcO_{Z', z}  \migisurj \mcO_{Z, z} / \mfm_{Y_b, y}^i \mcO_{Z, z}$, we obtain the resulting 
 %homomorphism
\begin{equation}
g' :  
%(\mcO_Z \otimes_{\mcO_{\mbA^{1\mid 0}}} \mcO_{\mbA^{1\mid 0}_Y, z} /\mfm_{Y_b, y}^i \mcO_{\mbA^{1\mid 0}_Y, z}  \cong) \ 
 \mcO_{Z'^\circledS, z}/\mfm^i_{Z'^\circledS, z} \migi  \mcO_{\mbA_{Y'^\circledS}^{1\mid 0}, f_b (z)}/ \mfm^i_{\mbA_{Y'^\circledS}^{1\mid 0}, f_b (z)}
%\mcO_{Z, z} / \mfm_{Y_b, y}^i \mcO_{Z, z} \migi \mcO_{\mbA^{1\mid 0}_Y, z} /\mfm_{Y_b, y}^i \mcO_{\mbA^{1\mid 0}_Y, z}
\end{equation}
becomes a split injection of $f^{\flat, i}_z$.
Thus, we have 
\begin{equation}
\mcO_{\mbA_{Y'^\circledS}^{1\mid 0}, f_b (z)}/ \mfm^i_{\mbA_{Y'^\circledS}^{1\mid 0}, f_b (z)} \cong  (\mcO_{Z'^\circledS, z}/\mfm^i_{Z'^\circledS, z}) \oplus \mr{Ker} (f^{\flat, i}_z),
\end{equation}
 which contradicts the fact that $f^{\flat, 1}_z$ is an isomorphism.
Consequently, $f^\flat_z$ is an isomorphism (for any $z$), that is to say, $f^\circledS$ is an isomorphism.
This completes the proof of the claim.
%such that $g \circ h$ coincides with identity of $\mcO_{\mbA^{1\mid 0}_Y, z} /\mfm_{Y_b, y}^i \mcO_{\mbA^{1\mid 0}_Y, z}$.
%This implies that $h$ is injective, and we obtain a contradiction.
%COnsequecnly, $\mr{Ker} () =0$, and hence, this proves the claim.

Finally, it follows immediately from the  condition $(I)_{\mfz}$ and a standard argument  that  $Z^\circledS$ is  isomorphic to $\mbA^{1\mid 0}_{Y^\circledS}$.
%Conversely, the object $\mbA_Y^{1\mid 0}$ satisfies the three properties $(H)_{\mbA_Y^{1\mid 0}}$, $(I)_{\mbA_Y^{1\mid 0}}$, $(J)_{\mbA_Y^{1\mid 0}}$.
This complies the proof of Proposition \ref{p1811}.
%Then, there exists a local ring $R''$ together with a faithfully flat $\mr{Spec}(R'') \migi \mr{Spec}(\mcCO_{Y_b, y})$ such that the closed point  of $\mr{Spec}(R'')$ is isomorphic to the scheme $W_0$ over $\mr{Spec}(\mcCO_{Y_b, y})$.
%After possibly 
%The assertion follows from the definition of $\tau_Y$.
\end{proof}
\vspace{3mm}
%-----------------------------------------------------------------------[end proposition]-------------------

Let $Y^\circledS$ be an object of $\mfS \mfc \mfh_{/X^\circledS}^\circledS$.
We shall define a functor 
\begin{equation}
(\mbG_m)_{Y^\circledS} : \mfS \mfc \mfh_{/Y^\circledS}^\circledS \migi \mfG \mfr \mfp
\end{equation}
 (where $\mfG \mfr \mfp$ denotes the category of groups)
to be the functor  which, to any object  $Z^\circledS$ of  $\mfS \mfc \mfh_{/Y^\circledS}^\circledS$, assigns the group of automorphisms of $\mbA_{Z^\circledS}^{1 \mid 0}$ over $Z^\circledS$ that are  compatible with $\sigma^{[0]}_{Z^\circledS} : Z^\circledS \migi \mbA_{Z^\circledS}^{1\mid 0}$.
It may be represented uniquely (up to a canonical isomorphism) by an object of $\mfS \mfc \mfh_{/Y^\circledS}^\circledS$, which we also  denote by $(\mbG_m)_{Y^\circledS}$ by abuse of notation.
(Indeed,  the open subsuperscheme $\mbA_{Y^\circledS}^{1 \mid 0} |_{\mbA_{Y_b}^{1\mid 0} \setminus \mr{Im} ((\sigma^{[0]}_{Y^\circledS})_b)}$ of $\mbA_{Y^\circledS}^{1\mid 0}$ represents this functor.)
%In particular, the group objects $(\mbG_m)_{Y^\circledS}$ in $\mfS \mfc \mfh_{/Y^\circledS}^\circledS$ may be reconstructed category-theoretically (up to isomorphism) from the data $(\mfS \mfc \mfh_{/X^\circledS}^\circledS, Y^\circledS)$.
Write
\begin{equation}
\mu_{Y^\circledS}^{\mbG_m} : (\mbG_m)_{Y^\circledS} \times_{Y^\circledS} (\mbG_m)_{Y^\circledS} \migi (\mbG_m)_{Y^\circledS}
\end{equation}
for the multiplication morphism of $(\mbG_m)_{Y^\circledS}$, and write
\begin{equation}
\mu^\dagger_{Y^\circledS} : (\mbG_m)_{Y^\circledS} \times_{Y^\circledS} \mbA^{1 \mid 0}_{Y^\circledS} \migi \mbA_{Y^\circledS}^{1\mid 0}
\end{equation}
for the natural action of $(\mbG_m)_{Y^\circledS}$ on $\mbA_{Y^\circledS}^{1\mid 0}$.
The morphism $\mu^\dagger_{Y^\circledS}$ induces a morphism
\begin{equation}
\nu_{Y^\circledS} \  (:= \mu_{Y^\circledS}^\dagger \circ (\mr{id}_{(\mbG_m)_{Y^\circledS}} \times \sigma^{[1]}_{Y^\circledS})) : (\mbG_m)_{Y^\circledS} \migi \mbA^{1\mid 0}_{Y^\circledS}
\end{equation}
which is an open immersion.
It follows from Proposition  \ref{p1811} that the group object $(\mbG_m)_{Y^\circledS}$   in $\mfS \mfc \mfh_{/Y^\circledS}^\circledS$  and the morphisms $\mu_{Y^\circledS}^\dagger$ and $\nu_{Y^\circledS}$ in $\mfS \mfc \mfh_{/Y^\circledS}^\circledS$  may be reconstructed (up to isomorphism) category-theoretically from the data $(\mfS \mfc \mfh_{/X^\circledS}^\circledS, Y^\circledS)$.
The following two lemmas will be used in the proof of Corollary \ref{p12open} below.

%-----------------------------------------------------------------------[begin proposition]------------------
\vspace{3mm}
\ble
%[{\bf Characterization of $\mbA_Y^{1\mid 0}$ for $Y \in \mr{Ob} (\mr{SSch}_X)$}]
\label{p1844}\leavevmode\\
 \ \ \ %Denote by $\nu_Y : (\mbG_m)_Y \migi \mbA^{1\mid 0}_Y$ the open immersion.
Denote by 
\begin{equation}
 \mu^{1\mid 0}_{Y^\circledS} : \mbA_{Y^\circledS}^{1 \mid 0} \times_{Y^\circledS} \mbA_{Y^\circledS}^{1 \mid 0} \migi \mbA_{Y^\circledS}^{1 \mid 0}
\end{equation}
the morphism corresponding to the  multiplication of
$\mcO_{Y_b}$ (via the functorial bijection (\ref{e5})).
% $\Gamma (Y_b, \mcO_{Y_b})$ ($\cong \mr{Map}_{\mfS \mfc \mfh_{/Y^\circledS}^\circledS}(Y^\circledS, \mbA_{Y^\circledS}^{1\mid 0})$ via the bijection (\ref{e5})).
%Now, let us fix an isomorphism $(\mbG_m)_Y \isom (\mbG_m)_Y^\dagger$ (cf. assertion (ii) above).
%In particular, it gives, by definition, an action
%\begin{equation} \label{e49}
%\mu_Y^\dagger  : (\mbG_m)_Y  \times_Y \mbA_Y^{1 \mid 0} \migi \mbA_Y^{1 \mid 0}
%\end{equation}
%of $(\mbG_m)_Y$ on $ \mbA_Y^{1 \mid 0}$.
Then, a   morphism  $\mu^\circledS :  \mbA_{Y^\circledS}^{1 \mid 0} \times_{Y^\circledS} \mbA_{Y^\circledS}^{1 \mid 0} \migi \mbA_{Y^\circledS}^{1 \mid 0}$ in  $\mfS \mfc \mfh_{/Y^\circledS}^\circledS$ coincides with $\mu_{Y^\circledS}^{1\mid 0}$
%The  may be characterized category-theoretically
% (up to an automorphism of the triple $(\mbA_Y^{1 \mid 0}, \sigma^{[0]}_Y, \sigma^{[1]}_Y)$)
% the triple $(\nu_Y,  \alpha_Y, \mu_Y)$
  if and only if 
 it satisfies the following condition:
%\vspace{2mm}
 \begin{itemize}
 \vspace{2mm}
 \item[$(J)_{\mu}$:]
%It holds the equality  
%\begin{equation}
%\nu = \mu^\dagger_Y \circ (\mr{id}_{(\mbG_m)_Y} \times \sigma^{[1]}_Y)
%\end{equation}
%of morphisms $(\mbG_m)_Y \migi \mbA^{1\mid 0}_Y$;
% \vspace{2mm}
%  \item[$(M)_{\nu, \alpha, \mu}$:]
the equality 
\begin{equation}
\mu^\circledS  \circ  (\nu_{Y^\circledS} \times \nu_{Y^\circledS})  = \nu_{Y^\circledS}  \circ \mu_{Y^\circledS}^{\mbG_m}
\end{equation}
of morphisms
 $(\mbG_m)_{Y^\circledS} \times_{Y^\circledS} (\mbG_m)_{Y^\circledS} \migi \mbA_{Y^\circledS}^{1 \mid 0}$ holds;
 \vspace{2mm}
 \end{itemize}
Consequently, the morphism $\mu^{1\mid 0}_{Y^\circledS}$ in $\mfS \mfc \mfh_{/Y^\circledS}^\circledS$ may be reconstructed category-theoretically (up to isomorphism) from the data $(\mfS \mfc \mfh_{/X^\circledS}^\circledS, Y^\circledS)$.
\ele
%-----------------------------------------------------------------------[begin proof]-------------------
\begin{proof}
Since the equality $\mu_{Y^\circledS}^{1\mid 0}  \circ  (\nu_{Y^\circledS} \times \nu_{Y^\circledS})  = \nu_{Y^\circledS}  \circ \mu_{Y^\circledS}^{\mbG_m}$ holds,
the assertion follows directly  from the fact that 
$\nu_{Y^\circledS} \times \nu_{Y^\circledS}$ is an epimorphism in $\mfS \mfc \mfh_{/Y^\circledS}^\circledS$.
\end{proof}
\vspace{3mm}
%-----------------------------------------------------------------------[end proposition]-------------------

%-----------------------------------------------------------------------[begin proposition]------------------
\vspace{3mm}
\ble
%[{\bf Characterization of $\mbA_Y^{1\mid 0}$ for $Y \in \mr{Ob} (\mr{SSch}_X)$}]
\label{p1833}\leavevmode\\
 \ \ \ %Denote by $\nu_Y : (\mbG_m)_Y \migi \mbA^{1\mid 0}_Y$ the open immersion.
Denote by 
\begin{equation}
\alpha^{1\mid 0}_{Y^\circledS} : \mbA_{Y^\circledS}^{1 \mid 0} \times_{Y^\circledS} \mbA_{Y^\circledS}^{1 \mid 0} \migi \mbA_{Y^\circledS}^{1 \mid 0} 
 %(\text{resp.},  \ \mu^{1\mid 0}_Y : \mbA_Y^{1 \mid 0} \times_Y \mbA_Y^{1 \mid 0} \migi \mbA_Y^{1 \mid 0})
\end{equation}
the morphism corresponding to the addition  of $\mcO_{Y_b}$  (via the functorial bijection (\ref{e5})).
%$\Gamma (Y_b, \mcO_{Y_b})$.
%Now, let us fix an isomorphism $(\mbG_m)_Y \isom (\mbG_m)_Y^\dagger$ (cf. assertion (ii) above).
%In particular, it gives, by definition, an action
%\begin{equation} \label{e49}
%\mu_Y^\dagger  : (\mbG_m)_Y  \times_Y \mbA_Y^{1 \mid 0} \migi \mbA_Y^{1 \mid 0}
%\end{equation}
%of $(\mbG_m)_Y$ on $ \mbA_Y^{1 \mid 0}$.
Then, a   morphism  $\alpha^\circledS  :  \mbA_{Y^\circledS}^{1 \mid 0} \times_{Y^\circledS} \mbA_{Y^\circledS}^{1 \mid 0} \migi \mbA_{Y^\circledS}^{1 \mid 0}$  in $\mfS \mfc \mfh_{/Y^\circledS}^\circledS$ coincides with
%The  may be characterized category-theoretically
% (up to an automorphism of the triple $(\mbA_Y^{1 \mid 0}, \sigma^{[0]}_Y, \sigma^{[1]}_Y)$)
$\alpha_{Y^\circledS}^{1\mid 0}$ if and only if 
 it satisfies the following two conditions $(K)_{\alpha^\circledS}$ and $(L)_{\alpha^\circledS}$:
%\vspace{2mm}
 \begin{itemize}
 \vspace{2mm}
%It holds the equality  
%\begin{equation}
%\nu = \mu^\dagger_Y \circ (\mr{id}_{(\mbG_m)_Y} \times \sigma^{[1]}_Y)
%\end{equation}
%of morphisms $(\mbG_m)_Y \migi \mbA^{1\mid 0}_Y$;
% \vspace{2mm}
%  \item[$(M)_{\nu, \alpha, \mu}$:]
 \item[$(K)_{\alpha^\circledS}$:]
  The square diagram
  \begin{equation}
  \begin{CD}
  (\mbG_m)_{Y^\circledS}  \times_{Y^\circledS}  \mbA_{Y^\circledS}^{1 \mid 0} \times_{Y^\circledS}  \mbA_{Y^\circledS}^{1 \mid 0} @> \mr{id}_{  (\mbG_m)_{Y^\circledS}} \times \alpha^\circledS >>  (\mbG_m)_{Y^\circledS}  \times_{Y^\circledS} \mbA_{Y^\circledS}^{1 \mid 0}
  \\
  @V (\mu^\dagger_{Y^\circledS} \times \mu^\dagger_{Y^\circledS}) \circ \lambda^\circledS VV @VV \mu^\dagger_{Y^\circledS} V
  \\
   \mbA_{Y^\circledS}^{1 \mid 0} \times_{Y^\circledS}  \mbA_{Y^\circledS}^{1 \mid 0} @> \alpha^\circledS  >>  \mbA_{Y^\circledS}^{1 \mid 0}
  \end{CD}
  \end{equation}
  is commutative, where $\lambda^\circledS$ denotes the morphism
  \begin{align}
    (\mbG_m)_{Y^\circledS}  \times_{Y^\circledS}  \mbA_{Y^\circledS}^{1 \mid 0} \times_{Y^\circledS}  \mbA_{Y^\circledS}^{1 \mid 0} &\migi   (\mbG_m)_{Y^\circledS}  \times_{Y^\circledS}    \mbA_{Y^\circledS}^{1 \mid 0} \times_{Y^\circledS} (\mbG_m)_{Y^\circledS}  \times_{Y^\circledS}   \mbA_{Y^\circledS}^{1 \mid 0} \\
    (g, a_1, a_2) \hspace{10mm} & \mapsto  \hspace{20mm} (g, a_1, g, a_2)  \notag
  \end{align}
 over $Y^\circledS$.
  \vspace{2mm}
 \item[$(L)_{\alpha^\circledS}$:]
We have the equalities 
\begin{equation}
\alpha^\circledS  \circ (\sigma^{[0]}_{Y^\circledS} \times \mr{id}_{\mbA^{1\mid 0}_{Y^\circledS}}) = \alpha^\circledS  \circ ( \mr{id}_{\mbA^{1\mid 0}_{Y^\circledS}} \times \sigma^{[0]}_{Y^\circledS}) = \mr{id}_{\mbA^{1\mid 0}_{Y^\circledS}}.
\end{equation}
of endomorphisms  of $\mbA_{Y^\circledS}^{1\mid 0}$.
 \end{itemize}
Consequently, the morphism $\alpha^{1\mid 0}_{Y^\circledS}$ in $\mfS \mfc \mfh_{/Y^\circledS}^\circledS$ may be reconstructed category-theoretically (up to isomorphism) from the data $(\mfS \mfc \mfh_{/X^\circledS}^\circledS, Y^\circledS)$.
\ele
%-----------------------------------------------------------------------[begin proof]-------------------
\begin{proof}
Let 
$\alpha^\circledS$ be a morphism satisfying the   conditions $(K)_{\alpha^\circledS}$ and  $(L)_{\alpha^\circledS}$.
We   write
$\alpha^\flat : \mcO_{Y^\circledS} [t] \migi \mcO_{Y^\circledS}[t] \otimes_{\mcO_{Y^\circledS}} \mcO_{Y^\circledS} [t]$
for the  homomorphism of superalgebra over $\mcO_{Y^\circledS}$ corresponding to $\alpha^\circledS$.
The condition  $(L)_{\alpha^\circledS}$ implies that
$\alpha^\flat$ is given by $t \mapsto a \cdot t \otimes 1 + b \cdot  1 \otimes t$ for some $a$, $b \in \Gamma (Y_b, \mcO_{Y_b})$.
But, the equalities in $(L)_{\alpha^\circledS}$ imply that $a = b = 1$, that is to say, $\alpha^\circledS = \alpha_{Y^\circledS}^{1|0}$.
Thus, we complete the proof of Lemma \ref{p1833}.
\end{proof}
%\vspace{3mm}
%-----------------------------------------------------------------------[end proposition]-------------------

\vspace{5mm}
%----------------------------------------------------------------------[begin subsection]-------------
\subsection{} \label{su1b36}
By combining the results in  \S\,\ref{sub29} and  \S\,\ref{sub36},
one may reconstruct  category-theoretically the schematic structure of $X_b$ as follows.

% the various pairs  $(Y |_{U}, \Gamma (U,  \mcO_{Y_b}))$ (i.e., pairs consisting of an object of $\mr{SSch}_Y$ and a commutative ring)  which are .
%Consequently, we have proved the following assertion.

%-----------------------------------------------------------------------[begin proposition]------------------
\vspace{3mm}
\bco[{\bf Characterization of $Y_b$ for  $Y^\circledS \in \mr{Ob}(\mfS \mfc \mfh_{/X^\circledS}^\circledS)$}] \label{p12open}\leavevmode\\
\ \ \ 
Let $Y^\circledS$ be an object of $\mfS \mfc \mfh_{/X^\circledS}^\circledS$.
Then, the  schematic structure of   $Y_b$ (i.e., a topological space together with a sheaf of rings on it) may be reconstructed category-theoretically (up to isomorphism) from the data $(\mfS \mfc \mfh_{/X^\circledS}^\circledS, Y^\circledS)$.
% (i.e., a category and an object of it).
Moreover, this reconstruction is functorial (in a natural sense) in $Y^\circledS \in \mr{Ob}(\mfS \mfc \mfh_{/X^\circledS}^\circledS)$; strictly speaking, if we are given a morphism $f^\circledS : Z^\circledS \migi Y^\circledS$  in $\mfS \mfc \mfh_{/X^\circledS}^\circledS$, then 
(the two schemes $Y_b$,  $Z_b$ and) its  underlying morphism $f_b : Z_b \migi Y_b$
 may be reconstructed category-theoretically.
 \eco
%-----------------------------------------------------------------------[begin proof]-------------------
\begin{proof}
By Corollary \ref{p19} and Proposition  \ref{p12}, one may reconstruct  (up to equivalence) category-theoretically the topological structure of $X_b$ and  the full subcategory  of $\mfS \mfc \mfh_{/X^\circledS}^\circledS$ whose  objects  are  
\begin{align}
\{ X^\circledS |_U \in \mr{Ob} (\mfS \mfc \mfh_{/X^\circledS}^\circledS)  \ | \ \text{$U$ is a quasi-compact  open subscheme of $X_b$} \}.
\end{align}
Moreover,  it follows from Proposition \ref{p1811}, Lemma \ref{p1844},  and Lemma \ref{p1833} that  one may reconstruct ring objects $\mbA^{1|0}_{X^\circledS |_U} \in \mr{Ob} (\mfS \mfc \mfh_{X^\circledS}^\circledS)$ (for each quasi-compact open $U$ in $X_b$) over $X^\circledS |_U$ corresponding to $\mcO_{U}$.
By considering the set of various  sections $X^\circledS |_U \migi \mbA^{1|0}_{X^\circledS |_U}$,
we obtain  the ring structure of  $\Gamma (U, \mcO_{X_b})$  that is  compatible with  
restriction to   open subschemes of $U$.
Consequently, 
%open immersions between open subschemes of $X_b$, and hence, 
the schematic structure of $X_b$ may be reconstructed, as desired.
The latter assertion follows from this reconstructing  procedure.
%By   Proposition 2.9, Proposition  \ref{p1811}, Lemma \ref{p1844}, and  Lemma \ref{p1833} proved above,
%the three lemmas obtained above, 
%one may reconstruct
%,  for each open subset $U$ of  $Y$,
% the morphisms  $Y |_U \migi Y$ (i.e., the open immersion determined by $U$) in $\mr{SSch}_Y$  and a commutative ring  $ \Gamma (U,  \mcO_{Y_b})$ that are compatible with various open immersions $U' \migiincl U$.
\end{proof}
%\vspace{3mm}
%-----------------------------------------------------------------------[end proposition]-------------------

\vspace{5mm}
%----------------------------------------------------------------------[begin subsection]-------------
\subsection{} \label{sub31}

In this subsection, we consider reconstructing the various $\mbA^{0|1}$-twists associated with  fermionic twists of $X^\circledS$, together with the multiplication and addition maps.
Consequently, one may reconstruct (cf. Corollary \ref{p218}) the schematic structure of  superschemes $Z^\circledS$ with $Z^\circledS \stackrel{f}{\sim} X^\circledS$.

Let us fix   an object $Y^\circledS$ of $\mfS \mfc \mfh_{/X^\circledS}^\circledS$.

%-----------------------------------------------------------------------[begin proposition]------------------
%\vspace{3mm}
\bpr[{\bf Characterization of  $\mbA^{0 \mid 1}$-twists}] \label{p188}\leavevmode\\
 \ \ \
%Suppose that an object  $Z$ of $\mr{SSch}_Y$ satisfies the following two properties:
Let $(Z^\circledS, \sigma^\circledS)$  be a pair consisting of an object  $Z^\circledS$ of  $\mfS \mfc \mfh_{/Y^\circledS}^\circledS$  (i.e., a morphism $f^\circledS : Z^\circledS \migi Y^\circledS$) and a morphism $\sigma^\circledS : Y^\circledS \migi Z^\circledS$ in $\mfS \mfc \mfh_{/Y^\circledS}^\circledS$ (i.e., a section $\sigma^\circledS$ of $f^\circledS$).
Then, the pair $(Z^\circledS,  \sigma^\circledS)$ forms 
an  $\mbA^{0 \mid 1}$-twists  over  $Y^\circledS$ if and only if it 
% may be characterized category-theoretically (up to isomorphism) as the pairs  (i.e., a section of the structure morphism $Z \migi Y$ of $Z$) 
 satisfies the following three conditions $(M)_{Z^\circledS, \sigma^\circledS}$-$(O)_{Z^\circledS, \sigma^\circledS}$:
\vspace{2mm}
\begin{itemize}
\item[]
\begin{itemize}
\item[$(M)_{Z^\circledS, \sigma^\circledS}$:]
The underlying morphism $f_b : Z_b \migi Y_b$  of schemes (which may be  reconstructed category-theoretically from the data $(\mfS \mfc \mfh_{/X^\circledS}^\circledS, f^\circledS)$) is finite (cf. Corollary \ref{p12open} for the category-theoretic characterization of this condition);
%\begin{itemize}
\item[$(N)_{Z^\circledS, \sigma^\circledS}$:]
%There exists a section  $\sigma_Z : Y \migi Z$ of the structure morphism $Y \migi Z$;
%\item[$(L)_Z$:]
For each minimal object $W^\circledS$ over $Y^\circledS$, the fiber product $Z^\circledS \times_{Y^\circledS} W^\circledS$ is  isomorphic to $\mbA_{W^\circledS}^{0\mid 1}$ (which may be  reconstructed category-theoretically from the data $(\mfS \mfc \mfh_{/Y^\circledS}^\circledS, W^\circledS)$ by Proposition \ref{p10});
%\vspace{2mm}
\item[$(O)_{Z^\circledS, \sigma^\circledS}$:]
Let $Y'^\circledS$ be an open subsuperscheme of $Y^\circledS$ (i.e., an  object $Y'^\circledS$ of $\mfS \mfc \mfh_{/Y^\circledS}^\circledS$  satisfying the condition $(F)_{Y'^\circledS, \overline{U}}$ in Proposition \ref{p12} for some open subscheme $\overline{U}$ of $Y_t$). 
Also, let $(Z'^\circledS, \sigma'^\circledS)$ be a pair, where  $Z'^\circledS$ denotes 
an object in $\mfS \mfc \mfh_{/Y'^\circledS}^\circledS$  and 
$\sigma'^\circledS$ denotes a morphism $Y'^\circledS \migi Z'^\circledS$  in 
$\mfS \mfc \mfh_{/Y'^\circledS}^\circledS$, satisfying the conditions 
$(M)_{Z'^\circledS, \sigma'^\circledS}$ and $(N)_{Z'^\circledS, \sigma'^\circledS}$.
Then,
%For each open subsuperscheme $Y^\circledS |_U$ (i.e., an  object $Z^\circledS$ of $\mfS \mfc \mfh_{/Y^\circledS}^\circledS$  satisfying the condition $(F)_{Z^\circledS}$) and 
%each  pair $(Z'^\circledS, \sigma'^\circledS)$ satisfying the condition  $(O)_{Z', \sigma'}^{Y |_U}$,
there exists an open subsuperscheme  $Y''^\circledS$  of $Y'^\circledS$ and a monomorphism
$h^\circledS: Z'^\circledS \times_{Y'^\circledS} Y''^\circledS \migi Z^\circledS \times_{Y^\circledS} Y''^\circledS$ in  $\mfS \mfc \mfh_{/ Y''^\circledS}^\circledS$ satisfying the equality of morphisms
\begin{align}
h^\circledS \circ (\sigma'^\circledS \times \mr{id}_{Y''^\circledS}) = \sigma^\circledS \times \mr{id}_{Y''^\circledS} : Y''^\circledS \migi Z^\circledS \times_{Y^\circledS} Y''^\circledS.
\end{align}
\end{itemize}
\end{itemize}
%\vspace{1mm}
Consequently, 
the set of objects in $\mfS \mfc \mfh_{/Y^\circledS}^\circledS$ which are isomorphic to $\mbA^{0|1}$-twists over $Y^\circledS$
%the morphism $\mu^{1\mid 0}_{Y^\circledS}$ in $\mfS \mfc \mfh_{/Y^\circledS}^\circledS$ 
may be reconstructed category-theoretically (up to isomorphism) from the data $(\mfS \mfc \mfh_{/X^\circledS}^\circledS, Y^\circledS)$.
 \epr
%-----------------------------------------------------------------------[begin proof]-------------------
\begin{proof}
 %the section $\sigma$,
Let $(Z^\circledS, \sigma^\circledS)$ be a pair satisfying the required three  conditions. 
%$(N)_{Z^\circledS, \sigma^\circledS}$ and $(O)_{Z^\circledS, \sigma^\circledS}$.
By the existence of  a section $\sigma^\circledS$ and the condition $(N)_{Z^\circledS, \sigma^\circledS}$,
the underlying continuous map of 
$f^\circledS$ 
%the structure morphism $Z^\circledS \migi Y^\circledS$ of $Z^\circledS$
  is a homeomorphism (hence, we consider $\mcO_{Z^\circledS}$ as a sheaf on the underlying topological space of $Y^\circledS$). 
 The conditions $(M)_{Z^\circledS, \sigma^\circledS}$  implies that $\mcO_{Z^\circledS}$ is a finite $\mcO_{Y^\circledS}$-module.
It follows from the condition  $(N)_{Z^\circledS, \sigma^\circledS}$ and 
Nakayama's lemma that
one may find,  locally on $Y_b$,   an isomorphism 
$\mcO_{Z^\circledS} \isom \mcO_{Y^\circledS} \oplus (\mcO_{Y^\circledS}/\mcJ)$ of  $\mcO_{Y^\circledS}$-superalgebras, where the multiplication of  the right-hand side is given by
$(a, \overline{b}) \cdot (c, \overline{d}) = (a c, a \overline{d} + c \overline{b})$.  
%The section $\sigma^\circledS$ gives 
 %$\mcO_Z$ admits, by means of the section $\sigma$, 
% a decomposition $\mcO_{Z^\circledS} \cong \mcO_{Y^\circledS} \oplus \mcA$ of the $\mcO_{Y^\circledS}$--module $\mcO_{Z^\circledS}$, where $\mcA$ denotes  some $\mcO_Y$-module.
%It follows from  and Nakayama's lemma that $\mcA$ is locally of the form $\mcO_Y /\mcI$ for some ideal sheaf $\mcI$ of $\mcO_Y$.
Moreover,  the universal property described in $(O)_{Z^\circledS, \sigma^\circledS}$ implies that $\mcI =0$.
Consequently, $(Z^\circledS, \sigma^\circledS)$ forms an $\mbA^{0 \mid 1}$-twist over $Y^\circledS$.
Since the reverse direction of this assertion  may be verified immediately, we complete the proof of Proposition \ref{p188}.
%the condition $(O_{Z, \sigma}$,
%$Z$ is, Zariski locally on $Y_b$, isomorphic to 
%The assertion follows from the definition of $\tau_Y$.
%The latter assertion follows directly  from the former assertion and Proposition \ref{p18} (since $\mbA_{Z/Y}^{1\mid 1} \cong \mbA_{Z/Y}^{0\mid 1} \times_Y \mbA^{1\mid 0}_Y$).
\end{proof}
\vspace{3mm}
%-----------------------------------------------------------------------[end proposition]-------------------

Next, let us fix an $\mbA^{0 \mid 1}$-twist $(Z^\circledS, \sigma_{Z^\circledS/Y^\circledS})$ over $Y^\circledS$.

%-----------------------------------------------------------------------[begin proposition]------------------
\vspace{3mm}
\ble
%[{\bf Characterization of  $\mbA^{0 \mid 1}$-twists}] 
\label{p198}\leavevmode\\
 \ \ \
%Let $Y^\circledS$ be an object of $\mfS \mfc \mfh_{/X^\circledS}^\circledS$ and  be 
We shall write 
\begin{equation}
\mr{Aut} (Z^\circledS, \sigma_{Z^\circledS/Y^\circledS}) : \mfS \mfc \mfh_{/Y^\circledS}^\circledS \migi \mfG \mfr \mfp
\end{equation}
 for  the functor which, to any  $W^\circledS \in \mr{Ob}(\mfS \mfc \mfh_{/Y^\circledS}^\circledS)$, assigns the group of autormophisms of the $\mbA^{0 \mid 1}$-twist $(Z^\circledS \times_{Y^\circledS} W^\circledS, \sigma_{Z^\circledS/Y^\circledS} \times \mr{id}_{W^\circledS})$ over $W^\circledS$.
Consider 
the   isomorphism 
\begin{align} \label{e987}
\eta_{Z^\circledS} : (\mbG_m)_{Y^\circledS} \isom \mr{Aut} (Z^\circledS, \sigma_{Z^\circledS/Y^\circledS})
%``\mcO_Y[t]"
\end{align}
which, to any automorphism in $(\mbG_m)_{Y^\circledS}(W^\circledS)$ (where  $W^\circledS \in \mr{Ob} (\mfS \mfc \mfh_{/Y^\circledS}^\circledS)$) corresponding to the automorphism of $\mcO_{W^\circledS} [t]$ determined by $t \mapsto g \cdot t$ (where $g \in \Gamma (W_b, \mcO_{W_b}^\times)$),
assigns the automorphism of $(Z^\circledS \times_{Y^\circledS} W^\circledS, \sigma_{Z^\circledS/Y^\circledS} \times \mr{id}_{W^\circledS})$ corresponding to the automorphism of $\mcO_{Z^\circledS \times_{Y^\circledS} W^\circledS}$ (which is locally isomorphic to $\mcO_{W^\circledS} [\psi]$) determined   by $\psi \mapsto g \cdot \psi$. 
 Then, an  isomorphism $\eta^\circledS : (\mbG_m)_{Y^\circledS} \isom \mr{Aut} (Z^\circledS, \sigma_{Z^\circledS/Y^\circledS})$ coincides with $\eta_{Z^\circledS}$ if and only  if it satisfies the following condition:
\vspace{2mm}
 \begin{itemize}
 \item[$(P)_{\eta^\circledS}$:]
 Let
 $W^\circledS$ be an object of  $\mfS \mfc \mfh_{/Y^\circledS}^\circledS$ and $h^\circledS \in (\mbG_m)_{Y^\circledS} (W^\circledS)$ whose  induced automorphism   of $\mcO_{W_b}[t]$ ($=\mcO_{(\mbA_W^{1\mid 0})_b}$) is  given by $t \mapsto g \cdot t$ for some 
% as multiplication by $g$ for some
  $g \in \Gamma (\mcO_{W_b}, \mcO_{W_b}^\times)$.
 (Such a pair $(W^\circledS, h^\circledS)$ may be characterized category-theoretically thanks to Corollary \ref{p12open}.)
 Here, note that
 the section 
 \begin{align}
 (\sigma_{Z^\circledS/Y^\circledS} |_{W^\circledS}, \sigma_{Z^\circledS/Y^\circledS} |_{W^\circledS}) : W^\circledS \migi Z_W^\circledS \times_{W^\circledS} Z^\circledS_W
 \end{align}
(where $Z_W^\circledS := Z^\circledS \times_{Y^\circledS} W^\circledS$)
 determines
 a decomposition $\mcO_{(Z_W^\circledS \times_{W^\circledS} Z^\circledS_W)_b} \cong \mcO_{W_b} \oplus \mcO_{W_b} \epsilon$, where the multiplication of the right-hand side is given by 
 $(a, b\epsilon) \cdot (c, d\epsilon) = (ac, (bc + ad) \epsilon)$.
 Then, the automorphism $\eta^\circledS (h^\circledS) \times \eta^\circledS (h^\circledS)$ of $Z_W^\circledS \times_{W^\circledS} Z^\circledS_W$ induces  the automorphism of $\mcO_{W_b} \oplus \mcO_{W_b} \epsilon$ given by assigning $(a, b \epsilon) \mapsto (a, g^2 b \epsilon)$.
 % Also, let $T^\circledS$ be  a   minimal object  over $W^\circledS$ with $T^\circledS \cong \mr{Spec}(k)$ for some field $k$.
% (The field $k$ and the decomposition $\mcO_{(\mbA^{0 \mid 1}_{T^\circledS} \times_{T^\circledS} \mbA^{0 \mid 1}_{T^\circledS})_b} = k \oplus k \psi \psi'$ corresponding to  the section $\tau_{\mbA^{0 \mid 1}_{T^\circledS}} \times \tau_{\mbA^{0 \mid 1}_{T^\circledS}} : T^\circledS \migi \mbA^{0 \mid 1}_{T^\circledS} \times_{T^\circledS} \mbA^{0 \mid 1}_{T^\circledS}$ may be reconstructed (up to isomorphism) category-theoretically from the data $(\mfS \mfc \mfh_{/X^\circledS}^\circledS, T^\circledS)$ by Corollary \ref{p19} and Corollary \ref{p12open}.)
 % Then,  the automorphism of $\mcO_{(\mbA^{0 \mid 1}_{T^\circledS} \times_{T^\circledS} \mbA^{0 \mid 1}_{T^\circledS})_b}$ ($= k \oplus k \psi \psi'$) corresponding to the automorphism $(h^\circledS |_{\mbA^{0 \mid 1}_T}\times h |_{\mbA^{0 \mid 1}_T})_b$
% of 
 %  $(\mbA^{0 \mid 1}_T \times_T \mbA^{0 \mid 1}_T)_b$  is given by $\psi^2 \mapsto g^2 \cdot \psi^2$.
 \end{itemize}
Consequently, 
the morphism 
$\eta_{Z^\circledS}$ in $\mfS \mfc \mfh_{/Y^\circledS}^\circledS$ 
%which are isomorphic to $\mbA^{0|1}$-twists over $Y^\circledS$
%the morphism $\mu^{1\mid 0}_{Y^\circledS}$ in $\mfS \mfc \mfh_{/Y^\circledS}^\circledS$ 
may be reconstructed category-theoretically (up to isomorphism) from  $(\mfS \mfc \mfh_{/X^\circledS}^\circledS, Y^\circledS, (Z^\circledS, \sigma_{Z^\circledS/Y^\circledS}))$, i.e., a collection of data consisting of a category $\mfS \mfc \mfh_{/X^\circledS}^\circledS$, an object $Y^\circledS$ of it, and a pair $(Z^\circledS, \sigma_{Z^\circledS/Y^\circledS})$ satisfying the conditions described in Proposition \ref{p188}.
 \ele
%-----------------------------------------------------------------------[begin proof]-------------------
\begin{proof}
The assertion follows from the various definitions involved.
\end{proof}
\vspace{3mm}
%-----------------------------------------------------------------------[end proposition]-------------------

%Let $Y^\circledS$ be an object of $\mfS \mfc \mfh_{/X^\circledS}^\circledS$ and $(Z^\circledS, \sigma_{Z^\circledS/Y^\circledS})$ an $\mbA^{0 \mid 1}$-twist over $Y^\circledS$. 
We shall write 
\begin{equation}
%\mu^\dagger_{Y} \times \eta_{(Z/Y, \sigma_{Z/Y})}
\mu\eta^\dagger_{Z^\circledS} : (\mbG_{m})_{Y^\circledS}  \times_{Y^\circledS}  \mbA^{1\mid 0}_{Z^\circledS}  \migi  \mbA^{1\mid 0}_{Z^\circledS} 
\end{equation}
for the action of $(\mbG_{m})_{Y^\circledS}$ on $\mbA^{1\mid 0}_{Z^\circledS}$ ($\cong \mbA^{1\mid 0}_{Y^\circledS} \times_{Y^\circledS} Z^\circledS$) defined by
\begin{align}
(\mbG_{m})_{Y^\circledS}  \times_{Y^\circledS}  \mbA_{Y^\circledS}^{1\mid 0} \times_{Y^\circledS} Z^\circledS &\migi  \hspace{10mm} \mbA_{Y^\circledS}^{1\mid 0} \times_{Y^\circledS}  Z^\circledS \\
(g,  \ \ a,  \ \ b)\hspace{10mm}&\mapsto \hspace{5mm}(\mu_{Y^\circledS}^\dagger (g, a), \eta_{Z^\circledS} (g, b)). \notag
\end{align}
According to Proposition \ref{p1811},  Lemma \ref{p198}, and the discussion preceding Lemma \ref{p1844}, this action may be reconstructed category-theoretically from $(\mfS \mfc \mfh_{/X^\circledS}^\circledS, Y^\circledS, (Z^\circledS, \sigma_{Z^\circledS/Y^\circledS}))$.

%-----------------------------------------------------------------------[begin proposition]------------------
\vspace{3mm}
\ble
 \label{p190}\leavevmode\\
 \ \ \
%Let $Y^\circledS$,  $(Z^\circledS, \sigma_{Z^\circledS/Y^\circledS})$ $\mu\eta_{Z^\circledS}^\dagger$, and $\lambda^\circledS$ be as in  Lemma \ref{p1888}.
Let  $\mu^\circledS : \mbA_{Z^\circledS}^{1\mid 0} \times_{Y^\circledS} \mbA_{Z^\circledS}^{1\mid 0} \migi \mbA_{Z^\circledS}^{1\mid 0}$ be a morphism  in $\mfS \mfc \mfh^\circledS_{/Y^\circledS}$ and consider the following condition concerning $\mu^\circledS$:
%\vspace{2mm}
\begin{itemize}
\item[$(Q)_{\mu^\circledS}$:]
There exists a fermionic twist  $W^\circledS$ of $Y^\circledS$
satisfying  that the $\mbA^{0|1}$-twist $(\mbA^{0|1}_{Y^\circledS \leadsto W^\circledS}, \sigma_{\mbA^{0|1}_{Y^\circledS \leadsto W^\circledS}})$ (cf. (\ref{e57})) over  $Y^\circledS$ associated with $W^\circledS$ coincides with
$(Z^\circledS, \sigma_{Z^\circledS/Y^\circledS})$ and the equality $\mu^\circledS = \mu_{Y^\circledS \leadsto W^\circledS}$ holds.
\end{itemize}
\vspace{2mm}
Then,  %$\mu^\circledS$ satisfies 
the above condition $(Q)_{\mu^\circledS}$ is equivalent that 
%if and only if 
% and an isomorphism $h : (\mbA^{0 \mid 1}_{Y \leadsto}, \sigma_{\mbA^{0 \mid 1}_{Y \leadsto}}) \isom (W/Y, \sigma_{W/Y})$ of $\mbA^{0 \mid 1}$-twists over $Y$ such that $\mu$ coincides with $\mu_{Y \leadsto W}$ via the isomorphism $\mr{id}_{\mbA^{1\mid 0}_Y} \times h : \mbA^{1\mid 1}_{Y \leadsto W}\isom \mbA^{1\mid 0}_Z$ if 
 $\mu^\circledS$ satisfies 
 the following four conditions $(R)_{\mu^\circledS}$-$(U)_{\mu^\circledS}$:
\vspace{2mm}
\begin{itemize}
\item[$(R)_{\mu^\circledS}$:]
The square diagram
\begin{equation}
\begin{CD}
(\mbG_{m})_{Y^\circledS}  \times_{Y^\circledS}\mbA^{1\mid 0}_{Z^\circledS}   \times_{Y^\circledS}  \mbA^{1\mid 0}_{Z^\circledS}  @> \mr{id}_{(\mbG_{m})_{Y^\circledS}} \times \mu^\circledS >> (\mbG_{m})_{Y^\circledS}  \times_{Y^\circledS} \mbA^{1\mid 0}_{Z^\circledS} 
\\
@V (\mu\eta_{Z^\circledS}^\dagger \times \mu\eta_{Z^\circledS}^\dagger) \circ \lambda^\circledS VV @VV \mu\eta_{Z^\circledS}^\ddagger V
\\
\mbA^{1\mid 0}_{Z^\circledS}   \times_{Y^\circledS} \mbA^{1\mid 0}_{Z^\circledS}    @> \mu^\circledS >>  \mbA^{1\mid 0}_{Z^\circledS} 
\end{CD}
\end{equation}
is commutative, where $\mu\eta_{Z^\circledS}^\ddagger$ denotes the action of $(\mbG_m)_{Y^\circledS}$ on $\mbA_{Z^\circledS}^{1\mid 0}$ given by $(g, a) \mapsto \mu\eta_{Z^\circledS}^\dagger (g^2, a)$ and  $\lambda^\circledS$ denotes the morphism
\begin{align} \label{DD017}
\lambda^\circledS : (\mbG_{m})_Y  \times_Y  \mbA^{1\mid 0}_{Z}  \times_Y \mbA^{1\mid 0}_{Z}   &  \migi  
(\mbG_{m})_Y  \times_Y  \mbA^{1\mid 0}_{Z}  \times_Y (\mbG_{m})_Y  \times_Y \mbA^{1\mid 0}_{Z} \\
%(\mbG_{m})_Y^\dagger \times_Y \mbA_Z^{1\mid 0} \times_Y   (\mbG_{m})_Y^\dagger \times_Y \mbA_Z^{1\mid 0} \\
(g,  \ \ a_1, \ \  a_2) \hspace{8mm} & \mapsto \hspace{15mm} (g, \ \  a_1, \ \  g, \ \  a_2); \notag 
\end{align}

\vspace{2mm}
\item[$(S)_{\mu^\circledS}$:]
The square diagrams
\begin{equation}
\begin{CD}
(\mbG_{m})_{Y^\circledS}  \times_{Y^\circledS}\mbA^{1\mid 0}_{Z^\circledS}   \times_{Y^\circledS}  \mbA^{1\mid 0}_{Z^\circledS}  @> \mr{id}_{(\mbG_{m})_{Y^\circledS}} \times \mu^\circledS >> (\mbG_{m})_{Y^\circledS}  \times_{Y^\circledS} \mbA^{1\mid 0}_{Z^\circledS} 
\\
@V \mu\eta_{Z^\circledS}^\dagger \times \mr{id}_{\mbA^{1\mid 0}_{Z^\circledS}}  VV @VV \mu\eta_{Z^\circledS}^\dagger  V
\\
\mbA^{1\mid 0}_{Z^\circledS}   \times_{Y^\circledS} \mbA^{1\mid 0}_{Z^\circledS}    @> \mu^\circledS >>  \mbA^{1\mid 0}_{Z^\circledS} 
\end{CD}
\end{equation}
and 
\begin{equation}
\begin{CD}
(\mbG_{m})_{Y^\circledS}  \times_{Y^\circledS}\mbA^{1\mid 0}_{Z^\circledS}   \times_{Y^\circledS}  \mbA^{1\mid 0}_{Z^\circledS}  @> \mr{id}_{(\mbG_{m})_{Y^\circledS}} \times \mu^\circledS >> (\mbG_{m})_{Y^\circledS}  \times_{Y^\circledS} \mbA^{1\mid 0}_{Z^\circledS} 
\\
@V (\mu\eta_{Z^\circledS}^\dagger \times \mr{id}_{\mbA^{1\mid 0}_{Z^\circledS}})  \circ \theta^\circledS VV @VV \mu\eta_{Z^\circledS}^\dagger   V
\\
\mbA^{1\mid 0}_{Z^\circledS}   \times_{Y^\circledS} \mbA^{1\mid 0}_{Z^\circledS}    @> \mu^\circledS >>  \mbA^{1\mid 0}_{Z^\circledS} 
\end{CD}
\end{equation}
are commutative, where $\theta^\circledS$ denotes the isomorphism
\begin{align}
(\mbG_{m})_{Y^\circledS}  \times_{Y^\circledS}\mbA^{1\mid 0}_{Z^\circledS}   \times_{Y^\circledS}  \mbA^{1\mid 0}_{Z^\circledS} & \migi (\mbG_{m})_{Y^\circledS}  \times_{Y^\circledS}\mbA^{1\mid 0}_{Z^\circledS}   \times_{Y^\circledS}  \mbA^{1\mid 0}_{Z^\circledS} \\
(g, \ \ a_1, \ \  a_2) \hspace{10mm}& \mapsto \hspace{15mm} (g,  \ \ a_2,  \ \ a_1); \notag
\end{align}

\vspace{2mm}
\item[$(T)_{\mu^\circledS}$:]
Let us write
\begin{equation}
p^\circledS := \sigma^{[1]}_{Y^\circledS} \times \sigma_{Z^\circledS/Y^\circledS}
%\tau_{\mbA_{Y \leadsto Z}^{1\mid 1}} 
 : Y^\circledS \migi \mbA^{1 \mid 0}_{Z^\circledS}, \ \ \ q^\circledS := \sigma^{[0]}_{Y^\circledS} \times \mr{id}_{Z^\circledS} : Z^\circledS \migi  \mbA^{1\mid 0}_{Z^\circledS}.
\end{equation}
Then, the following equalities hold:
\begin{align} 
\mu^\circledS \circ (p^\circledS \times p^\circledS) = p^\circledS &: Y^\circledS \migi \mbA_{Z^\circledS}^{1\mid 0}; \label{e01} \\
%\end{align}
%\begin{align}
\mu^\circledS \circ (p^\circledS \times q^\circledS) = q^\circledS &:   Z^\circledS \migi  \mbA^{1\mid 0}_{Z^\circledS};  \label{e02} \\
\mu^\circledS \circ (q^\circledS \times p^\circledS) = q^\circledS &:   Z^\circledS \migi  \mbA^{1\mid 0}_{Z^\circledS}. \label{e03}
\end{align}
Also, it holds the equality 
\begin{align}
%(\pi_{Y} \times \mr{id}_{Z}) \circ
 \mu^\circledS \circ (q^\circledS \times q^\circledS) 
 = \sigma^{[0]}_{Z^\circledS}\circ \sigma_{Z^\circledS/Y^\circledS} \circ  (h^\circledS \times h^\circledS) \label{e04} 
%\end{align}
%\begin{align}
%\pi \circ (\sigma \times \sigma \circ \sigma^{[1]})
\end{align}
of morphisms $Z^\circledS \times_{Y^\circledS} Z^\circledS \migi \mbA^{1|0}_{Z^\circledS}$, where $h^\circledS$ denotes the structure morphism $Z^\circledS  \migi Y^\circledS$ of $Z^\circledS$;

\vspace{2mm}
\item[$(U)_{\mu^\circledS}$:]
The morphism 
\begin{equation} \label{e05}
((\mr{id}_{\mbA_{Y^\circledS}^{1\mid 0}} \times h^\circledS) \circ \mu^\circledS \circ (q^\circledS \times q^\circledS))_b :   (Z^\circledS \times_{Y^\circledS} {Z^\circledS})_b  \migi (\mbA_{Y^\circledS}^{1\mid 0})_b
\end{equation}
is a closed immersion of  schemes.
%There exists an object in $\mr{SSch}_Y$ of the form $\mr{Spec}(K)$ for some field (which may be reconstructed category-theoretically by Propostion ?) and a morphism  $s : \mbA^{0 \mid 1}_k \migi Y_b$ such that the composite $\pi (\mu) \circ  s$ does not fartor through $ \mbA^{0 \mid 1}_k \migi  \mbA^{0 \mid 0}_k$ ($= \mr{Spec}(k)$).
\end{itemize}
%Then, $Z$ is, Zariski  locally,  isomorphic to $\mbA_Y^{0 \mid 1}$.
Moreover, if these equivalent conditions are satisfied, then such a fermionic twist $W^\circledS$ in $(Q)_{\mu^\circledS}$ is uniquely determined up to isomorphism.

Consequently,  the objects $\mbA^{0|1}_{Y^\circledS \leadsto W^\circledS}$ (where $W^\circledS$ is any fermionic twist of $Y^\circledS$) together with morphisms $\sigma_{\mbA^{0|1}_{Y^\circledS \leadsto W^\circledS}}$  and $\mu_{Y^\circledS \leadsto W^\circledS}$ may be reconstructed (up to isomorphism) category-theoretically from the data $(\mfS \mfc \mfh_{/X^\circledS}^\circledS, Y^\circledS)$.
%$(\mbA^{0|1}_{Y^\circledS \leadsto W^\circledS}, \sigma_{\mbA^{0|1}_{Y^\circledS \leadsto W^\circledS}})$
 \ele
%-----------------------------------------------------------------------[begin proof]-------------------
\begin{proof}
Let
$\mu^\circledS$ be a morphism satisfying the required four conditions.
It corresponds, Zariski locally on $Y_b$,  to a homomorphism
\begin{equation}
\mu^\flat : \mcO_{Y^\circledS} [t, \psi] \migi \mcO_{Y^\circledS} [t, \psi]\otimes_{\mcO_{Y^\circledS}} \mcO_{Y^\circledS} [t, \psi]
\end{equation}
of  $\mcO_{Y^\circledS}$-superalgebras.
By the conditions $(T)_{\mu^\circledS}$ and  $(U)_{\mu^\circledS}$,  $\mu^\flat$ may be  given by
\begin{align}
t \mapsto a_1 \cdot t \otimes t + a_2 \cdot \psi \otimes \psi + b_1 \cdot \psi \otimes t + b_2 \cdot t \otimes \psi
\end{align}
and
\begin{align}
\psi  \mapsto b_3 \cdot t \otimes t + b_4 \cdot  \psi \otimes \psi + a_3 \cdot \psi \otimes t + a_4  \cdot t \otimes \psi,
\end{align}
where $a_i \in \Gamma (Y_b, \mcO_{Y_b})$ and $b_i \in \Gamma (Y_b, \mcO_{Y_f})$ ($1 \leq i \leq 4$).
The equality (\ref{e01}) implies that $a_1 = 1$ and $b_3 =0$.
The equality (\ref{e02}) implies that $b_2 =0$ and $a_4 =1$.
The equality (\ref{e03}) implies that $b_1 =0$ and $a_3 =1$.
The equality (\ref{e04}) implies that $b_4 =0$.
%The assertion follows from the definition of $\tau_Y$.
Hence, the morphism (\ref{e05}) corresponds to the homomorphism $\mcO_{Y_b} [t] \migi \mcO_{Y_b}\oplus (\mcO_{Y_b} \cdot \psi \otimes \psi)$ of $\mcO_{Y_b}$-algebras  given by $t \mapsto a_2 \cdot \psi \otimes \psi$.
But, the condition $(W)_{\mu^\circledS}$ implies that $a_2 \in \Gamma (Y_b, \mcO_{Y_b}^\times)$.
Thus, there exists a Zariski open covering $\{ U_\alpha \}_{\alpha \in I}$ of $Y_b$ such that
the pair $(\mbA_{Z^\circledS}^{1|0}, \mu^\circledS)$ may be obtained by gluing the pairs $(\mbA_{Y^\circledS |_{U_\alpha}}^{1|1}, \mu^\circledS_\alpha)$ together, where
$\mu^\circledS_\alpha$ denotes the morphism $\mbA_{Y^\circledS |_{U_\alpha}}^{1|1} \times_{Y^\circledS |_{U_\alpha}} \mbA_{Y^\circledS |_{U_\alpha}}^{1|1} \migi  \mbA_{Y^\circledS |_{U_\alpha}}^{1|1}$ corresponding to the homomorphism
 \begin{align}
 \mcO_{Y^\circledS |_{U_\alpha}}[t, \psi]  & \migi \mcO_{Y^\circledS |_{U_\alpha}}[t, \psi] \otimes_{\mcO_{Y^\circledS |_{U_\alpha}}} \mcO_{Y^\circledS |_{U_\alpha}}[t, \psi] \\
t \hspace{10mm}& \mapsto  \hspace{10mm} t \otimes t + s_\alpha \cdot \psi \otimes \psi \notag \\
\psi \hspace{10mm} & \mapsto  \hspace{10mm} t \otimes \psi + \psi \otimes t. \notag
 \end{align}
%  given by assigning
%$t \mapsto t \otimes t + s_\alpha \cdot \psi \otimes \psi$ and $\psi \mapsto t \otimes \psi + \psi \otimes t$ 
(for some $s_\alpha \in \Gamma (U_\alpha, \mcO^\times_{U_\beta})$).
If $U_{\alpha, \beta} := U_\alpha \cap U_{\beta} \neq \emptyset$, then the gluing automorphism $\xi_{\alpha, \beta}^\circledS$ of  $\mbA_{Y^\circledS |_{U_{\alpha, \beta}}}^{1|1}$ (over $\mbA_{Y^\circledS |_{U_{\alpha, \beta}}}^{1|0}$) is given by $\psi \mapsto t_{\alpha, \beta} \cdot  \psi$ for some $t_{\alpha, \beta} \in \Gamma (U_{\alpha, \beta}, \mcO^\times_{U_{\alpha, \beta}})$.
Since  $\xi_{\alpha, \beta}^\circledS$ is compatible with $\mu^\circledS_\alpha$ and $\mu^\circledS_\beta$, we have the equality $s_\alpha = t_{\alpha, \beta}^2 \cdot s_\beta$.
Hence, we obtain a collection of data $(\{ U_\alpha \}_\alpha, \{ s_\alpha \}_\alpha, \{ t_{\alpha, \beta} \}_{\alpha, \beta})$ representing an element $a$ of $H^1_{\text{\'{e}t}} (Y_b, \mu_2)$.
One verifies immediately that $W^\circledS := {^a Y}^\circledS$ becomes the required fermionic twit of $Y^\circledS$. 
%Consequently,  $\mu^\circledS$ coincides (locally on $Y_b$)  with $\mu_{Y  \leadsto W_{a_2}}$ associated with the fermion-multiplicative deformation $W_{a_2}$ by means of $a_2$.
This completes the proof of Lemma \ref{p190}.
\end{proof}
%\vspace{3mm}
%-----------------------------------------------------------------------[end proposition]-------------------

%-----------------------------------------------------------------------[begin proposition]------------------
%\vspace{3mm}
\ble
%[{\bf Characterization of  $\alpha$}]
 \label{p1888}\leavevmode\\
 \ \ \
We shall assume that there exist a fermionic  twist $W^\circledS$ of $Y^\circledS$
and  an isomorphism $h^\circledS : (\mbA^{0 \mid 1}_{Y^\circledS \leadsto W^\circledS}, \sigma_{\mbA^{0 \mid 1}_{Y^\circledS \leadsto W^\circledS}}) \isom (Z^\circledS, \sigma_{Z^\circledS/Y^\circledS})$ of $\mbA^{0|1}$-twists.
(This assumption may be characterized category-theoretically thanks to Lemma \ref{p190}.)
Let  $\alpha^\circledS : \mbA^{1\mid 0}_{Z^\circledS}   \times_{Y^\circledS}  \mbA^{1\mid 0}_{Z^\circledS}  \migi   \mbA^{1\mid 0}_{Z^\circledS}$ be a morphism  in $\mfS \mfc \mfh_{/Y^\circledS}^\circledS$.
%there exist a fermionic  twist  $W^\circledS$ of $Y^\circledS$
%and an isomorphism $h : (\mbA^{0 \mid 1}_{Y \leadsto W}, \sigma_{\mbA^{0 \mid 1}_{Y \leadsto W}}) \isom (W/Y, \sigma_{W/Y})$ of $\mbA^{0 \mid 1}$-twists over $Y$
%such that 
Then,   $\alpha^\circledS$ coincides with $\alpha_{Y^\circledS \leadsto W^\circledS}$ (cf. (\ref{e10})) via the isomorphism $ h^\circledS \times \mr{id}_{\mbA^{1|0}}: \mbA^{1 \mid 1}_{Y^\circledS \leadsto W^\circledS} \isom \mbA^{1\mid 0}_{Z^\circledS}$
  if and only if 
$\alpha^\circledS$  satisfies  the following two conditions $(V)_{\alpha^\circledS}$ and $(W)_{\alpha^\circledS}$:
\vspace{2mm}
\begin{itemize}
\item[$(V)_{\alpha^\circledS}$:]
The square diagram
\begin{equation}
\begin{CD}
(\mbG_{m})_{Y^\circledS}  \times_{Y^\circledS}  \mbA^{1\mid 0}_{Z^\circledS}  \times_{Y^\circledS} \mbA^{1\mid 0}_{Z^\circledS} @> \mr{id}_{(\mbG_{m})_{Y^\circledS}} \times \alpha^\circledS >> (\mbG_{m})_{Y^\circledS}   \times_{Y^\circledS}\mbA^{1\mid 0}_{Z^\circledS}
\\
@V (\mu\eta^\dagger_{Z^\circledS} \times \mu\eta^\dagger_{Z^\circledS}) \circ \lambda^\circledS VV @VV \mu\eta^\dagger_{Z^\circledS} V
\\
 \mbA^{1\mid 0}_{Z^\circledS} \times_{Y^\circledS} \mbA^{1\mid 0}_{Z^\circledS}  @> \alpha^\circledS  >>  \mbA^{1\mid 0}_{Z^\circledS}
\end{CD}
\end{equation}
is commutative, where $\lambda^\circledS$ is as defined  in (\ref{DD017}).
%  denotes the morphism
%\begin{align}
%(\mbG_{m})_Y  \times_Y  \mbA^{1\mid 0}_{Z}  \times_Y \mbA^{1\mid 0}_{Z}   &  \migi  
%(\mbG_{m})_Y  \times_Y  \mbA^{1\mid 0}_{Z}  \times_Y (\mbG_{m})_Y  \times_Y \mbA^{1\mid 0}_{Z} \\
%(\mbG_{m})_Y^\dagger \times_Y \mbA_Z^{1\mid 0} \times_Y   (\mbG_{m})_Y^\dagger \times_Y \mbA_Z^{1\mid 0} \\
%(g, a_1, a_2) \hspace{15mm} & \mapsto \hspace{20mm} (g, a_1, g, a_2). \notag 
%\end{align}
%\vspace{1mm}
\item[$(W)_{\alpha^\circledS}$:]
We have the equalities
\begin{equation}
\alpha^\circledS \circ ((\sigma_{Z^\circledS}^{[0]} \circ\sigma_{Z^\circledS/Y^\circledS}) \times  \mbA^{1\mid 0}_{Z^\circledS})  =  \alpha^\circledS \circ ( \mbA^{1\mid 0}_{Z^\circledS} \times (\sigma_{Z^\circledS}^{[0]} \circ\sigma_{Z^\circledS/Y^\circledS})) = \mr{id}_{ \mbA^{1\mid 0}_{Z}}.
\end{equation}
of endomorphisms of $ \mbA^{1\mid 0}_{Z^\circledS}$.
\end{itemize}
Consequently,  the objects $\mbA^{0|1}_{Y^\circledS \leadsto W^\circledS}$ (where $W^\circledS$ is any fermionic twist of $Y^\circledS$) together with morphisms $\sigma_{\mbA^{0|1}_{Y^\circledS \leadsto W^\circledS}}$  and $\alpha_{Y^\circledS \leadsto W^\circledS}$ may be reconstructed (up to isomorphism) category-theoretically from the data $(\mfS \mfc \mfh_{/X^\circledS}^\circledS, Y^\circledS)$.
%$(\mbA^{0|1}_{Y^\circledS \leadsto W^\circledS}, \sigma_{\mbA^{0|1}_{Y^\circledS \leadsto W^\circledS}})$

 \ele
%-----------------------------------------------------------------------[begin proof]-------------------
\begin{proof}
The assertion follows from an argument similar to the argument in the proof of Lemma \ref{p1833}.
\end{proof}
\vspace{3mm}
%-----------------------------------------------------------------------[end proposition]-------------------

%-----------------------------------------------------------------------[begin proposition]------------------
\vspace{3mm}
\bco[{\bf Characterization of fermionic twists over $Y^\circledS$}] \label{p218}\leavevmode\\
 \ \ \ 
%Let $Y^\circledS$ be an object of $\mr{SSch}_X$.
The collection of  fermionic twists over $Y^\circledS$ (i.e., a collection of topological spaces together with a sheaf of superrings)  are reconstructed  category-theoretically (up to isomorphism) from the data $(\mfS \mfc \mfh_{/X^\circledS}^\circledS, Y^\circledS)$.
Moreover, this reconstruction is functorial (in a natural sense) in $Y^\circledS \in \mfS \mfc \mfh_{/X^\circledS}^\circledS$.
\eco
%-----------------------------------------------------------------------[begin proof]-------------------
\begin{proof}
The assertion follows from Proposition \ref{p12}, Lemma \ref{p190}, Lemma \ref{p1888}, and the discussion in \S\,\ref{sub5} (especially, the isomorphism (\ref{e09})).
%It follows from Proposition \ref{p1811}, Proposition \ref{p188}, Lemma \ref{p1888}, and  Lemma \ref{p190}
%that for each open subsuperscheme $Y |_U$ of $Y$ (which may be reconstructed category-theoretically by  Proposition \ref{p12}),
%the collection of data 
%\begin{equation}
%\big\{ (\mbA_{Y |_U \leadsto W}^{1\mid 1}, \alpha_{Y |_U \leadsto W}, \mu_{Y |_U \leadsto W}) \ \big| \ \text{$W$ is an {\it fmt} of $Y |_U$} \big\}
%\end{equation}
%    in $\mr{SSch}_Y$  may be reconstructed (up to isomorphism) category-theoretically from the data $(\mr{SSch}_X, Y, Y|_U)$.
%Also, note that  the isomorphism $\Gamma (W_b, \mcO_W) \isom \mr{Map}_{\mr{SSch}_Y} (Y, \mbA_{Y \leadsto W}^{1\mid 1})$ of  super-rings  described in  (\ref{e09})   is functroial with respect to restriction to open immersions $Y |_V \migi Y |_U$.
%Hence, by taking $Y |_U$ equal to various open subsuperschemes, we obtain (by gluing) the  collection of  {\it fmt}s over $Y$.
\end{proof}
%\vspace{3mm}
%-----------------------------------------------------------------------[end proposition]-------------------

\vspace{5mm}
%----------------------------------------------------------------------[begin subsection]-------------
\subsection{} \label{sub331}

We turn to the proof of the main result of the present paper, i.e., Theorem A.
Before beginning the proof,  
let us first mention  the following rigidity property concerning $\mfS \mfc \mfh_{/X^\circledS}^\circledS$.

%-----------------------------------------------------------------------[begin proposition]------------------
%\vspace{3mm}
\bpr
%[{\bf Characterization of fermionic twists over $Y^\circledS$}] 
\label{pg218}\leavevmode\\
 \ \ \ 
Let  
$X^\circledS$ and $X'^\circledS$  be two locally noetherian superschemes.
Let 
\begin{align}
\mfI \mfs \mfo \mfm (\mfS \mfc \mfh_{/X^\circledS}^\circledS, \mfS \mfc \mfh_{/X'^\circledS}^\circledS)
\end{align}
 denotes the category of equivalences  $\mfS \mfc \mfh_{/X^\circledS}^\circledS \isom \mfS \mfc \mfh_{/X'^\circledS}^\circledS$ and 
 \begin{align}
 \mr{Isom} (\mfS \mfc \mfh_{/X^\circledS}^\circledS, \mfS \mfc \mfh_{/X'^\circledS}^\circledS)
 \end{align}
 denotes  the set of isomorphism classes of equivalences  $\mfS \mfc \mfh_{/X^\circledS}^\circledS \isom \mfS \mfc \mfh_{/X'^\circledS}^\circledS$  (i.e., the set of isomorphism classes of objects of the category $\mfI \mfs \mfo \mfm (\mfS \mfc \mfh_{/X^\circledS}^\circledS, \mfS \mfc \mfh_{/X'^\circledS}^\circledS)$).
Also, let
\begin{align}
\mr{Isom} (X'^\circledS, X^\circledS)
\end{align}
 denotes  the set of isomorphisms of superschemes  $X'^\circledS \isom X^\circledS$.
Consider the map of sets
\begin{align} \label{R010}
\mr{Isom} (X'^\circledS, X^\circledS)  \migi \mr{Isom} (\mfS \mfc \mfh_{/X^\circledS}^\circledS, \mfS \mfc \mfh_{/X'^\circledS}^\circledS)
\end{align}
which, to any isomorphism $f^\circledS : X'^\circledS \isom X^\circledS$, assigns (the isomorphism class of) the equivalence $\mfS \mfc \mfh_{/X^\circledS}^\circledS \isom \mfS \mfc \mfh_{/X'^\circledS}^\circledS$ given by base-change via $f^\circledS$.
Then, this map  (\ref{R010}) is injective.
\epr
%-----------------------------------------------------------------------[begin proof]-------------------
\begin{proof}
%The injectivity of (\ref{R010}) 
The assertion  follows  immediately from  the functorial bijection (\ref{e5}) and the various reconstructing procedures involved.
%We shall consider the surjectivity.
%Let $\phi$ be an autoequivalence of  $\mfS \mfc \mfh_{/X^\circledS}^\circledS$.
%It follows from Corollary \ref{p12open} that $\phi$ induces an automorphism $\phi_{X_b} : X_b \isom X_b$.
%The restrictions (to various quasi-compact open subschemes $U \subseteq X_b$)  of the superring object $\mbA^{1|1}_{X^\circledS}$ over $X^\circledS$ are sent by $\phi$  to the various restrictions of a certain superring object over $X^\circledS$.
%Thus, by the functorial bijection (\ref{e05}), we have an isomorphism $\phi^{-1}(\mcO_{X^\circledS}) \isom \mcO_{X^\circledS}$, that is to say, an automorphism $f^\circledS$ of $X^\circledS$.
%By considering the various reconstructing procedures involved, 
%the image of $f^\circledS$ via the map (\ref{R010}) is isomorphic to $\phi$.
%This completes  the proof of the surjectivity, and hence, Proposition \ref{pg218}.
\end{proof}
\vspace{1mm}
%-----------------------------------------------------------------------[end proposition]-------------------

%-----------------------------------------------------------------------[end definition]-------------------
\begin{rema} \label{rg90df} \leavevmode\\
 \ \ \ 
Unlike the case of schemes proved in  ~\cite{Mzk1}, Theorem 1.7 (ii),  the map (\ref{R010}) may  not  be surjective.
 Indeed, suppose that $X^\circledS = X'^\circledS = Y$ for some scheme $Y$ and there exists  a nonzero  element $a \in H^1_{\text{\'{e}t}} (Y, \mu_2)$.
 Then,  the  assignment  $Z^\circledS \mapsto {^a Z}^\circledS$ defines an
  autoequivalence ${^a \phi} : \mfS \mfc \mfh_{/Y}^\circledS \isom \mfS \mfc \mfh_{/Y}^\circledS$.
 Since $Z^\circledS$ is, in general,  not isomorphic to ${^a Z}^\circledS$, ${^a \phi}$ is not isomorphic to the identity functor.
  But, one may verifies immediately that ${^a \phi}$ cannot arise  from the base-change via any  automorphism of $Y$.
  This implies that 
  the isomorphism class of  ${^a \phi}$ does not lie in the image of 
  the map (\ref{R010}).
   \end{rema}
%-----------------------------------------------------------------------[end remark]-------------------
\vspace{5mm}

Finally, by applying the results obtained so far, we  prove the remaining portion of  Theorem A (cf. Proposition \ref{p23u}) as follows:

\vspace{1mm}

\begin{proof}[Proof of Theorem A]
%Thanks to Propsotion ?, it suffices to prove that
%$\mr{SSch}_{X} \cong \mr{SSch}_{X'}$ implies $X \sim X'$.
Suppose that we are given an equivalence of categories:
\begin{equation}
\phi : \mfS \mfc \mfh_{/X^\circledS}^\circledS \isom \mfS \mfc \mfh_{/X'^\circledS}^\circledS.
\end{equation}
Let us take   a Zariski open  covering $\{ U_\alpha \}_{\alpha \in I}$   of $X_b$,  where each $U_\alpha$ is  quasi-compact, i.e., $X^\circledS |_{U_\alpha} \in \mr{Ob} (\mfS \mfc \mfh_{/X^\circledS}^\circledS)$.
The image   $\phi (X^\circledS |_{U_\alpha})$ of  $X^\circledS |_{U_\alpha}$ (for each $\alpha \in I$)    is isomorphic (as an object of $\mfS \mfc \mfh_{/X'^\circledS}^\circledS$) to  
$X'^\circledS |_{U'_\alpha}$ for some quasi-compact open subscheme $U'_\alpha$ of  $X_b'$ (cf. Proposition \ref{p12}). 
%a quasi-compact open subsuperscheme of $X'^\circledS$ , and hence $\phi (X^\circledS |_{U_X}) \cong Y^\circledS |_{U_Y}$ for some quasi-compact open $U_Y$ in $Y_b$.
% a quasi-compact open subscheme $U_X$  of $X_b$ (hence, $X^\circledS |_{U_X} \in \mr{Ob} ( \mfS \mfc \mfh_{/X^\circledS}^\circledS)$).
%Then,  $\phi (X^\circledS |_{U_X})$ is isomorphic to  a quasi-compact open subsuperscheme of $X'^\circledS$ (cf. Proposition \ref{p12}), and hence $\phi (X^\circledS |_{U_X}) \cong Y^\circledS |_{U_Y}$ for some quasi-compact open $U_Y$ in $Y_b$.
It follows from  Corollary \ref{p218}  (and the various reconstructing procedures involved) that 
one may find  an isomorphism  $\iota_\alpha^\circledS  : Z_\alpha^\circledS \isom X^\circledS |_{U_\alpha}$ of superschemes, where $Z_\alpha^\circledS$ denotes  a fermionic twist   of  $X'^\circledS |_{U'_\alpha}$;
 such an isomorphism  $\iota_\alpha^\circledS$ is uniquely determined (thanks to Proposition \ref{pg218}) by the condition that the  functor $ \mfS \mfc \mfh_{/X^\circledS |_{U_\alpha}}^\circledS \isom \mfS \mfc \mfh_{/Z_\alpha^\circledS}^\circledS$ given by base-change via $\iota_\alpha^\circledS$ is isomorphic to 
the composite functor 
\begin{align}
\iota^{\mfS \mfc \mfh}_\alpha : \mfS \mfc \mfh_{/X^\circledS |_{U_\alpha}}^\circledS \stackrel{ \phi |_{U_\alpha}}{\migi} \mfS \mfc \mfh_{/X'^\circledS |_{U'_\alpha}}^\circledS \stackrel{(\ref{R02})}{\migi} \mfS \mfc \mfh_{/Z_\alpha^\circledS}^\circledS,
\end{align}
 where the first arrow denotes the restriction of $\phi$ to $\mfS \mfc \mfh_{/X^\circledS |_{U_\alpha}}^\circledS$.
For any pair $(\alpha, \beta) \in I \times I$  with $U_{\alpha, \beta} := U_\alpha \cap U_\beta \neq \emptyset$, we obtain an isomorphism $\iota_{\alpha, \beta}^\circledS := (\iota^\circledS_\beta)^{-1} \circ \iota^\circledS_\alpha : Z_\alpha^\circledS |_{U_{\alpha, \beta}} \isom Z_\beta^\circledS |_{U_{\alpha, \beta}}$.
 Proposition \ref{pg218} implies that
 the collection of isomorphisms $\{ \iota^\circledS_{\alpha, \beta} \}_{\alpha, \beta}$ satisfies the cocycle condition (in an evident  sense), and hence,
 the  superschemes $\{ Z^\circledS_\alpha \}_{\alpha \in I}$ may be glued (by means of $\{ \iota^\circledS_{\alpha, \beta} \}_{\alpha, \beta}$) together to a superscheme $Z^\circledS$.
By construction, $Z^\circledS$ is a fermionic twist of $X'^\circledS$ and the isomorphisms $\{ \iota^\circledS_\alpha \}_{\alpha \in I} $ may be glued together to an isomorphism $\iota^\circledS : Z^\circledS \isom X^\circledS$.
Consequently,  we have $X^\circledS \stackrel{f}{\sim} X'$.
This completes the proof of Theorem A.
\end{proof}

%%%%%%%%%%%%%%%%%%%%%%%%%%%%%%%%--[ begin  section1]---%%%%%%
\vspace{10mm}
\section{Further rigidity properties} \vspace{3mm}

In this final section, we  propose further rigidity properties concerning the category of superschemes.
%Let us consider an automorphism of the identity of $\mr{id}_{\mr{SSch}_X}$, i.e., the 
%For each super-scheme $Y := (Y_b, \mcO_Y)$, we shall write $(-1)_Y$ for the automorphism
%\begin{equation}
%(-1)_Y := (\mr{id}_{Y_b}, (-1)^\sharp_Y) : Y \isom Y
%\end{equation}
%of $Y$, where $(-1)^\sharp_Y$ denotes the automorphism of $\mcO_Y = \mcO_{Y_b}\oplus \mcO_{Y_f}$ given by
%$(a, b) \mapsto (a, -b)$.
%Note that $(-1)_Y \circ (-1)_Y = \mr{id}_Y$ and for each morphism $f : Y \migi X$, we have the equality $f \circ (-1)_Y = (-1)_{X} \circ f$.
%If $Y$ is bosonic, then $(-1)_Y = \mr{id}_Y$.
%Since  the  fiber products in $\mr{SSch}_X$ exist (cf. ?, Corollary 10.3.9),
%one may define, for each morphism $f : Y \migi X$ of super-schemes,
%a functor
%\begin{equation}
%\mr{SSch}_f : \mr{SSch}_X \migi \mr{SSch}_Y
%\end{equation}
%induced by base-change by $f$. 

%-----------------------------------------------------------------------[begin proposition]------------------
\vspace{3mm}
\bpr \label{p756}\leavevmode\\
 \ \ \ 
Let
$X^\circledS$ and $Y^\circledS$ be two locally noetherian superschemes.
% both of whose underlying schemes are locally noetherian.
Also, let $f^\circledS \ (:= (f_b, f^\flat)) : Y^\circledS \migi X^\circledS$  be  a morphism of superschemes such that $f_b$ is quasi-compact.
We shall write 
\begin{equation}
\mfS \mfc \mfh^\circledS_{f^\circledS} : \mfS \mfc \mfh_{/X^\circledS}^\circledS \migi \mfS \mfc \mfh_{/Y^\circledS}^\circledS
\end{equation} 
for the functor induced by base-change via $f^\circledS$.
Then, the following properties are satisfied.

%Let us write $(W)_f$, $(X)_f$ for the condition described as follows:
%\begin{itemize}
%\item[$(W)_f$]
%$f^\sharp : f^*_b (\mcO_X) \migi \mcO_Y$ is surjective.
%$f_{b*} (\mcO_{Y_b})$ is nontrivial on every connected component of  $Y$;
%\item[$(X)_f$]
%$Y$ is bosonic.
%\end{itemize}
%Then, the following assertions are hold.
\begin{itemize}
\item[(i)]
If  
 %$f^\sharp : f^*_b (\mcO_X) \migi \mcO_Y$ is surjective and 
$f_{b*} (\mcO_{Y_b})\neq 0$,
% on every  connected component of  $X_b$ is not trivial.
 %be a bosonic morphism of super-schemes whose underlying morphism of schemes is quasi-compact.
%quasi-compact  and bosonic morphism of super-schemes.
then the functor $\mfS \mfc \mfh^\circledS_{f^\circledS}$
 has no nontrivial automorphisms.
\item[(ii)]
If   $X^\circledS$ is a scheme (i.e., $\mcO_{X_f} =0$), 
then a nontrivial automorphism of $\mfS \mfc \mfh^\circledS_{f^\circledS}$ is uniquely determined as the automorphism   given by 
the collection of automorphisms $\{ (-1)_{Z^\circledS \times_{X^\circledS} Y^\circledS} \}_{Z^\circledS \in \mr{Ob}(\mfS \mfc \mfh_{/X^\circledS}^\circledS)}$.
% (cf. Remark \ref{r4gg0d}).
% of  the superschemes $Y \times_X Z$
 % ( $Z^\circledS \in \mr{Ob}(\mfS \mfc \mfh_{/X^\circledS}^\circledS)$).
\end{itemize}

 \epr
%-----------------------------------------------------------------------[begin proof]-------------------
\begin{proof}
%First, we verify assertion (i).
%Before beginning the proof of Proposition \ref{p756},
First,  let us  make the following observation.
Let $\zeta$ be an automorphism of $\mfS \mfc \mfh^\circledS_{f^\circledS}$, which consists of  autormorphisms
\begin{equation}
\zeta^\circledS_{Z^\circledS} \ (:= (\zeta_{Z, b}, \zeta_Z^\flat)) : Y^\circledS \times_{X^\circledS} Z^\circledS \isom Y^\circledS \times_{X^\circledS} Z^\circledS
\end{equation}
in $\mfS \mfc \mfh_{/Y^\circledS}^\circledS$  (for   $Z^\circledS \in \mr{Ob} (\mfS \mfc \mfh_{/X^\circledS}^\circledS)$) that are functorial in $Z^\circledS$.
If $\mfS \mfc \mfh_{f_b} : \mfS \mfc \mfh_{/X_b} \migi \mfS \mfc \mfh_{/Y_b}$ denotes the functor defined by base-change via $f_b : Y_b \migi X_b$, then it makes 
 the following  square diagram
\begin{align}
\xymatrix{
\mfS \mfc \mfh_{/X_b}\ar[r]^{\mfS \mfc \mfh_{f_b}}  \ar[d]_{\mfS \mfc \mfh_{\beta_{X^\circledS}}} & \mfS \mfc \mfh_{/Y_b} \ar[d]^{\mfS \mfc \mfh_{\beta_{Y^\circledS}}}  
\\
\mfS \mfc \mfh^\circledS_{/X^\circledS} \ar[r]_{\mfS \mfc \mfh^\circledS_{f^\circledS}} & \mfS \mfc \mfh^\circledS_{/Y^\circledS}
}
\end{align}
commute, where the left-hand and right-hand vertical arrows arise from base-change  via $\beta_{X^\circledS}$ and $\beta_{Y^\circledS}$ respectively.
Since  $(W \times_{X_b} Y^\circledS)_b = W \times_{X_b} Y_b$ (for any $W \in \mr{Ob} (\mfS \mfc \mfh_{/X_b})$), the autormorphism $\zeta$ restricts to an automorphism $\zeta |_{\mfS \mfc \mfh_{/X_b}}$ of  
$\mfS \mfc \mfh_{f_b}$, which is given by $\{ \zeta_{W \times_{X_b} X^\circledS, b} \}_{W \in \mr{Ob} (\mfS \mfc \mfh_{/X_b})}$.
By ~\cite{Mzk1}, Theorem 1.7, (i), we have $\zeta_{W \times_{X_b} X^\circledS, b} = \mr{id}_{W \times_{X_b} Y_b}$ for any $W \in \mr{Ob} (\mfS \mfc \mfh_{/X_b})$.
%$\zeta_{Z, b} = \mr{id}_{Y \times_X Z}$.
%any functorial automorphsm of the objects of the essential image $\mr{Im} (\mr{SSch}_f)$ of $\mr{SSch}_{X'}$ is the identity on the underlying schemes.
%Let us consider the case where $Z = \mbA_X^{0 \mid 1}$ and $Z = \mbA_X^{0 \mid 1}$.
In particular, the equality $\zeta_{\mbA^{1|0}_{X^\circledS},b} = \mr{id}_{\mbA^{1|0}_{Y_b}}$ implies the equality 
\begin{align} \label{DD019}
\zeta_{\mbA^{1|0}_{X^\circledS}}^\circledS = \mr{id}_{\mbA^{1|0}_{Y^\circledS}}.
\end{align}
Next, let us denote by $\gamma^\circledS_1$ (resp., $\gamma^\circledS_2$)  the morphism $\mbA_{Y^\circledS}^{0 \mid 2} \migi \mbA_{Y^\circledS}^{0 \mid 1}$  in $\mfS \mfc \mfh^\circledS_{/Y^\circledS}$ corresponding to the homomorphism $\mcO_{Y^\circledS} [\psi] \migi \mcO_{Y^\circledS} [\psi]\otimes_{\mcO_{Y^\circledS}} \mcO_{Y^\circledS} [\psi]$ given by $\psi \mapsto \psi \otimes 1$ (resp., $\psi \mapsto 1 \otimes \psi$).
Note that  the automorphism $\zeta^\circledS_{\mbA_{X^\circledS}^{0 \mid 1}}$ of $\mbA^{0|1}_{Y^\circledS}$ is given by $\psi \mapsto g \cdot \psi$ for some $g \in \Gamma (Y_b, \mcO_{Y_b}^\times)$.
Since
$\gamma^\circledS_\Box : \mbA_{Y^\circledS}^{0 \mid 2} \migi  \mbA_{Y^\circledS}^{0 \mid 1}$ (for each $\Box = 1, 2$) is compatible with $\zeta^\circledS_{\mbA_{X^\circledS}^{0\mid 2}}$ and $\zeta^\circledS_{\mbA_{X^\circledS}^{0\mid 1}}$ 
%, for each $\Box = 1, 2$, the square diagram
%\begin{align}
%\xymatrix{
%\mbA_{Y^\circledS}^{0 \mid 2} \ar[r]^{\mr{id}_Y \times \gamma_\Box} \ar[d]_{\zeta_{\mbA_X^{0\mid 2}} } & \mbA_{Y^\circledS}^{0 \mid 1} \ar[d]^{\zeta_{\mbA_X^{0\mid 1}} }
%\\
%\mbA_{Y^\circledS}^{0 \mid 2} \ar[r]_{\mr{id}_Y \times \gamma_\Box} & \mbA_{Y^\circledS}^{0 \mid 1}
%}
%\end{align}
%is commutative
%Thus, 
%the equalities
%\begin{align}
%\zeta_{\mbA_X^{0\mid 2}} \circ (\mr{id}_Y \times \gamma_1) &= (\mr{id}_Y \times \gamma_1) \circ \zeta_{\mbA_X^{0\mid 1}},  \\
%\zeta_{\mbA_X^{0\mid 2}} \circ (\mr{id}_Y \times \gamma_2) &= (\mr{id}_Y \times \gamma_2) \circ \zeta_{\mbA_X^{0\mid 1}} \notag
%\end{align}
 (due to the functoriality of $\zeta^\circledS_{Z^\circledS}$ with respect to  $Z^\circledS$), $\zeta_{\mbA_X^{0\mid 2}}$ is given by
 $\psi \otimes 1 \mapsto g \cdot \psi \otimes 1$ and $1 \otimes \psi \mapsto g \cdot 1 \otimes \psi$ (hence $\psi \otimes \psi \mapsto g^2 \cdot \psi \otimes \psi$).
 %$\{\zeta_Z \}_{Z \in \mr{Ob}(\mr{SSch}_X)}$)
Here, for any superscheme $Z^\circledS$, we shall write $\mbA^{\epsilon |0}_{Z^\circledS} := Z^\circledS \times \mr{Spec} (\mbZ [\frac{1}{2}][\epsilon]/\epsilon^2)$.
Since $\mbA^{\epsilon |0}_{Y^\circledS}$ lies in the essential image of the composite $\mfS \mfc \mfh_{\beta_{Y^\circledS}}\circ \mfS \mfc \mfh_{f_b}$, we have $\zeta_{\mbA^{\epsilon |0}_{X^\circledS}, b} = \mr{id}_{(\mbA^{\epsilon |0}_{Y^\circledS})_b}$
But, 
a morphism $\gamma_\epsilon^\circledS : \mbA^{0|2}_{Y^\circledS} \migi \mbA^{\epsilon |0}_{Y^\circledS}$ over $Y^\circledS$ given by 
assigning  $\epsilon \mapsto \psi \otimes \psi$ is compatible with $\zeta^\circledS_{\mbA_{X^\circledS}^{0\mid 2}}$  and $\zeta_{\mbA^{\epsilon |0}_{X^\circledS}}^\circledS$.
This implies that $g^2 =1$, equivalently, $g =1$ or $-1$.
Since we have obtained   the equality (\ref{DD019}), $\zeta^\circledS_{\mbA_{X^\circledS}^{1\mid 1}} $ coincides with either  $\mr{id}_{\mbA^{1|1}_{Y^\circledS}}$ or  $\mr{id}_{\mbA^{1|0}_{Y^\circledS}} \times (-1)_{\mbA^{0|1}}$.
%\begin{align} \label{DD034}
%\zeta^\circledS_{\mbA_{X^\circledS}^{1\mid 1}} = \mr{id}_{\mbA^{1|1}_{Y^\circledS}} \ \  \text{or} \ \  \mr{id}_{\mbA^{1|0}_{Y^\circledS}} \times (-1)_{\mbA^{0|1}}.
%\end{align}
Hence, by  the discussion in \S\,\ref{S15} (especially, the composite  bijection  \ref{e5}) and the functoriality of $Z^\circledS \mapsto \zeta_{Z^\circledS}^\circledS$,
$\zeta$ coincides with either  the identity morphism or the automorphism  given by $\{ (-1)_{Z^\circledS \times_{X^\circledS} Y^\circledS} \}_{Z^\circledS \in \mr{Ob}(\mfS \mfc \mfh_{/X^\circledS}^\circledS)}$.
%, which  we shall denote by  $(-1)_{(-) \times_{X^\circledS} Y^\circledS}$.
%But, since we have seen that $\zeta_{\mbA_X^{0\mid 2}, b} = \mr{id}_{\mbA_Y^{0\mid 2}}$,
%  $g^2$ must be equal to $1$ (i.e., $g =1$ or $g =-1$).
%For an arbitrary  $Z \in \mr{Ob} (\mr{SSch}_X)$, 
%the  isomorphism 
%\begin{equation}
%\kappa : \Gamma ((Y \times_X Z)_b, \mcO_{(Y \times_X Z)_f}) \isom \mr{Map}_{\mr{SSch}_Y} (Y \times_X Z, \mbA_Y^{0\mid 1})
%\end{equation}
% induced from  (\ref{e5})  fits into the commutative diagram
%\begin{equation} \label{c1}
%\begin{CD}
%\Gamma ((Y \times_X Z)_b, \mcO_{(Y \times_X Z)_f}) @> \kappa >> \mr{Map}_{\mr{SSch}_Y} (Y \times_X Z, \mbA_Y^{0\mid 1})
%\\
%@V \zeta_{Z}^\sharp VV @VV \zeta_{\mbA_Y^{0 \mid 1}} \circ (-) V
%\\
%\Gamma ((Y \times_X Z)_b, \mcO_{(Y \times_X Z)_f}) @> \kappa >> \mr{Map}_{\mr{SSch}_Y} (Y \times_X Z, \mbA_Y^{0\mid 1}).
%\end{CD}
%\end{equation}
%This diagram is functorial in $Z$ in a natural sense.
%Hence, if $g =1$, then $\zeta$ is the identity of $\mr{SSch}_f$,  and if  $g =-1$ (i.e., $\zeta_{\mbA_X^{0\mid 1}} = (-1)_{\mbA_Y^{0\mid 1}}$), then $\zeta$ is given by 
% $ (-1)_{Z \times_X Y}$  ($Z \in \mr{Ob}(\mr{Sch}_X)$).

Now, we shall prove assertion (i) and (ii).
Since assertion (ii) follows directly from the above discussion, it suffices to consider only assertion (i).
Since $f_{b*} (\mcO_{Y_{f}}) \neq 0$, there exists an open subscheme  $U$ of $X_b$ such that
$\Gamma (f_b^{-1}(U), \mcO_{Y_f}) \neq 0$.
But, $\zeta_{X^\circledS |_U}^\circledS$ must be the identity morphism of  $Y^\circledS |_{f_b^{-1}(U)}$.  (in particular, the fermionic part of $\zeta_{X^\circledS |_U}^\flat$ coincides with the identity morphism of $\mcO_{Y_f} |_{f_b^{-1}(U)}$.
Hence, $\zeta$ must be equal to the identity morphism.
%the commutative diagram (\ref{c1}) implies that $g =1$.
%Thus, assertion (i) follows from the above discussion, and w
This completes the proof of Proposition \ref{p756}.
\end{proof}
\vspace{3mm}
%-----------------------------------------------------------------------[end proposition]-------------------

%-----------------------------------------------------------------------[begin proposition]------------------
%\vspace{3mm}
\bpr
%[{\bf Further rigidity property}] 
\label{p4FRP}\leavevmode\\
 \ \ \ 
%\begin{itemize}
%\item[(i)]
%\item[(ii)]
Let
$X^\circledS$ be a locally noetherian superscheme.
Suppose that for any $Y^\circledS \in \mr{Ob} (\mfS \mfc \mfh^\circledS_{/X^\circledS})$, one has  an automorphism
$\zeta_{Y^\circledS}^\circledS$ of $Y^\circledS$ (which is not necessarily over $X^\circledS$) 
and
%satisfying the following property:
%with the property that
for any morphism $f^\circledS : Y^\circledS_1 \migi Y^\circledS_2$ in  $\mfS \mfc \mfh^\circledS_{/X^\circledS}$, one has a commutative square diagram:
\begin{equation}
\xymatrix{
Y^\circledS_1 \ar[r]^{\zeta^\circledS_{Y^\circledS_1}} \ar[d]_{f^\circledS} &  Y^\circledS_1 \ar[d]^{f^\circledS} \\
Y^\circledS_2  \ar[r]_{\zeta^\circledS_{Y^\circledS_2}} & Y^\circledS_2.
}
%\begin{CD}
%Y^\circledS_1 @> \zeta^\circledS_{Y^\circledS_1} >> Y^\circledS_1
%\\
%@V f^\circledS VV @VV f^\circledS V
%\\
%Y^\circledS_2 @> \zeta^\circledS_{Y^\circledS_2} >> Y^\circledS_2.
%\end{CD}
\end{equation}
Then, all of $\zeta^\circledS_{Y^\circledS}$  are either the identity morphisms or $(-1)_{Y^\circledS}$.
%  \end{itemize}
  \epr
%-----------------------------------------------------------------------[begin proof]-------------------
\begin{proof}
The assertion follows immediately  from an argument similar to the argument  in the proof of Proposition \ref{p756}.
%Let $\zeta_Y := (\zeta_{Y, b}, \zeta^\sharp_Y)$
% ($Y \in \mr{Ob}(\mr{SSch}_X)$) are as in the statemant.
%It follows from  ~\cite{Mzk1}, Theorem 1.8, that  $\zeta_{Y, b} = \mr{id}_Y$ for every $Y$ (in particular,  $\zeta_{Y, b}$ is over  $X_b$).
%Hence, one may apply
%the assertion follows from
% the discussion in the proof of Proposition \ref{p756} to complete the proof the assertion.
\end{proof}
\vspace{10mm}
%-----------------------------------------------------------------------[end proposition]-------------------

%%%%%%%%%%%%%%%%%%%%%%%%%%%%%%%%%%%%%%%%%%%%%%%%%%%
\end{document}